%% file: main.tex
\documentclass[a4paper,10pt,reqno]{amsart}

\input{Settings/packages}
\input{Settings/commands}
\input{Settings/authors}

\title[Conditioned stochastic stability of equilibrium states]{Conditioned stochastic stability of equilibrium states on uniformly expanding repellers} 
\date{\today}

\begin{document}

\begin{abstract}
   We propose a notion of conditioned stochastic stability of invariant measures on repellers: we consider whether quasi-ergodic measures of absorbing Markov processes, generated by random perturbations of the deterministic dynamics and conditioned upon survival in a neighbourhood of a repeller, converge to an invariant measure in the zero-noise limit.
   Under suitable choices of the random perturbation, we find that equilibrium states on uniformly expanding repellers are conditioned stochastically stable.
   In the process, we contribute to the rigorous foundation for the existence of ``natural measures'', which were proposed by Kantz and Grassberger in 1984 to aid the understanding of chaotic transients.
\end{abstract}

\keywords{Stochastic stability; absorbing Markov processes; conditioned random dynamics; thermodynamic formalism; transient dynamics}
\subjclass[2020]{60J05; 37D20; 37D35; 37D45}

\maketitle
{\thispagestyle{empty}\tableofcontents\newpage}

\section{Introduction}
Understanding how typical trajectories evolve in a dynamical system  and describing its relevant statistics
is a central topic in Dynamical Systems theory. 
This question is commonly addressed from an ergodic theoretical point of view, stating that each ergodic invariant measure $\mu$ provides the distribution of the trajectory starting at a point $x$, $\mu$-almost surely.
Dynamical systems often admit infinitely many ergodic invariant measures, so it is natural to ask which ones are the most meaningful or relevant to study. To tackle this, Kolmogorov and Sinai, proposed the notion of \emph{stochastic stability} of invariant measures~\cite{Kifer1974, Alves2007}.

Stochastic stability concerns the stationary measures of Markov processes generated by small bounded random perturbations of a deterministic dynamical system and their limit as the amplitude of the perturbation vanishes~\cite{Kifer1974, Alves2007}. When a stationary measure converges to an invariant measure of the original deterministic system we say that the limiting measure is stochastically stable. These measures have been recognised to highlight the statistics of (Lebesgue) typical trajectories~\cite{Young2002}. Note that stochastically stable invariant measures sit on attractors.
 
In transient dynamics~\cite{Lai2011}, trajectories that remain for a long time near a repeller have been observed to have well-defined statistics. 
While there is also an abundance of invariant ergodic measures on repellers, so-called \emph{natural measures} have been heuristically identified as the relevant invariant measures that represent observed long time behaviour of trajectories near a repeller, and provide important information regarding the statistics of transient dynamics~\cite{Kantz1985}. Despite the fact that such measures feature at the heart of the intuitive understanding of transient dynamics, their existence and mathematical properties remain to be rigorously established. 

Like stochastic stability successfully provides relevant measures on attractors, we seek a strategy to establish persistence of measures on repellers under random perturbations. 
The strategy of Kolmogorov and Sinai fails since stationary measures of the Markov process generated by random perturbations of the original system do not converge to invariant measures supported on repellers in the deterministic limit. 

In this paper, we propose a novel notion of stochastic stability for repellers referring to \emph{quasi-ergodic measures} rather than stationary measures.  
Quasi-ergodic measures originate from the theory of absorbing Markov processes and capture the typical average behaviour of trajectories conditioned upon remaining in a certain region of the state space for asymptotically long times. By conditioning the Markov process generated by random bounded perturbations of the original map upon survival in a suitable neighbourhood of the repeller, the associated quasi-ergodic measure provides the conditioned statistics of (Lebesgue) typical trajectories that stay close to the repeller for asymptotically long times. When these quasi-ergodic measures converge to an invariant measure of the deterministic system, we say that the limiting measure is \emph{conditioned stochastically stable}. Note that while stochastically stable invariant measures are supported on attractors, conditioned stochastically stable invariant measures may be supported on repellers.

We show that uniformly expanding repellers admit a unique conditioned stochastically stable invariant measure, which corresponds to the equilibrium state associated with the geometric potential~\cite[Section 1.2.2]{Climenhaga2017} in the framework of thermodynamic formalism~\cite{Ruelle1978}.
More generally, we establish that any equilibrium state from the thermodynamic formalism on repellers\footnote{This result also applies to attractors.} is approximated by quasi-ergodic measures of so-called \emph{weighted Markov processes}, which originate from the theory of Feynman-Kac path distributions (see~\cite{delMoral2004, Champagnat2018, Champagnat2023, Kim2024} and references therein), and thus show that equilibrium states are conditioned stochastically stable in a broader sense.

\subsection{Conditioned stochastic stability}\label{sec:cond_stoch_stab}

The notion of conditioned stochastic stability that we propose is based on ideas from the theory of absorbing Markov processes~\cite{Collet1994} and
conditioned random dynamical systems~\cite{Zmarrou2007, Engel2019, Castro2021, Castro2022, Castro2023}.
As mentioned above, the statistical behaviour of a Markov process $X_n$ on a state space $M$ conditioned upon remaining outside of a subset $\partial \subset M$ is captured by its quasi-ergodic measure $\nu$ on $M \setminus \partial$~\cite{Darroch1965, Breyer1999, Zhang2014, Colonius2021}. This object describes the limiting distribution of the conditioned Birkhoff averages of $X_n$, i.e.~given an observable $h: M \to \R$ it holds that for $\nu$-almost every $x \in M\setminus \partial$,
\begin{align*}
\E_x\left[\frac{1}{n}\sum_{i= 0}^{n-1}h \circ X_i \, \bigg|\, \tau >n \right] \coloneqq \frac{1}{\P_x\left[\tau > n\right]} \E_x\left[\1_{\{\tau > n\}}\frac{1}{n}\sum_{i= 0}^{n-1} h \circ X_i\right]   \xrightarrow[]{n \to \infty} \int h(x) \nu(\dx),
\end{align*}
where conditioning upon $\tau \coloneqq \min\{i \in {\mathbb N}; X_i \in \partial\} > n$ ensures the process has not been absorbed by time $n$. 

Given a map $T:M \to M$ on a manifold $M$ and a subset $\partial \subset M$, consider the Markov process $X_n^\e$ on $M$ generated by $\e$-bounded random perturbations of $T$. Conditioned stochastic stability concerns the quasi-ergodic measures of $X_n^\e$ on $M\setminus \partial$, and their limit as the amplitude of the perturbation $\e$ goes to 0. When these quasi-ergodic measures converge to a $T$-invariant measure $\nu_0$ (in the weak$^*$ topology), we say that the limiting measure is \emph{conditioned stochastically stable} on $M \setminus \partial$.
Observe that this notion depends on the choice of random perturbation generating $X_n^\e$, which is also true for (classical) stochastic stability. 
As is common in the study of (classical) stochastic stability, we only consider random bounded diffusive perturbations~\cite{Benedicks1992, Baladi1996, Araujo2000, Benedicks2006, Alves2007} (see Section~\ref{sec:setup_notation} for the precise details), \textcolor{black}{and prove that all such random perturbations give rise to the same notion.}

In the context of uniformly expanding repellers, it is natural to assume that the repelling set is characterised by
\begin{equation}\label{eq:lambda}
\Lambda = \bigcap_{n \geq 0} T^{-n}(M \setminus \partial),
\end{equation}
where $\partial$ could be, for example, the complement of a neighbourhood of the repeller (see Section~\ref{sec:local} and Section~\ref{ex:complex_quadratic}) or
a small open neighbourhood of the attractors of $T$ (see Section~\ref{sec:global} and Section~\ref{ex:logistic_map}). In this paper, we prove the following result (see Theorem~\ref{thm:global} for a more precise and more general result):

\begin{mytheorem}[A]\label{thm:thmA}
    Given a $\mathcal C^2$ map $T$ on $M$ and a suitable open set $\partial \subset M$, with $\Lambda$ as in equation \eqref{eq:lambda}, assume that 
       \begin{enumerate}
            \item $\left. T \right|_{\Lambda}:\Lambda \to \Lambda$ is uniformly expanding, 
            \item $\Lambda \subset \mathrm{Int} (M \setminus \partial)$, and
            \item $T:\Lambda \to \Lambda$ admits a unique invariant measure $\nu_0$ known as the \emph{equilibrium state}\footnote{This is the equilibrium state associated with the (geometric or natural) potential $- \log |\det \d T|$, see Section~\ref{sec:tfandwmp}.}, which is mixing\footnote{Recall that a $T$-invariant probability measure $\mu$ is said to be \emph{mixing} if for any measurable sets $A,B$ it holds that \[\lim_{n \to \infty} \mu(T^{-n}(A) \cap B)) - \mu(A) \mu(B) = 0.\]} (see e.g.~\cite[Section 7.1]{Viana2016}).
       \end{enumerate}
    Then $\nu_0$ is conditioned stochastically stable on $M \setminus \partial$.
\end{mytheorem}

Importantly, we highlight that the closure of the set of periodic points of $\Lambda$ may not be topologically transitive (see Lemma~\ref{lem:dyndec}) and so, it is not clear whether perturbation arguments on the spectrum of the transfer operator~\cite{Keller1999, Gouezel2008} are applicable.
Moreover, this result allows for the study of stochastic stability in open systems as it provides a new perspective based on conditioned random dynamics and circumventing the lack of continuity mentioned in~\cite[Section 8.1.2]{Gouezel2008}. 

As mentioned in Theorem~\ref{thm:thmA} (3), it turns out that $\nu_0$ is a well-known object in the theory of thermodynamic formalism~\cite{Ruelle1978, Climenhaga2017} and corresponds to the unique equilibrium state on the set $R~\coloneqq~\overline{\{p \in \Lambda;\ p \text{ is } T-\text{periodic}\}}$ associated with the potential $-\log|\det \d T|$, i.e.~$\nu_0$ is the unique $T$-invariant measure satisfying
\[h_{\nu_0}(T) - \int \log|\det \d T| \, \d \nu_0 = \sup_{\mu \in \mathcal I(T, R)}\left(h_{\mu}(T) -\int \log|\det \d T| \, \d\mu\right),\]
where $h_\mu$ is the Kolmogorov-Sinai (or metric) entropy~\cite{Kolmogorov1958, Sinai1959} and $\mathcal I (T, R)$ is the set of $T$-invariant probability measures on $R$. This result has its parallel in the (classical) theory of stochastic stability. Indeed, given a uniformly hyperbolic transformation $T:M \to M$ on a compact metric space $M$, it is well known that stochastically stable invariant measures on attractors correspond to the equilibrium states from the thermodynamic formalism associated with the potential $- \log | \det \d T \left. \right|_{E^u}|$, where $E^u$ denotes the unstable expanding direction of $T$~\cite{Young2002}.

In this paper, we uncover a stronger connection between conditioned stochastic stability and the thermodynamic formalism, establishing the approximation of any equilibrium state by quasi-ergodic measures of \emph{weighted Markov processes}.

\subsection{Thermodynamic formalism and weighted Markov processes}  \label{sec:tfandwmp}
The thermodynamic formalism is a powerful framework for the analysis of statistical properties of dynamical systems. Pioneered by Sinai, Ruelle and Bowen~\cite{Sinai1972, Bowen1975, BowenRuelle1975, Ruelle1976, Ruelle1978} and motivated by the field of statistical physics, this theory aims to describe properties of equilibrium states, such as the measure of maximal entropy and other invariant Gibbs measures~\cite{Climenhaga2017, Baladi2018}.

Given a $T$-invariant set $\Lambda \subset M$, an equilibrium state on $\Lambda$ is defined for each given potential $\psi:\Lambda \to \R$ as an invariant measure $\nu^\psi$ on $\Lambda$ whose \emph{metric pressure} is equal to the \emph{topological pressure} $P(T,\psi,\Lambda)$ of the system on $\Lambda$, i.e.~$\nu^\psi$ satisfies
\begin{equation}\label{eq:top_press}
h_{\nu^\psi}(T) +\int \psi\, \d \nu^\psi = \sup_{\mu \in \mathcal I(T,\Lambda)}\left(h_{\mu}(T) +\int \psi\, \d\mu\right)\eqqcolon P(T,\psi,\Lambda).
\end{equation}
In particular, observe that when $\psi =0$ the equilibrium states associated with this potential correspond to the measures of maximal entropy~\cite[Section 10.5]{Viana2016}. Moreover, a classical result of Ruelle (see~\cite[Lemma 1.4]{Ruelle1989} or Lemma~\ref{lem:dyndec} below) provides the existence and uniqueness of equilibrium states for H\"older potentials on uniformly expanding repellers~\cite[Theorem 12.1]{Viana2016}. 

It is natural to ask whether the definition of conditioned stochastic stability can be extended to approximate other equilibrium states of $T$.
This question appears not to have been raised in the literature, even for stochastic stability of equilibrium states on attractors.
Here, we show that equilibrium states on uniformly expanding repellers are approximated by quasi-ergodic measures of weighted Markov processes~\cite{delMoral2004, Champagnat2018, Champagnat2023, Kim2024}, providing a general notion of {conditioned stochastic stability}.

Given a Markov process $X_n$ on $M$, consider a non-positive \emph{weight} function\footnote{The weight function is sometimes referred to as a ``potential'' in the literature~\cite{delMoral2004, Zhang2014}. Here, we reserve the term ``potential'' for the symbol $\psi$ in equation~\eqref{eq:top_press}.} $\phi:M \to \mathbb R_{\leq 0}$ and define the new process $X_n^\phi$ by
\begin{equation}\label{eq:weighted_proc}
X_{n+1}^{\phi} = \begin{cases}
    X_{n+1}, & \text{with probability } e^{\phi(X_n)},\\*
    \partial, & \text{with probability } 1-e^{\phi(X_n)},
\end{cases}
\end{equation}
where $\partial$ is a cemetery state. 
If $X_n$ is already an absorbing Markov process killed at $\partial^\prime$, we may (and do) set $\partial =\partial ^\prime$. We refer to the new Markov process $X_n^\phi$ as a \emph{$e^\phi$-weighted Markov process}. A quasi-ergodic measure $\nu^\phi$ provides the statistical behaviour of the process when conditioned upon survival on $M \setminus \partial$.

\begin{definition}
    We say that $\nu^\phi$ is a \textit{quasi-ergodic measure} for the $e^\phi$-weighted Markov process $X^\phi_n$ if for any observable $h:M \to \mathbb R$,
\begin{align*}
\E_x^\phi \left[\frac{1}{n}\sum_{i= 0}^{n-1}h \circ X_i^\phi \, \bigg|\, \tau^\phi > n \right] \xrightarrow[]{n \to \infty} \int h(x) \nu^{\phi}(\dx),
\end{align*}
for $\nu^\phi$-almost every $x \in M \setminus \partial$, where $\tau^{\phi}\coloneqq \min \{n \in {\mathbb N};\, X_n^{\phi} \in \partial\}$ and $\E^\phi_x$ is the expectation with respect to the weighted process $X_n^\phi$ with $X_0^\phi = x$.
\end{definition}

The random variable $\tau^\phi$ denotes the time at which the process is killed, either by dynamically entering $\partial$ (hard killing) or due to the weight $e^\phi$ (soft killing). When both are present, the conditioned Birkhoff averages simplify to (see Section~\ref{sec:setup_notation} for precise details)
\begin{equation}
\begin{split}
\E_x^\phi \left[\frac{1}{n}\sum_{i= 0}^{n-1}h \circ X_i^\phi \, \bigg|\, \tau^\phi > n \right]  &  {\coloneqq \frac{1}{\E_x[e^{S_n\phi}\1_{\{\tau > n\}}]}{\E_x\left[ e^{S_n\phi}\1_{\{\tau > n\}} \frac{1}{n}\sum_{i= 0}^{n-1}h \circ X_i \right]}}\label{eq:phiqem}\\*
&{ \xrightarrow[]{n \to \infty} \int h(x) \nu^{\phi}(\dx),}
\end{split}
\end{equation}
where $\tau = \min \{n;\, X_n \in \partial\}$ relates to hard killing and $S_n\phi \coloneqq \sum_{i= 0}^{n-1} \phi \circ X_i$ relates to soft killing. Note that when $\phi = 0$, we recover the setting introduced in the previous section.

Observe that the right-hand side of equation \eqref{eq:phiqem} is also well-defined as long as $\phi$ is measurable and bounded, even if it is occasionally positive. Indeed,
$$\frac{\E_x\left[ e^{S_n\phi}\1_{\{\tau > n\}} \frac{1}{n}\sum_{i= 0}^{n-1}h \circ X_i \right]}{\E_x[e^{S_n{\phi}}\1_{\{\tau > n\}}]} = \frac{\E_x\left[ e^{S_n\bar{\phi}}\1_{\{\tau > n\}} \frac{1}{n}\sum_{i= 0}^{n-1}h \circ X_i \right]}{\E_x[e^{S_n\bar{\phi}}\1_{\{\tau > n\}}]} $$
where $\bar{\phi} = \phi - \sup \phi_+$, with $\phi_+(x)\coloneqq\max\{\phi(x),0\}$. Defining the $e^\phi$-weighted process $X_n^\phi$ to be equal to $X_n^{\bar{\phi}}$, we recover the interpretation from equation~\eqref{eq:weighted_proc}.

Recall that $X_n^{\e}$ is a Markov process on $M$ generated by $\e$-bounded random perturbations of the map $T$ and absorbed on $\partial \subset M$, and denote by $X_n^{\e, \phi}$ the appropriate weighting of $X_n^{\e}$ as in equation~\eqref{eq:phiqem}.

\begin{definition}\label{def:css_nu0}
    We say that a $T$-invariant measure $\nu_0^\phi$ is \emph{conditioned $e^\phi$-weighted stochastically stable} if the quasi-ergodic measures $\nu^\phi_\e$ on $M \setminus \partial$ of the weighted process $X_n^{\e,\phi}$ converge to $\nu_0^\phi$ in the weak$^*$ topology as $\e$ goes to $0$.
\end{definition}

We generalise Theorem~\ref{thm:thmA} to allow for soft killing and show conditioned $e^\phi$-weighted stochastic stability with the following result (see Theorem~\ref{thm:global} for a more precise and more general statement):
\begin{mytheorem}[B]\label{thm:thmB}
    \textbf{(Main Theorem)} Given a $\mathcal C^2$ map $T$, a H\"older weight function $\phi$, and a suitable open set $\partial \subset M$, with $\Lambda$ as in equation \eqref{eq:lambda}, assume that 
       \begin{enumerate}
            \item $\left. T \right|_{\Lambda}:\Lambda \to \Lambda$ is uniformly expanding, 
            \item $\Lambda \subset \mathrm{Int} (M \setminus \partial)$, and
            \item $T:\Lambda \to \Lambda$ admits a unique equilibrium state $\nu^\psi$ associated with the potential $\psi = \phi - \log|\det \d T|$, which is mixing (see e.g.~\cite[Section 7.1]{Viana2016}).
       \end{enumerate}
    Then $\nu^\psi$ is conditioned $e^\phi$-weighted stochastically stable on $M \setminus \partial$, \textcolor{black}{i.e.~$\nu^\psi = \nu_0^\phi$ from Definition~\ref{def:css_nu0}.}
\end{mytheorem}

\begin{remark}
    \textcolor{black}{For each H\"older weight function $\phi$, we emphasise that every choice of random perturbation presented below, i.e. every family $\{T_\w\}_{\w\in \Omega}$ in Section~\ref{sec:randperturb}, generating the Markov process $X_n^\e$ identifies the same invariant measure $\nu_0^\phi$ of $T$ as the conditioned stochastically stable one.}
\end{remark}

In particular, for every repeller $R^1, \dots, R^k \subset R$ of the dynamical decomposition of $T$ (see Lemma~\ref{lem:dyndec}), 
Theorem~\ref{thm:thmB} also holds for the process conditioned upon survival in $R^{i}_\delta = M \setminus \partial, i =1, \dots, k$, where $R^i_\delta$ is a sufficiently small $\delta$-neighbourhood of the repeller $R^i$. This gives rise to the following Corollary for each $i = 1, \dots, k$ (see Theorem~\ref{thm:local} for more precise details):
\begin{mycorollary}[B1] \label{cor:corB1}
    Given a H\"older weight function $\phi$, there exists a unique $T$-invariant measure $\nu^\psi$ on $R^i$ which is conditioned $e^\phi$-weighted stochastically stable on every sufficiently small neighbourhood $R^i_\delta$. Moreover, $\nu^\psi$ is the unique equilibrium state associated with the potential $\psi = \phi - \log |\det \d T|$ on $R^i$, i.e.~$\nu^\psi = \nu_0^\phi$ \textcolor{black}{from Definition~\ref{def:css_nu0}}.
\end{mycorollary}

Observe that the main difference between both results relies on the choice of $\partial$. On the one hand, Theorem~\ref{thm:thmB} corresponds to the so-called \emph{global problem} and requires $\Lambda \subset M \setminus \partial$. On the other hand, Corollary~\ref{cor:corB1} refers the \emph{local problem} as the process is conditioned upon survival locally around some (i.e.~in a small neighbourhood of) $R^i$. 
We mention that Corollary~\ref{cor:corB1} as stated above can be derived directly from spectral stability arguments~\cite{Keller1999}. In this paper, however, we first prove Corollary~\ref{cor:corB1} using a Hilbert cone contraction argument (Section~\ref{sec:local}) which provides a more precise description of the quasi-ergodic measure later used in the proof of Theorem~\ref{thm:thmB} (Section~\ref{sec:global}). 
For the latter, we identify a graph structure representing the dynamical behaviour of $X_n^{\e,\phi}$ conditioned upon staying on $M\setminus \partial$. This construction resembles the graphs built via chain recurrence and filtration methods~\cite{Crovesier2015, deLeo2021a, deLeo2021b} and allows us to recover the setting of Corollary~\ref{cor:corB1}.

\subsection{Context of the results}
The results in this paper 
relate to the theories of open systems and spectral stability.

The theory of conditioned random dynamical systems, which builds 
on the theory of absorbing Markov processes,  essentially addresses transient properties of random dynamics. In particular, the quasi-stationary and quasi-ergodic measures of the associated absorbing Markov process provide relevant statistical properties of the random system conditioned upon remaining in a given region: the former is associated with the rate of escape of the system and the latter provides the Birkhoff averages of the process conditioned upon survival. 

This has clear parallels with the theory of deterministic open dynamical systems~\cite{Demers2006}, where conditionally invariant measures are similar to quasi-stationary measures, and invariant measures on the survival set resemble quasi-ergodic measures. 

It is important to note that a conditioned random dynamical system with diffusive-like noise typically has 
a single quasi-stationary measure and a single quasi-ergodic measure. 
This differs from the setting of deterministic open systems, where uncountably many conditionally invariant measures and invariant measures may exist.
We would like to emphasise that the quasi-ergodic measure is not a stationary measure of the conditioned random dynamical system. It is the relevant object to study ergodic theory from a conditioned point of view~\cite{Darroch1965, Breyer1999, Champagnat2016, Engel2019, Castro2022}. Indeed, quasi-ergodic measures are similar to invariant measures of open systems but their support on the state space 
usually has non-empty interior, while 
invariant measures of open systems are often supported on a Cantor-like set.

Taking the limit of noise amplitude to zero, we obtain an alternative perspective to
transient dynamics in deterministic systems, which naturally 
aligns with the theory of open dynamical systems. In addition, this conditioned random dynamics approach provides an elegant, natural means to approximate also other equilibrium states on repellers. 

In hyperbolic systems, equilibrium states may be constructed from the combination of right and left eigenfunctions of maximal eigenvalue, associated with a particular transfer operator displaying a spectral gap~\cite{Baladi2000, Baladi2018, Demers2021, Gouezel2008}. Spectral stability of such operators then ensures that properties of the peripheral spectrum are stable under suitable (abstract) perturbations~\cite{Keller1999, Gouezel2006}.
Exploiting spectral stability of transfer operators has proven to be useful in several settings, including the proof of the existence of absolutely continuous conditioned invariant measures for systems with holes~\cite{Demers2006}, the study of decay of correlations for small random perturbations of hyperbolic systems~\cite{Gouezel2006}, and the study of linear response~\cite{Baladi2014, Baladi2018}, to name a few.

In this paper, we consider Birkhoff averages of a canonical weighted Markov process whose explicit dynamics are meaningful and can be well understood. 
We then study the quasi-ergodic measure of this process, a well-established object from the theory of absorbing Markov processes, which provides relevant dynamical and statistical properties of the random system.
For the problem at hand, we show that quasi-ergodic measures can be constructed using a functional analytical approach (see Appendix~\ref{sec:appendix}). In particular, we use an elementary Hilbert cone technique~\cite{Liverani1995, Viana1997} from which we obtain a detailed description of the quasi-ergodic measure and its constituents. Finally, we show that quasi-ergodic measures converge to equilibrium states. 

{\color{black}We mention that there are several results in the literature on the zero–noise limit for absorbing Markov processes, where one studies the limiting behaviour of the quasi-stationary distribution as the noise strength tends to zero (see~\cite{FaureSchreiber14,SchreiberHuangJiangWang21,HeningQiZhongweiYi24,ProdhommeStrickler24} and the references therein). In previous works, the deterministic system is assumed to have an attractor on which the limiting quasi-stationary measures concentrate. In contrast, the present paper focuses on repellers: we consider small random perturbations of the dynamics near uniformly expanding repellers and obtain zero–noise limits supported on these repelling sets. Moreover, we treat both quasi-stationary and quasi-ergodic limits within the same framework, showing that quasi-ergodic measures generalise the role of stationary measures in the study of stochastic stability.}

It should be noted that the local (conditioned) stochastic stability results in Section~\ref{sec:local} (Theorem~\ref{thm:local}) 
align with well-known 
spectral stability arguments as mentioned above. However, the global results in Section~\ref{sec:global} (Theorem~\ref{thm:global}) require the fine description of the quasi-ergodic measure provided by Section~\ref{sec:local}, which is not accessible through standard spectral stability results. 

\subsection{Outline} This paper is organised as follows.
In Section~\ref{sec:setup_notation}, we introduce the objects of interest from the theory of conditioned random dynamics. We also lay out the required technical conditions (Hypotheses~\ref{hyp:H1} and~\ref{hyp:H2}) regarding the deterministic systems considered, their random perturbations and present the two main theorems (Theorems~\ref{thm:local} and~\ref{thm:global}). Three examples are presented in Section~\ref{sec:examples}. In Section~\ref{sec:consequences}, we explore the direct implications of the hypotheses. In Section~\ref{sec:local}, we analyse the local problem (i.e.~conditioning the random dynamics on a small neighbourhood of a repeller) and prove Theorem~\ref{thm:local}. In Section~\ref{sec:global}, we consider the global picture (i.e.~conditioning upon not escaping from a general neighbourhood of the repeller) and prove Theorem~\ref{thm:global}.
We provide examples in Section~\ref{sec:examples} where these theorems are applicable. Finally, we devote Appendix~\ref{sec:appendix} to a general proof for the existence of (weighted) quasi-ergodic measures, simplifying previous techniques.

\section{Setup, main results and examples}\label{sec:setup_notation}

We begin with a brief recollection of the basic concepts in the theory of conditioned random dynamics as introduced in~\cite{Castro2021, Castro2023}. Consider a Markov chain $X_n$ evolving in a metric space $(E,d)$ and let $Y \subset E$ be a compact subset. We are interested in studying the behaviour of a Markov chain as it evolves in $Y$, {we} condition upon remaining in $Y$, and kill the process as soon as it leaves this subset. We thus identify $E\setminus Y$ with a ``cemetery state'' $\partial$ and consider the space $E_Y \coloneqq Y \sqcup \partial$ with the induced topology. Throughout this paper, we assume that 
\[X\coloneqq \left(\Omega, \{\mathcal F_n\}_{n\in \N_0}, \{X_n\}_{n\in \N_0}, \{\mathbf P^n\}_{n\in \N_0}, \{\P_x\}_{x\in E_Y}\right)\]
is a Markov chain with state space $E_Y$, in the sense of~\cite[Definition III.1.1]{Rogers1994}. Hard killing, or absorption, on $\partial$ means that $\mathbf P(\partial, \partial) = 1$. {We} define the {(dynamical)} stopping time
$\tau \coloneqq \inf\{n \in \N;\, X_n \in \partial\}.$ 

{Consider an $\alpha$-H\"older weight function $\phi: Y \to \mathbb R$, \textcolor{black}{$\alpha > 0$}, and define the $e^\phi$-weighted process $X_n^\phi$ as in Section~\ref{sec:tfandwmp}, equation~\eqref{eq:weighted_proc}. For this process, we define the stopping time $\tau^\phi \coloneqq \min\{n \in \mathbb N; \, X_n^\phi \in \partial\}$, providing the time at which $X_n^\phi$ enters $\partial$ either dynamically (hard killing) or due to the weight function $\phi$ (soft killing).}

{Observe that the weighted process $X_n^\phi$ has transition \textcolor{black}{kernels} given by $\mathbf P^\phi(x, \d y) = e^{\bar{\phi}(x)}\mathbf{P}(x, \d y)$ for all $x  \in Y$, recall that $\bar{\phi} = \phi-\sup \phi_+.$ Moreover, \eqref{eq:weighted_proc} naturally induces a filtered space $(\Omega^\phi,\{\mathcal F_n^\phi\}_{n\in\mathbb N_0})$ and a family of probability measures $\{\mathbb P_x^\phi\}_{x \in E_Y}$ which makes $X_n^\phi$ a Markov process (see~\cite[Section III.7]{Rogers1994} for such a construction). Finally, we denote by $\mathbb E_x$ and $\mathbb E_x^\phi$ the expectation with respect to $\mathbb P_x$ and $\mathbb P_x^\phi$, respectively.}

Under an irreducibility condition of $X_n$ on $Y$~\cite{Castro2021}, the process almost surely escapes this set, implying that the system's long-term behaviour is characterised by a stationary delta measure sitting on the cemetery state. To understand the dynamics of the process before escaping from $Y$ {one} generalise{s} the notion of stationary measures to that of quasi-stationary measures~\cite{Darroch1965, Breyer1999, Collet1994, Cattiaux2009}.

\begin{definition}\label{def:quasistat}
    Given a bounded and measurable function $\phi:Y\to \mathbb R$, {we say that a} Borel probability measure $\mu$ on $Y$ is {a \emph{quasi-stationary measure} of} the {weighted Markov process $X_n^\phi$} if 
    $$\int_Y e^{\phi(y)}\mathbf P(y,\d x)\mu(\d x) = \lambda^\phi \mu(\d x)$$
    and $\lambda^\phi =\int_Y e^{\phi(x)}\mathbf P (x, Y)\mu(\d x)>0$ is the \emph{growth rate} of $\mu$ for $X_n^\phi$ on $Y$. Observe that when $\phi=0$ we recover the classical definition of quasi-stationary measure~\cite[Definition 2.1]{Collet2013}.
\end{definition}

\begin{remark}
    Note that in the usual setting of absorbed Markov processes with no weight function, i.e.~$\phi = 0$, and only hard killing, $\lambda^\phi\leq1$ is called the \emph{survival rate} and denotes the probability that the process is not killed in the next iterate when distributed according to $\mu$.
\end{remark}

We recall that quasi-stationary measures are not the {relevant measures} to consider when studying conditioned Birkhoff averages~\cite{Darroch1965, Breyer1999, Champagnat2016, Castro2021}, as these measures do not perceive how likely it is for a point to remain indefinitely in $Y$. Instead, this information is provided by the so-called quasi-ergodic measure. 

\begin{definition}
    {A probability measure $\nu$ on $Y$ is a \emph{quasi-ergodic measure} of the $e^\phi$-weighted Markov process $X_n^\phi$ if for any bounded measurable function $h:Y\to \mathbb R$ it holds that
    \[\lim_{n\to \infty} \E_x^\phi \left[\frac{1}{n}\sum_{i = 0}^{n-1}  h \circ X_i^\phi \, \bigg|\, \tau^\phi >n \right]\ = \int_Y h(y)\nu(\d y)\quad \text{for $\nu$-almost every }x \in Y.\]}
\end{definition}
{If $X_n^\phi$ has both hard and soft killing, then for every $n \in \mathbb N$
\begin{align*}
\E_x^\phi \left[\frac{1}{n}\sum_{i= 0}^{n-1}h \circ X_i^\phi \, \bigg|\, \tau^\phi > n \right]  &  {= \frac{1}{\E_x[e^{S_n\phi}\1_{\{\tau > n\}}]}{\E_x\left[ e^{S_n\phi}\1_{\{\tau > n\}} \frac{1}{n}\sum_{i= 0}^{n-1}h \circ X_i \right],}}
\end{align*}
where $S_n \phi \coloneqq \sum_{i=0}^{n-1} \phi \circ X_i$ is the Birkhoff sum.}

While showing the existence of quasi-stationary measures relates to solving an eigenfunctional equation and can be approached using fixed point arguments (see~\cite[Theorem 4]{timur2005} and~\cite[Proposition 2.10]{Collet2013}), this is not the case for quasi-ergodic measures and proving their existence and uniqueness is not straightforward. Indeed, this involves characterising the non-trivial limit of a conditional expectation that requires rigorous techniques in functional analysis and probability theory~\cite{Champagnat2016, Zhang2014, Castro2021}. We devote the Appendix~\ref{sec:appendix} to address this question in our setup. 

From here onwards, let $(M,\langle \cdot, \cdot \rangle)$ be an orientable Riemannian compact manifold, possibly with boundary and let $U\subset M$ be an open subset. Without loss of generality, we may assume that $M$ is embedded in an orientable boundaryless compact manifold $E$ of the same dimension and endowed with a Riemannian metric whose restriction to $M$ coincides with $\langle\cdot,\cdot\rangle$ (in the case that $M$ is without boundary, we assume that $E=M$). Since this will be clear by context, we may also write the Riemannian metric of $E$ as $\langle\cdot,\cdot\rangle$. The manifold $E$ should be thought of as an ambient space for $M$ and a mere theoretical artefact since it does not play a major role in applications, while $U$ {may} be interpreted as an open hole in the system.

\begin{notation}\label{not:notation}
Throughout this paper, we use the following notation:
\begin{enumerate}[label = (\roman*)]
    \item  Given $x \in E$ and $v\in T_x E$, define $\|v\|_x \coloneqq \sqrt{\langle v,v\rangle_x}$ as the natural norm on $T_x M$.
    \item We denote by $\mathrm{dist}(\cdot, \cdot)$ the distance on $E$ induced by the Riemannian metric $\langle \cdot,\cdot  \rangle$.
    \item As usual, we write $\rho$ for a Borel measure on $E$ induced by a smooth volume form $V_{E}$ compatible with $\langle \cdot,\cdot \rangle$.
    \item We denote by $\mathcal C^k(E)$ the space of continuous functions with $k$ continuous derivatives on $E$ and use $\mathcal M(E)$ to denote the space of signed Borel finite measures on a $E$. Given a non-negative measure $\rho \in \mathcal M(E)$, we denote by $L^k(E, \rho)$ the space of functions with finite $k$-th $\rho$-moment (although $\rho$ may be omitted when it is the reference measure). $\mathcal C^k_+(E), L^k_{+}(E)$ and $\mathcal M_+(E)$ denote the respective subsets of non-negative functions and measures on $E$.
    \item Given a $\mathcal C^1$ function $G:E\to E$ and $x\in E,$ we denote its determinant by 
    \[\det \d G(x) =  \frac{V_E(G(x))(\d G(x)v_1,\ldots, \d G(x) v_{\dim E})}{V_E(x)(v_1,\ldots,v_{\dim E})} ,\]
    for any (and therefore all) orthonormal basis $\{v_1,\ldots,v_m\}$ of $T_x M.$
    \item\label{not:exp-1} Given a set $A\subset E$ we denote its closed neighbourhood of radius $\delta>0$ by $A_\delta = \overline{B_\delta (A)}\coloneqq\{x \in E;\,\dist(x, a) \leq \delta \text{ for some } a \in A\}$.
\end{enumerate}
\end{notation}

\subsection{The deterministic dynamics}
\textcolor{black}{Let $T: E\to E$ be a map such that $\left.T\right|_{E\setminus \partial}$ is $\mathcal C^2$, for a suitable choice of $\partial$, e.g. $\partial$ may be an open subset $U \subset M \subset E$ where the process is killed or an artificial cemetery state in the absence of hard killing. For an invariant set $\Lambda$ as in equation~\eqref{eq:lambda}, i.e.
\[\Lambda = \bigcap_{n \geq 0}T^{-n}(M \setminus \partial),\]
we consider the following hypothesis.}


\begin{hypothesis}[H1]\label{hyp:H1} There exists a compact $T$-invariant set $\Lambda\subset E$ that is uniformly hyperbolic expanding, i.e.~there exists $r >0$ such that for all $x \in \Lambda$,
\begin{align}
   |\d T^{n}(x)^{-1}\|< \frac{1}{(1+r)^n} \quad \text{for every } n\geq 1,\label{eq:uh} 
\end{align}
and there exists a neighbourhood $V$ of $\Lambda$ in $E$ such that $T^{-1}(\Lambda) \cap V = \Lambda$. \textcolor{black}{We call $\Lambda$ a (uniformly expanding) repeller.}
\end{hypothesis}
\begin{remark}
\textcolor{black}{Observe that if there exists \(C>0\) and $r>0$ such that, for all $x\in\Lambda$ and every $n\ge1$,
\begin{align}
   \|\d T^{n}(x)^{-1}\|\le C\frac{1}{(1+r)^n} \quad\text{for every }n\ge1.\label{eq:uh1}
\end{align}
Then, after a suitable change of the Riemannian metric on $M$, one can arrange that \eqref{eq:uh} holds (see \cite[Proposition 4.2]{Shub1987}).}
\end{remark}

\textcolor{black}{Uniformly expanding sets admit a well-known ``spectral'' or ``dynamical decomposition''~\cite[Theorem~19.3.6]{HasselblattKatok2003}, providing the fundamental sets on which we shall perform our local analysis of conditioned stochastic stability. We denote these by $R^i, i = 1, \dots, k$, and recall that they are given by the following result.}
\begin{lemma}\label{lem:dyndec}
Let $T$ and $\Lambda$ satisfy Hypothesis~\ref{hyp:H1} and consider the set $$R \coloneqq \overline{\mathrm{Per}(T)} \coloneqq \overline{\{p \in \Lambda;\ p \text{ is a }T\text{-periodic point}\}}.$$
Then there exists a (finite) partition of $R$ in non-empty compact sets $R^{i,j},$ with $1\leq i \leq k$ and $1\leq j \leq m(i)$, such that
\begin{enumerate}[label = (\arabic*)]
    \item $R^{i} = \cup_{j=1}^{m(i)} R^{i,j}$ is a $T$-invariant set for every $i$,
    \item $T(R^{i,j}) = R^{i,j+1\ (\mathrm{mod}\ m(i))}$ for every $i,j$,
    \item $T:R^i\to R^i$ is uniformly hyperbolic and topologically transitive, and
    \item each $T^{m(i)}:R^{i,j}\to R^{i,j}$ is uniformly hyperbolic and topologically exact.
\end{enumerate}
Furthermore, the number $k$, the numbers $m(i)$ and the sets $R^{i,j}$ are unique up to renumbering.
\end{lemma}
\begin{proof}
See, for example,~\cite[Theorem~18.3.1]{Katok1995} or~\cite[Theorem~11.2.15]{Viana2016}.
\end{proof}


\textcolor{black}{Throughout, we let $\phi: M \to \mathbb R$ be a H\"older function, which we may call the \emph{weight function}, and say that the triple $(T,\phi,\Lambda)$ satisfies Hypothesis~\ref{hyp:H1}, although $\phi$ does not play a role in this assumption. We recall in the following theorem a well-known result of Ruelle~\cite[Lemma~1.4]{Ruelle1989} (see also~\cite[Chapters 7.26-7.31]{Ruelle1978}) which provides the existence and uniqueness of equilibrium states associated with the potential $\psi = \phi - \log |\det \d T|$ on each $R^i, i = 1, \dots, k$ (recall Section~\ref{sec:tfandwmp}), for triples $(T, \phi, R^i)$ satisfying Hypothesis~\ref{hyp:H1}.}

\begin{theorem}[Ruelle]\label{thm:ruele} Let $T$ satisfy Hypothesis~\ref{hyp:H1}. Let $R^1,\ldots, R^k$ be as in Lemma~\ref{lem:dyndec} and fix $i\in\{1,\ldots,k\}$. For every $\alpha$-H\"older potential $\psi: R^i\to \R$, $\alpha> 0$, {consider the operator} 
\begin{align*}
    \mathcal L:\mathcal C^0(R^i)&\to \mathcal C^0(R^i)\\*
    f&\mapsto \sum_{T(y)= x} e^{\psi(y)}f(y),
\end{align*}
and denote by $\mathcal L^*:\mathcal M(R^i) \to \mathcal M(R^i)$ its dual (recall Notation~\ref{not:notation}). {Then there exist a unique $m\in \mathcal C_+^0(R^i)$, $\gamma \in \mathcal M_+(R^i)$ and $\lambda>0$ satisfying}
\begin{enumerate}[label = (\arabic*)]
    \item $\ker(\calL - \lambda ) = \s(m)$,
    \item $\ker(\calL^* -\lambda) = \s(\gamma)$ and $\int_{R^i} m(x) \gamma(\dx) = 1$, and
    \item $\log \lambda= \log r(\calL) = h_{\nu} + \int \phi(x)\nu(\d x)$, {where $\nu(\d x) = m(x) \gamma(\d x)$.}
\end{enumerate}
In this context, $\nu$ is the unique $T$-invariant equilibrium state for the potential $\psi$ on $R^i$.
\end{theorem}

\textcolor{black}{In the study of local conditioned stochastic stability, we shall set $\Lambda = R^i$ for a fixed $i \in \{1, \dots, k\}$ and consider small random perturbations of the dynamics as long as the process remains in a suitable $\delta$-neighbourhood of this invariant set, i.e.~$\partial = M \setminus R^i_\delta$. As mentioned above, Theorem~\ref{thm:ruele} provides the existence and uniqueness of equilibrium states, each associated to a particular H\"older potential $\psi$, on every repeller $R^i$, $i \in \{1, \dots, k\}$. These equilibrium states are the objects we prove to be conditioned stochastically stable.}

\textcolor{black}{The study of \emph{global} conditioned stochastic stability, as motivated in the Introduction, begins by considering a more general open cemetery state, hole, or absorbing region, $U \subset M$. The set of points that never enters this region $U$ under the dynamics, similarly to equation~\eqref{eq:lambda}, is given by
\begin{equation}\label{eq:survival}
    \Lambda \coloneqq \bigcap_{n \geq 0} T^{-n}(M\setminus U).
\end{equation}
}
\textcolor{black}{The main working assumption in this global setting ensures that $T|_{\Lambda}$, with $\Lambda$ as in~\eqref{eq:survival}, admits a unique equilibrium state for each $\alpha$-H\"older weight function $\phi$. This is detailed in the following hypothesis.}

\begin{hypothesis}[H2]\label{hyp:H2}
We say that $(T, \phi, \Lambda)$ satisfies Hypothesis~\ref{hyp:H2} for the open hole $U$ if the following holds:
\begin{enumerate}[label = (\roman*)]
    \item\label{it:hyp2-1} $\Lambda \coloneqq \bigcap_{n \geq 0} T^{-n}(M\setminus U)$ is a uniformly expanding set,
    \item\label{it:hyp2-2} $T$ admits a unique equilibrium state for the potential $\psi = \phi - \log |\det \d T|$ on $\Lambda,$ and
    \item\label{it:hyp2-3} there exists $\delta>0$ such that $T^{-1}(\Lambda_\delta)\cap M_\delta \subset \Lambda_\delta$ and $M_\delta \setminus M$ has no $T$-invariant subsets.
    \end{enumerate}
We set the cemetery state $\partial \coloneqq  U \cup (E \setminus M_\delta)$.
\end{hypothesis}

\begin{remark}\label{rmk:AxiomA}
   If $M = [0,1]$, items~\ref{it:hyp2-1} and~\ref{it:hyp2-2} of Hypothesis~\ref{hyp:H2} are equivalent to the Axiom A (see e.g.~\cite[Chapter 3.2.b]{deMelo_vanStrien1993}).
\end{remark}
\begin{remark} 
   In the case that $M$ is a manifold without boundary, then $M_\delta = M = E$ and $\Lambda_\delta \subset M.$ Observe that item~\ref{it:hyp2-3} always holds true in dimension one. In fact, recall that $M \subset E$, where $E$ acts as an ambient space. With item~\ref{it:hyp2-3}, we ensure that there are no invariant subsets near $M$, which could trap the perturbed dynamics.
\end{remark}
\begin{remark}
       Observe that Hypothesis~\ref{hyp:H2} implies Hypothesis~\ref{hyp:H1}.
\end{remark}

\subsection{The random perturbations}\label{sec:randperturb}
When $T$ satisfies Hypothesis~\ref{hyp:H2}, given a finite $\e>0$ we consider the {random} perturbation of the {form}
\(F_\e:[-\e,\e]^m \times E\to E, \)
where $F_\e(\w, \cdot) \in \mathcal C^2(E\setminus U; E)$ and $\partial_{\w} F_\e(\w, x)$ is surjective for all $\w \in [-\e,\e]^m$. Moreover, we assume that $\mathrm{dist}_{\mathcal C^2}(F_\e(\w,\cdot), T)\leq {C}\|\w\|$ for some  $C>0$, where $\mathrm{dist}_{\mathcal C^2}$ denotes the metric on $\mathcal C^2(E \setminus U,E)$ which generates the $\mathcal C^2$-Whitney topology~\cite[Chapter 1.2]{Palis1982}. In particular, surjectivity of $\partial_\w F_\e(\w, x)$ implies $m \geq \dim E$. We note that this type of random perturbation is natural and commonly considered~\cite{Benedicks1992, Baladi1996, Araujo2000, Benedicks2006, Alves2007}.

Let $\Omega_\e \coloneqq ([-\e,\e]^m)^{\mathbb N}$ be the space of semi-infinite sequences of elements in $[-\e,\e]^m$ endowed with the probability measure $\P_\e \coloneqq (\left.\Leb\right|_{[-\e, \e]^m}/(2\e)^m)^{\otimes \N}$, {and let} $\E_\e$ {denote} the corresponding expectation with respect to $\P_\e$. For every $\w \in \Omega_\e$, $\w = (\w_0 \w_1 \ldots)$, we define $T_{\w}(x)\coloneqq T_{\w_0}(x) \coloneqq F_\e(\w_0,x)$ and $T^n_\w(x) \coloneqq T_{\w_{n-1}} \circ \dots \circ T_{\w_{n}}(x)$ for every $n \in \mathbb N$, \textcolor{black}{e.g.~we may consider by abuse of notation $X_{n+1} = F_{\w_n}(X_n)$ for an identically, independently distributed sequence of random variables $(\w_n)_{n \geq 0}$ following a uniform law on $[-\e,\e]^m$.}

As mentioned in the Introduction, Corollary~\ref{cor:corB1} applies to suitable $\delta$-neighbourhoods of each repeller $R^i$, $1 \leq i \leq k$, in the dynamical decomposition of Lemma~\ref{lem:dyndec}. In particular, \textcolor{black}{by the definition of $F_\e$ there exist} $\delta,\e_0>0$ such that the following holds true. 
\begin{lemma}\label{lem:dyndec_short}
Under Hypothesis~\ref{hyp:H1}, for every $\delta > 0$ small enough there exists $\e_0 \coloneqq \e_0(\delta) > 0$ such that for every $0 \leq \e \leq \e_0$ we have that:
\begin{enumerate}[label = (\arabic*)]
    \item\label{it:dyndecshort1} $R_\delta = R^1_\delta \sqcup \ldots \sqcup R_\delta^k$, and
    \item\label{it:dyndecshort2}  \(\sup_{x\in R^i_\delta}\P_\e[\w \in \Omega_\e;\  T_\w (x) \in R^j_\delta] = 0\) for every $i\neq j \in \{1,\ldots, k\}.$
\end{enumerate}
\end{lemma}
\begin{proof}
    \textcolor{black}{We omit this proof as it follows from standard arguments.} 
\end{proof}

\textcolor{black}{To establish the existence of quasi-ergodic measures, we exploit the properties of a stochastic analogue to (the dual of\footnote{This may become clearer in the following section, particularly after introducing Notation~\ref{not:operators}.}) Ruelle's transfer operator $\mathcal L$ presented in Theorem~\ref{thm:ruele}. For each $i \in \{1, \ldots, k\}$ and every $\alpha$-H\"older function $\phi:R^i_\delta \to \mathbb R$ (see Notation~\ref{not:notation} item~\ref{not:exp-1}), we define the annealed Koopman operator
\begin{align*}
 \calP_{\e} : f \mapsto e^{\phi(x)}\E_{\e}[ f\circ T_\w (x)\cdot \1_{R^i_\delta}\circ T_\w (x)],
\end{align*}}
\textcolor{black}{for $f$ in a suitable domain. In other words, for the absorbing and weighted Markov process $X^\phi$ starting at $x \in \mathbb R^i_\delta$ and defined by
\[
X_{n+1}^\phi(\w, x) = \begin{cases}
X_{n+1} \coloneqq T^{n+1}_\w(x), & \text{with probability } e^{\phi(X_n)},\\
\partial, & \text{with probability } 1 - e^{\phi(X_n)},
 \end{cases}
 \]
with $\partial$ the complement of $R^i_\delta$, we have that
\[\calP_{\e} : f \mapsto e^{\phi(x)}\E_{x}[ f\circ X^\phi_1\cdot \1_{R^i_\delta}\circ X^\phi_1].
\]
}

\textcolor{black}{For the global picture, we consider H\"older weights $\phi:M \to \mathbb R$ and define the (global) annealed Koopman operator given by
\begin{align*}
 \calP_{\e}f(x) \mapsto e^{\phi(x)}\E_{\e}[ f\circ T_\w (x)\cdot \1_{M\setminus U}\circ T_\w (x)],
\end{align*}
for $f$ in a suitable domain, where $U$ is the open hole from Hypothesis~\ref{hyp:H2}.}

\subsection{Main results}
The main results of this paper are as follows:

\begin{restatable}[]{theorem}{localthm}\label{thm:local}
    Assume Hypothesis~\ref{hyp:H1} and let $\delta>0$ be small enough. Given $1\leq i\leq k$ and an  $\alpha$-H\"older function $\phi:R^i_\delta \to \R$, the following properties hold for $\e>0$ sufficiently small:
    \begin{enumerate}[label = (\arabic*)]
        \item\label{it:local-1} the $e^\phi$-weighted Markov process $X_n^{\e,\phi}$ admits a unique quasi-stationary measure $\mu_\e$ on $R_\delta^i$ such that $\Lambda\subset \supp \mu_\e$,
        \item\label{it:local-2} let $\lambda_\e$ be the growth rate of $X_n^{\varepsilon,\phi}$ on $R_\delta^i$, then $\lambda_\e$ is equal to the spectral radius of  $\calP_\e:L^\infty(R^i_\delta,\rho) \to L^\infty(R^i_\delta,\rho),$ and $\log(\lambda_\e) \to P(T,\phi -\log |\det \d T|, R^i)\ \text{as}\ \e \to 0,$
        \item\label{it:local-3} the operator ${\calP_\e: L^{\infty}(R_{\delta}^i, \rho) \to L^{\infty}(R_{\delta}^i, \rho)}$ admits a unique positive eigenfunction $g_\e\in L^{\infty}(R_{\delta}^i, \rho)$ associated with the eigenvalue $\lambda_\e$,
        \item\label{it:local-4} \textcolor{black}{the $e^\phi$-weighted process $X_n^{\e,\phi}$ on $\{g_\e >0\}$ admits a unique quasi-ergodic measure $\nu_\e(\d x)$ on $\{g_\e > 0\}$,}
        \item\label{it:local-5} \textcolor{black}{$\nu_\e(\d x)\to \nu_0(\d x)$ in the weak$^*$ topology as $\e \to 0$, and}
        \item\label{it:local-6} $\nu_0$ is the unique $T$-invariant equilibrium state for the potential $\phi -\log|\det \d T|$ on $R^i$.
    \end{enumerate}
If the measure $\nu_0$ is mixing for the map $T: R^i\to R^i$, then the measure $\nu_\e$ is also a quasi-ergodic measure of the $e^\phi$-weighted Markov process $X_n^\e$ on $R_\delta^i$.
\end{restatable}

\textcolor{black}{The \emph{local} Theorem~\ref{thm:local} under Hypothesis~\ref{hyp:H1} is a version of the following \emph{global} theorem under Hypothesis~\ref{hyp:H2}:}

\begin{restatable}[]{theorem}{globalthm}\label{thm:global}
 Assume Hypothesis~\ref{hyp:H2}, let  $\nu_0$ be the unique $T$-invariant equilibrium state for the potential $\phi -\log|\det \d T|$ on $\Lambda$, and $\delta>0$ be small enough. Given an  $\alpha$-H\"older function $\phi:M_\delta \setminus U \to \R$,  the following properties hold for $\e>0$ sufficiently small:
    \begin{enumerate}[label = (\arabic*)]
        \item\label{it:thm-exp1} the $e^\phi$-weighted Markov process $X_n^{\e,\phi}$ admits a unique quasi-stationary measure $\mu_\e$ on $M_\delta\setminus U$ such that $\supp \nu_0\subset \supp  \nu$,
        \item\label{it:thm-exp2} let $\lambda_\e$ be the growth rate of $X_n^{\varepsilon,\phi}$ on $M_\delta\setminus U$, then $\lambda_\e$ is equal the spectral radius of $\calP_\e:L^\infty(M_\delta\setminus U) \to L^\infty(M_\delta\setminus U) $, and $\log(\lambda_\e) \to P(T,\phi -\log |\det \d T|, \Lambda)$ as $\e \to 0,$
        \item\label{it:thm-exp3} the operator $\calP_\e: L^{\infty}(M_\delta\setminus U, \rho) \to L^{\infty}(M_\delta\setminus U, \rho)$ admits a unique positive eigenfunction $g_\e\in L^{\infty}(M_\delta\setminus U,\rho)$ associated with the eigenvalue $\lambda_\e$,
        \item\label{it:thm-exp4} \textcolor{black}{the $e^\phi$-weighted process $X_n^{\e,\phi}$ on $\{g_\e >0\}$ admits a unique quasi-ergodic measure, $\nu_\e(\d x)$ on $\{g_\e > 0\} \cap \supp \mu_\e$, and}
        \item\label{it:thm-exp5} $\nu_\e(\d x) \to \nu_0(\d x)$ in the weak$^*$ topology as $\e \to 0$.
    \end{enumerate}
If $\nu_0$ is mixing for the map $T: R\to R$, then the conclusions of the above theorem remain true when changing the set $\{g_\e >0\}\cap \supp \mu_\e$ by $M_\delta \setminus U$. Additionally, if $\nu_0$ is mixing and $\supp \nu_0 \subset \mathrm{Int}(M\setminus U)$, then (4) is also true on the set $M \setminus U$.
\end{restatable}

\subsection{Examples}\label{sec:examples}
\textcolor{black}{Let us provide three examples to illustrate Hypothesis \ref{hyp:H1} and \ref{hyp:H2}, along with the main results obtained from Theorem~\ref{thm:local} and Theorem~\ref{thm:global}.}

\subsubsection{The logistic map}\label{ex:logistic_map}
Consider the Markov process $X_{n+1}^\e = T(X_{n}^\e) + \w_n$, $n \in \N,$ with $T(x) = ax(1-x)$ and $\w_n \sim \text{Unif}(-\e, \e)$. Fix $a = 3.83$ so that the deterministic dynamical system (with $\e = 0$) has an almost sure global three-periodic attractor~\cite{Crutchfield1982, Sato2018}, i.e.~Lebesgue almost every initial condition in $[0,1]$ is attracted to the unique three-periodic hyperbolic attractor $\mathcal A = \{p, T(p) , T^{2}(p)\}$, with $p \approx 0.1456149$ (see Remark~\ref{rmk:AxiomA}). 

The dynamical decomposition of Lemma~\ref{lem:dyndec} yields two invariant sets: the origin $R^1 = \{0\}$, and a hyperbolic Cantor set $R^{2}$ consisting of the closure of all periodic points in $(0,1)$ that are not in the basin of attraction $B(T)$ of $\mathcal A$~\cite{Smale1976}. Let $\Lambda \coloneqq [0,1] \setminus B(T)$ and $U \supset \mathcal A$ be a small enough neighbourhood of the attractor such that $U \cap \Lambda = \emptyset$.
We consider the family of $\alpha$-H\"older weight functions $\phi_t:[0,1] \to \mathbb R, x \mapsto (-t+1)\log |T'(x)|$ for $t \geq 0$. Recall that an equilibrium state $\nu_i$ associated with the potential $\phi_t - \log|a(1-2x)|$ for $T$ on $R^i$ is a measure maximising
\[\mu \mapsto  h_{\mu}(T) + \int (\phi_t  - \log|T'|)\d \mu = h_{\mu}(T) -t \int \log|T'|\d \mu,\]
where $h_\mu$ is the metric entropy and $\mu \in \mathcal I(T,R^i)$, the set of $T$-invariant measures on $R^i$.

It is well known that $\Lambda$ is a hyperbolic (uniformly expanding) invariant set~\cite{deMelo_vanStrien1993} and $T$ admits a unique equilibrium state associated with the potential $\phi_t(x) -\log|T'(x)|$ on $\Lambda$ (see e.g.~\cite[Chapters 11 and 12]{Viana2016}). Therefore, Hypothesis~\ref{hyp:H2} is satisfied and we can apply the theory developed above. For $R^1$, it is clear that $\nu_1 = \delta_0$ and $P(T, \phi_t - \log|T'|, R^1) = -t \log |a|$. For $R^2$, 
\[P(T, \phi_t - \log|T'|, R^2) = h_{\nu_2}(T) - t \int \log|a(1 - 2x)| \nu_2(\dx) > -t\log |a|,\]
since $-\log |a(1-2x)|$ reaches its minimum at $0$ and $h_{\nu_2} = (1 + \sqrt{5})/2$ (see~\cite{Smale1976} for precise details). Therefore,
\[\log \lambda_\e = \log r(\mathcal P_\e) \xrightarrow{\e\to0} P(T, \phi_t - \log|T'|, \Lambda) = P(T, \phi_t - \log|T'|, R^2) = \log \lambda_2,\]
with $\calP_\e: L^{\infty}([0,1] \setminus U) \to L^{\infty}([0,1] \setminus U)$ the global annealed Koopman operator. It follows from Theorem~\ref{thm:global} that the unique equilibrium state sits on the invariant Cantor set repeller $R^2$ and can be approximated by quasi-ergodic measure $\nu_\e$ of the $\phi_t$-weighted Markov process $X_n^\e$ on a neighbourhood of $R^2$, as $\e \to 0$.

For the particular choice of $t = 1$, i.e.~$\phi_{t=1} = 0$, the system is no longer spatially weighted, and we recover the so-called ``natural measure'' of the repeller~\cite{Kantz1985}. We note that the relationship between limiting quasi-ergodic measures and natural measures in the case of the zero weighting was previously discussed in~\cite{Bassols2023} for this example.

Finally, consider the potential $\phi_0(x)= \log|T'|$. The topological pressure of the deterministic system on $\Lambda$ is given by \(P(T, 0 , \Lambda) = h_{\nu}(T),\) where $\nu$ is the unique equilibrium state. Since this measure maximises $P(T, 0, \Lambda)$, it coincides with the measure of maximal entropy of the system.

\subsubsection{The complex quadratic map} \label{ex:complex_quadratic}
Similarly to the previous example, let us consider random perturbations of iterates of the complex quadratic map $p_c(z) = z^2 + c$, $c \in \C$, acting on the Riemann sphere $\widehat {\mathbb C} = \mathbb C \sqcup \{\infty\}$. As before, we study the Markov process $X_{n+1}^\e = p_c(X_n^\e) + \w_n$, where $\{\w_n\}_n$ are i.i.d. random variables uniformly distributed on $\{a+i b\in\mathbb C;\, (a,b)\in[-\e,\e]^2\}$, with $\e>0$ small enough. 

Consider the Julia set \( J \subset \widehat{\mathbb{C}} \) associated with the polynomial \( p_c \). Recall that \(J\) is the closure of the set of repelling periodic points~\cite[Theorem 11.1]{Milnor2006}. The set \( J \) is non-empty, compact, and totally invariant, meaning that \( J = p_c(J) = p_c^{-1}(J) \) (see~\cite[Lemma 3.1]{Milnor2006}). Now, let \( c \) be a hyperbolic complex number within the Mandelbrot set, which ensures that \(J\) is hyperbolic, i.e., $J$ is connected and satisfies \( \|p_c'(z)\| = \|2z\| > 1 \) for every \( z \in J \). 

In this context, it is readily verified that $p_c$ admits a finite attractor  \(\mathcal{A}\subset \mathbb C\). Moreover, for any \(\alpha\)-H\"older weight function \(\phi: \mathbb{C} \to \mathbb{R}\), \(p_c\) satisfies Hypothesis~\ref{hyp:H2} with \(T= p_c\), \(\Lambda= J\) and \(E= M= \widehat{\mathbb C}\). Furthermore, notice that the unique equilibrium state of \(p_c\) for the potential \(\phi - \log|\det \d p_c|\) on \(J\) is mixing.

Finally, from Theorem~\ref{thm:global}, for any \(\alpha\)-H\"older potential \(\psi: \widehat{\mathbb C} \to \mathbb{R}\), the unique \(p_c\)-invariant equilibrium state for the potential \(\psi\) on \(J\), can be approximated in the weak\(^*\) topology by quasi-ergodic measures of the $(\psi + \log|\det \d p_c|)$-weighted Markov process \(X_n^\e\) on \(\widehat{\mathbb{C}}\setminus U\), where \(U\) is a neighbourhood of \(\mathcal{A}\) such that \(U\cap J = \emptyset.\)

\subsubsection{\textcolor{black}{The Boole map}} \textcolor{black}{Consider the Boole map $T:\mathbb S^1 \to \mathbb S^1$ on the circle given by~\cite{Boole}
$$T(x) = \begin{cases}
    \frac{x(1-x)}{1-x-x^2},& x \in [0,1/2) \\
    1 - T(1-x),& x \in [1/2, 1),
\end{cases}$$
which consists of two branches
This map is not uniformly expanding since $T'(0) =1$, hence does not satisfy either of the Hypotheses~\ref{hyp:H1} nor \ref{hyp:H2} for $\partial = \emptyset$. Nevertheless, this is the only neutral fixed point and $T$ is expanding elsewhere.}

\textcolor{black}{Consider the open hole $U_s = \{x\in\mathbb S^1;x\in [0,s)\ \text{or }x\in(1-s,1] \}$ for $s\in(0,1/8)$, and let
$$ \Lambda^s = \bigcap_{n\in \mathbb N} T^{-n}(\mathbb S^1\setminus U_s)$$
be the set of points in $\mathbb S^1 \setminus U_s$ that are never mapped into $U_s$.
Observe that for each $s\in (0,1/8)$, $\Lambda^s\neq \emptyset$, since it is easy verify that there exists a $2$-periodic orbit lying in $\mathbb S^1 \setminus U_s$ for every $s\in (0,1/8)$. Moreover, $\Lambda^s \subset \mathbb S^1\setminus U_s$ and therefore $T$ is uniformly hyperbolic when restricted to $\Lambda_s$. Therefore $(T, \phi, \Lambda^s)$ satisfies Hypothesis~\ref{hyp:H1} for any H\"older weight function $\phi$ for any $s\in(0,1/8)$ and $E = V = \mathbb S^1.$}

\textcolor{black}{We show that for suitable $\phi:\Lambda\to \mathbb R$ and for $\mathrm{Leb}$-almost every $s\in (0,1/8)$, $(T,\phi,\Lambda^s)$ satisfies Hypothesis \ref{hyp:H2}. It is well-known (see e.g.~\cite{Aaronson1997, ZweimullerNotes, Bonanno2018}) that $T$ admits an infinite invariant $\mu\ll \mathrm{Leb}$ on $\mathbb S^1$, and that $T$ is $\mu$-conservative, i.e. for any measurable set $\mu(A)>0$, $\mu$-a.e. $y\in \mathbb S^1$ returns to $A$ so that $\# \{n\in \mathbb N; T^n(y)\in A\} = \infty.$ In particular, the set  
\begin{align}
    B := \bigcap_{q \in (0,1/8)\cap \mathbb Q}\{x\in \mathbb S^1;\text{ there exists }n \in \mathbb N\ \text{such that }T^n(x) \in U_q\} \label{eq:B}
\end{align} 
has full $\mu$-measure, i.e. $\mu(\mathbb S^1 \setminus B) = 0$. Define the set of parameters
$$S := \{s\in (0,1/8);\ \text{there exists }n_1,n_2\in\mathbb N,\ \text{such that } T^{n_1}(s), T^{n_2}(1-s)\in U_s\}.$$
From \eqref{eq:B} it follows that $\Leb(S) = 1/8.$ }

\textcolor{black}{We therefore obtain that for each $s\in \mathbb S$, $\Lambda^s\cap \overline{U}_s = \emptyset$, yielding Hypothesis \ref{hyp:H2}~\ref{it:hyp2-3}. Since $T:\Lambda_s\to \Lambda_s$ is uniformly hyperbolic, \ref{hyp:H2}~\ref{it:hyp2-1} is also satisfied. The last item to verify is \ref{it:hyp2-2} which, for each $s \in \mathbb S$, already holds generically for $\phi \in \mathcal C^\alpha$, $\alpha > 0$.}

\textcolor{black}{As in the previous examples, we obtain conditioned stochastic stability of equilibrium states on $\Lambda^s$.}


\section{Some direct consequences of Hypothesis\texorpdfstring{~\ref{hyp:H1}}{ H1}}\label{sec:consequences}
This section contains several dynamical and topological results that follow from Hypothesis~\ref{hyp:H1} rather immediately and are exploited later in the paper. \textcolor{black}{We also formally introduce the transfer operators $\calL_\e$ and its dual $\calP_\e$, we show that their iterates are compact, and prove that $\calP_\e$ is strong Feller.}

\begin{lemma}\label{lem:invbranches} Let $T$ satisfy Hypothesis~\ref{hyp:H1}. Consider $\delta>0$ small enough. Then there exists $\e_0 \coloneqq \e_0(\delta)$ and $\sigma_1 \coloneqq \sigma_1(\delta) < 1$ such that for every $x,y \in \Lambda_\delta$ satisfying $T(y) = x$ and for every $0<\e < \e_0$, there exists a $\mathcal C^2$ function $h:[-\e,\e]^m\times C_x \to E$, where $C_x$ is the connected component of $x$ in $\Lambda_\delta$, with the following properties holding for every $\w \in \Omega_\e$:
\begin{enumerate}
    \item the map $z\mapsto h(\w,z)$ is a diffeomorphism onto its image,
    \item $T_\w \circ h(\w,z) = z$ for every $z\in C_x$ and $h(0,x)= y,$
    \item $\mathrm{dist}(h(\w,x_1), h(\w,x_2)) \leq \sigma_1 \dist(x_1,x_2)$ for every $x_1,x_2\in C_x$, and
    \item there exists $K_0 = K_0(\delta)>0$ uniform on $\e \in (0,\e_0)$, $x\in \Lambda_\delta$ and $y\in T^{-1}(x)\cap \Lambda_\delta$,  such that $ \sup\{\|\partial_\w h(\w,z)\|
    ;\ \w\in \Omega_\e,\ z\in C_x\}\leq K_0$.
\end{enumerate}
{All statements in this lemma also hold true replacing $\Lambda_\delta$ by $R^i_\delta$, $1 \leq i \leq k$.}
\end{lemma}
\begin{proof}
Take $\delta_0,\e_0>0$ small enough such that
\begin{equation}\label{eq:expansion}
    \sigma_1 \coloneqq \sup\left\{\|\d T_\w(x)^{-1}\|;\ x\in \Lambda_{\delta_0},\ \w\in\Omega_{\e_0} \right\}<1, 
\end{equation}
and such that the exponential map $\exp_{z}: B_{\delta_0}(0) \subset T_z E\to E$ is well defined for every $z\in \Lambda_{\delta_0}$. 

{Observe that there exists $\delta_2>0$ such that for $\e_0>0$ small enough and for every $\w \in [-\e_0,\e_0]^m$, if $\dist(x_1,x_2)<\delta_2$ then we obtain that $\dist(T_\w(x_1),T(x_2))<\delta_0.$} 
Consider the map
\begin{align*}
   G= G_{x,y}:[-\e,\e]^m\times B_{\delta}(x) \times B_{\delta_2}(y) &\to T_x E\\*
    (\w,z_1,z_2)&\mapsto \exp_{x}^{-1}(z_1) - \exp_{x}^{-1}(T_{\w}(z_2)).
\end{align*}    
 Observe that $G(0,x,y)= 0$. Since $ \partial_{y} G(0,x,y)$ is surjective, by means of the implicit function theorem, there exists a $\mathcal C^2$ function $h:[-\e_0(y),\e_0(y)]^m\times B_{\tau(y)}(x) \to E$ such that $T_{\w} (h(\w,z)) = z$ for every $z\in B_{\tau(y)}(x)$ and $h(0,x)=y$. Notice, as well, that $\e_0(y),\tau(y)$ can be taken uniformly since $\Lambda_\delta$ is compact and therefore we can $\mathcal C^2$-extend $h$ to the domain $[-\e,\e]^m\times C_x.$ Finally, from \eqref{eq:expansion} we obtain that the function $h$ satisfies all the desirable properties.
 Replacing $\Lambda_\delta$ by $R_\delta^i$ follows from Lemma~\ref{lem:dyndec_short} item~\ref{it:dyndecshort2}. 
\end{proof}

\begin{lemma}\label{lem:preims_list}
Let $T$ satisfy Hypothesis~\ref{hyp:H1}.
There exists $\delta_0> 0$ small enough satisfying Lemma~\ref{lem:dyndec_short} such that
\begin{enumerate}[label = (\arabic*), series = LP1]
    \item\label{it:preims1} there exists $\sigma_0\coloneqq \sigma_0(\delta_0) \in (0,1)$ such that $T^{-1}(\Lambda_\delta) \cap \Lambda_\delta \subset \Lambda_{\sigma_0 \delta}$ for all $0<\delta<\delta_0$.
\end{enumerate}
Moreover, there exists $\varepsilon_0 \coloneqq \varepsilon_0(\delta)$ satisfying Lemma~\ref{lem:invbranches} such that for every $0 < \varepsilon < \varepsilon_0$ we have that:
\begin{enumerate}[label = (\arabic*), resume*=LP1]
    \item\label{it:preims2} there exists $\sigma\coloneqq\sigma(\delta,\e) \in (0,1)$, such that $T^{-1}_\w (\Lambda_\delta) \cap \Lambda_\delta \subset \Lambda_{\sigma \delta}$ for every $\w \in \Omega_\e$, and
    \item\label{it:preims3} for all $x,y$ lying in the same connected component of $\Lambda_\delta$ and all $\w\in \Omega_\e$ we have
    \(\#\{T^{-1}(x)\cap \Lambda_\delta\} =  \#\{T_\w^{-1}(y)\cap \Lambda_\delta\}.\)
\end{enumerate}    
{All statements in this lemma also hold true replacing $\Lambda_\delta$ by $R^i_\delta$, $i \in \{1, \dots, k\}$.}
\end{lemma}

\begin{proof}
We prove~\ref{it:preims1}. 
Take $\delta_0$ and $\e_0$ small enough such that $\|\d T_{\w_0}(x)^{-1}|_{\Lambda_{\delta_0}}\| <1$ for {all} $\w_0\in [-\e_0,\e_0]^m$ and $\Lambda_{2\delta_0}\subset V$, where $V$ is {as} in Hypothesis~\ref{hyp:H1}. Let $0<\delta<\delta_0$ and $0<\e<\e_0$.

Given $x\in E$ and $v\in T_x E$ such that $\|v\|_x=1$, let $\gamma_{x,v}:(-\delta_0,\delta_0)\to E$ be a geodesic on $E$ such that $\gamma_{x,v}(0) = x$ and $\gamma_{x,v}'(0)= v\in T_x E$. 
From Hypothesis~\ref{hyp:H1}, the fact that $T$ is a $\mathcal C^2$ function and $\Lambda$ is compact, we have that $\gamma_{y,w}(\delta)$ is well-defined for every $y \in \Lambda$ and $w\in T_y E,$ and that there exists $r_0>0$ such that
\begin{equation}\label{eq:Taylor}
   \mathrm{dist}(T\circ \gamma_{y,w}(t), T(y)) \geq  |t|(1+r_0)\ \text{ for every $|t|\leq\delta_0$}.
\end{equation}
From the above equation and the fact that $\d T(x)$ is a surjective linear operator, we obtain that $T(B_{\delta/(1+r_0)}(y)) \supset B_\delta(T(y))$. Take $y\in \Lambda_\delta$, {then} there exists $x\in \Lambda$, $v\in T_x E$ and $h\in [-\delta,\delta]$ such that $y = \gamma_{x,v}(h)$. Let $x_1,\ldots,x_{\ell} \in \Lambda$ be all pre-images of $x$. From~\eqref{eq:Taylor} there exist $h_1,\ldots,h_\ell \in [-\delta/(1+r_0),\delta/(1+r_0)]$, and unit vectors $v_1 \in T_{x_1}E, \ldots , v_{\ell} \in T_{x_\ell} E$ such that all $y_i = \gamma_{x_i,v_i}(h_i), 1 \leq i \leq \ell,$ are pre-images of $y$. Note that $y_i \in \Lambda_{\delta/(1+r_0)}$. We claim that these are precisely the only pre-images of $y$ in $\Lambda_\delta$. Suppose there exists $y' \in \Lambda_\delta\setminus\Lambda_{\delta/(1+r_0)}$ such that $T(y') = y$. Since $x \in B_\delta(y)$, from~\eqref{eq:Taylor} there is $h'\in [-\delta/(1+r_0),\delta/(1+r_0)]$ and $v'\in T_x E$ such that $T(\gamma_{y',v'}(h')) = x$. This contradicts Hypothesis~\ref{hyp:H1} as $\gamma_{y,v}(h') \in \Lambda_{2\delta}\setminus \Lambda \subset V\setminus \Lambda$. Therefore, $T^{-1}(\Lambda_{\delta})\cap \Lambda_\delta \subset \Lambda_{\delta/(1+r_0)}.$ Set $\sigma_0 \coloneqq 1/(1 + r_0)$.

We prove~\ref{it:preims2}.
For every $y\in \Lambda_\delta$, let $\{y_1,\ldots,y_\ell\}\coloneqq T^{-1}(y)\cap \Lambda_\delta$. Let $h_1,\ldots,h_\ell :[-\e,\e]^m\times C_y\to E$ be the inverse branch functions defined in Lemma~\ref{lem:invbranches}, such that $h_i(0,y)=y_i.$ Since $\mathrm{dist}(h_i(\w,y), y_i) =\mathrm{dist}(h_i(\w,y), h_i(0,y)) < K_0\|\w\| \leq K_0 \e$, and $y_i\in \Lambda_{\delta/(1+r_0)}$ from item~\ref{it:preims1}, we obtain that $h_i(\w,y) \in \Lambda_{\delta(K_0\e + 1/(1 +r_0))}$. Choosing $\e$ small enough, there exists $\sigma \in (0,1)$ for which $h_i(\w,y) \in \Lambda_{\sigma \delta}$ for all $\w \in \Omega_\e$.

To finish the proof, we show that for $\e>0$ small enough, 
\(\#\{T^{-1}(x)\cap \Lambda_\delta\} = \#\{T_\w^{-1}(x)\cap \Lambda_\delta\},\)
for every $x\in \Lambda_\delta$ and $\w\in \Omega_\e.$
From the construction above, we obtain that $\#\{T^{-1}(x)\cap \Lambda_\delta\} \leq \#\{T_\w^{-1}(x)\cap \Lambda_\delta\}.$
Suppose for a contradiction that there exist sequences $\{x_n\}_{n\in\mathbb N}\subset \Lambda_\delta$ and $\{\w_n\}_{n\in\mathbb N}\subset [-\e_0,\e_0]^m$, such that  $\#\{T^{-1}(x_n)\cap \Lambda_\delta\} <\#\{T_{\w_n}^{-1}(x_n)\cap \Lambda_\delta\}$ and $\w_n\to 0.$ 
From the compactness of $\Lambda_\delta$ and the pigeonhole principle, the above assumption implies that there exist sequences $\{y^1_n\}_{n\in\mathbb N}$ and $\{y^2_n\}_{n\in\mathbb N}$ such that: (a) $y^1_n \neq y^2_n$  and $T_{\w_n}(y_n^1) = T_{\w_n}(y_n^2)$ for every $n\in\mathbb N$; and (b) $y^1_n,y^2_n\xrightarrow[]{n\to\infty} y^* \in \Lambda_\delta$.
From the continuity of $(\w,x) \mapsto T_\w (x)$, we obtain that $\dist(T(y_n^1),T(y_n^2))\xrightarrow[]{n\to\infty}0$, which contradicts the fact that $\d T(y^*)$ is invertible and completes the proof.

We prove~\ref{it:preims3}.
From the last part in the proof of item~\ref{it:preims2} we obtain that  $\#\{T^{-1}(x) \cap \Lambda_\delta\} = \#\{T_\w^{-1}(x) \cap \Lambda_\delta\},$ for every $x\in \Lambda_\delta$ and $\w\in \Omega_\e$ for $\e$ sufficiently small. Therefore, it is sufficient to show that the map \( x\in \Lambda_\delta \mapsto \#\{T^{-1}(x) \cap \Lambda_\delta\}\)
is locally constant. This is a direct consequence of $\|\d T^{-1}(x)\|<1$ for {all} $x\in \Lambda_\delta$ and the inverse function theorem (see e.g.~the proof of~\cite[Lemma~11.1.4]{Viana2016}). 

{The last statement follows from replacing $\Lambda_\delta$ by $R_\delta^i$ in every argument above, and from item~\ref{it:dyndecshort2} of Lemma~\ref{lem:dyndec_short}}. Note that for~\ref{it:preims1}, we have that $R^i$ is open in $T^{-1}(R^i)$ (see the proof of~\cite[Corollary 11.2.16]{Viana2016}).
\end{proof}



\textcolor{black}{To approximate the equilibrium states in Theorem~\ref{thm:ruele} and establish conditioned stochastic stability, we propose using quasi-ergodic measures, which we construct from the principal eigenfunctions of the following annealed transfer operators. }

\begin{notation}\label{not:operators}
For each $i \in \{1, \ldots, k\}$ and every $\alpha$-H\"older function $\phi:R^i_\delta \to \mathbb R$, we define the annealed Ruelle-Perron-Frobenius operator
\begin{align*}
 \calL_{\e} : L^1 (R^i_\delta,\rho) &\to L^1(R^i_\delta,\rho) \\*  
 f&\mapsto \E_{\e}\left[ \sum_{T_\w(y)=x} \frac{e^{\phi(y)} f(y) \1_{R_\delta^i}(y)}{|\det \d T_\w(y)|}\right]
\end{align*}
and the annealed Koopman operator
\begin{align*}
 \calP_{\e} : L^\infty (R^i_\delta,\rho) &\to L^\infty(R^i_\delta,\rho) \\*  
f &\mapsto e^{\phi(x)}\E_{\e}[ f\circ T_\w (x)\cdot \1_{R^i_\delta}\circ T_\w (x)],
\end{align*}
which are well-posed from Lemmas~\ref{lem:dyndec_short},~\ref{lem:invbranches} and~\ref{lem:preims_list}. Moreover, given $x\in R^i_\delta$ and $n\in\mathbb N$ we refer to the measure $\mathcal P^n_{\e}(x,\cdot)$ as the unique measure on $R_\delta^i$ such that $\calP_{\e}^n(x,A) = \calP_{\e}^n \1_A (x)$ for every measurable subset $A$ of $R_\delta^i.$

Observe that given an $\alpha$-H\"older weight function $\phi:R^i_\delta \to \R$, then $\calL_{\e}^* = \calP_{\e}$. Indeed, for any $f\in L^1(R^i_\delta)$ and $g\in L^\infty(R^i_\delta)$, a change of variable shows that \textcolor{black}{(see~\cite[Equation 1.1.1]{Demers2021})}
\begin{equation} \label{eq:dual} 
\begin{split}
    \int_{R_\delta^i} f(x) \calP_{\e} g(x) \rho(\d x) &= \mathbb E_\e \left[\int_{R_\delta^i} e^{\phi(x)} f(x) g\circ T_\w(x) \1_{R^i_\delta}\circ T_\w(x) \rho(\d x)\right]\\ 
    &= \int_{R_\delta^i}\E_\e\left[\sum_{T_\w(y) = x} \frac{e^{\phi(y)}f(y)\1_{R^i_\delta}(y)}{|\det \d T_\w(y)|}\right] g(x)\rho(\d x)\\ 
    &=  \int_{R_\delta^i} \calL_{\e} f(x) g(x)  \rho(\d x),
\end{split}
\end{equation}
where (if needed) we assume that $f:E \to \R$ vanishes outside of $R^{i}_\delta$. 

\begin{remark}\label{rmk:correction}
    Not{e} that the Ruelle-Perron-Frobenius operator $\mathcal L$ introduced in Theorem~\ref{thm:ruele} differs from the operator $\mathcal L_\e$ since the latter is divided by $|\det \d T_\w|$. This {causes} the correction $- \log| \det \d T|$ for the limiting potential in Theorems~\ref{thm:local} and~\ref{thm:global}. This choice provides a more interpretable expression for $\mathcal P_\e$ and its eigenfunctions as quasi-stationary measures for the $e^\phi$-weighted Markov process $X_n^{\e,\phi}$.
\end{remark}
\end{notation}

The following proposition {establishes} that $X_n^{\varepsilon,\phi}$ is a strong Feller absorbing Markov process. This is constantly exploited {throughout} the paper.

\begin{proposition}\label{prop:kernel1}
    For every $\alpha$-H\"older function $\phi:R^i_\delta\to \mathbb R$, the operator $\calP_\e: L^\infty(R_\delta^i,\rho) \to  L^\infty(R_\delta^i,\rho)$  is strong Feller, i.e.~given a bounded measurable function $f:M\to \mathbb R$ we have $\mathcal P_\e f \in \mathcal C^0(M).$ In particular{,} $\mathcal P_\e^2$ is a compact operator.
\end{proposition}

\begin{proof}

Let $\{x_n\}_{n\in\mathbb N}\subset R^i_\delta$ be a sequence converging to $x \in R^i_\delta$. Write $[-\e,\e]^m = [-\e,\e]^{e}\times [-\e,\e]^{m-e}$, where $e = \dim E$, and let $F \coloneqq F_\e:[-\e,\e]^{e} \times [-\e,\e]^{m-{e}} \times R_\delta^i \to E$. \textcolor{black}{By the means of the rank theorem~\cite[Theorem 4.12]{Lee2013}} we can assume without loss of generality that $ {\partial}_{\w_0}F( \w_0,\w_1,x_n)$ is surjective for every $n\in\mathbb N$. Let $F^{-1}_{(\w_1, x)}$ denote the inverse of $F$ for fixed $(\omega_1, x)$. Then, for any bounded and measurable function $f:M\to\R$, we obtain that

\begin{align*}
  \mathcal P_\e f(x_n) &= \frac{e^{\phi(x_n)}}{(2\e)^m}\int_{[-\e, \e]^{m-e} \times [-\e, \e]^{e}}   (\1_{R^i_\delta} f)\circ F(\omega_0,\omega_1,x_n) \d\omega_0 \d\w_1\\
    &= \frac{e^{\phi(x_n)}}{(2\e)^{m}} \int_{[-\e,\e]^{m-{e}}} \int_{F([-\e,\e]^{e},\w_1,x_n) \cap R_{\delta}^i}  f(y) \left|\det \d F^{-1}_{(\w_1,x)}(y)\right| \rho(\d y) \d \w_1\\
    &=  \int_{R_\delta^i} \left[\frac{e^{\phi(x_n)}}{(2\e)^{m}} \int_{[-\e,\e]^{m-{e}}}\1_{F([-\e,\e]^{e},\w_1,x_n)}(y) \left|\det \d F^{-1}_{(\w_1,x_n)}(y)\right| \d \w_1 \right] f(y)  \rho(\d y).
\end{align*}
Defining $\kappa$ as
$$ \kappa(x_n,y) \coloneqq \frac{e^{\phi(x_n)}}{(2\e)^{m}} \int_{[-\e,\e]^{m-{e}}}\1_{F([-\e,\e]^{e},\w_1,x_n)}(y) \left|\det \d F^{-1}_{(\w_1,x_n)}(y)\right| \d \w_1 ,$$
it is clear that $ \kappa(x_n,y) \xrightarrow[]{n\to\infty} \kappa(x,y) \ \text{for }\rho\text{-a.e. $y\in R^i_\delta$}.$ 
\textcolor{black}{Since $\kappa$ is bounded on \(R^i_\delta\times R^i_\delta\), then for any
$f\in L^\infty(R^i_\delta,\rho)\subset L^1(R^i_\delta,\rho)$, we have that $\left|\kappa(x_n,\cdot)\,f\right|\leq \|\kappa\|_{\infty}\,|f|
\in L^1(R^i_\delta,\rho).$} Therefore, by the Lebesgue dominated convergence theorem, we have that $\lim_{n\to\infty} \mathcal P_\e f(x_n)=\mathcal P_\e f(x),$ so $\mathcal P_\e h(x)$ is continuous and thus $\calP_\e$ is strong Feller. From {\cite[Chapter 1, Theorem 5.11]{Revuz1984} (which we recall in Lemma~\ref{lem:SF})}, we have that $\mathcal P_\e^2$ is a compact operator.
\end{proof}

From Proposition~\ref{prop:kernel1} and {equation~\eqref{eq:dual}}, we obtain that $\calL_\e^2:L^1(R_\delta^i)\to L^1(R_\delta^i)$ is also a compact operator.

\section{The local problem}\label{sec:local}
In this section, we focus on a {single} repeller $R^i$ of $T$ {from} the dynamical decomposition of Lemma~\ref{lem:dyndec} and establish the stochastic stability of equilibrium states associated with the restricted transformation $\left.T\right|_{R^i}$. To achieve this, we condition the process $X_n^\e$ upon remaining
within a $\delta$-neighbourhood of the repeller $R^i$. For a given $\alpha$-H\"older weight function $\phi$, we begin by showing that there exists a unique quasi-stationary measure $\mu_\e$ for the $e^\phi$-weighted Markov process $X_n^{\e,\phi}$ on $R^i_\delta$ absorbed in $\partial \coloneqq M_\delta \setminus R^i_\delta$. To do so, we {adapt} the analysis of conditionally invariant probability measures provided by Pianigiani and Yorke~\cite{Pianigiani1979} as fixed points of the (normalised) Ruelle-Perron-Frobenius operator, $\widehat{\calL_\e}$. 
We continue with a detailed study of the operator $\calP_\e$ to obtain the (unique) eigenfunction $g_\e$ of maximal eigenvalue $\lambda_\e$.
Finally, we prove the existence and uniqueness of a quasi-ergodic measure of the $e^\phi$-weighted Markov process $X_n^\e$ conditioned upon not escaping the support of $g_\e$, and characterise its limiting behaviour as the noise strength $\e$ vanishes. This measure follows from the pointwise product of $\mu_\e$ and $g_\e$. As previously mentioned, we show that the limiting object as $\e \to 0$ corresponds to an ergodic invariant measure sitting on the repelling set $R^i$ that corresponds to the unique equilibrium state for the potential $\phi-\log |\det \d T |$.

Throughout this section, we assume Hypothesis~\ref{hyp:H1} holds true and employ the notation introduced in Section~\ref{sec:setup_notation}. In particular, we use ``$\e$ small enough'' and ``$\delta$ small enough'' to refer to $\e$ and $\delta$ as in Lemmas~\ref{lem:dyndec_short},~\ref{lem:invbranches} and~\ref{lem:preims_list}. 
All arguments in this section hold for each $1 \leq i\leq k$ and every  $\alpha$-H\"older weight function $\phi :R^i_\delta\to \R$, which we fix once and for all. {To improve readability we drop the super-index $\phi$ of the weighted Markov process $X_n^{\e,\phi}$ and simply write $X_n^{\e}$.}

\subsection{Quasi-stationary measures on \texorpdfstring{$R_\delta^i$}{R i \_δ}}
Denote by $\widehat{\calL_\e}$ the $L^{1}$-normalised operator $\calL_\e$, i.e.~
\[\widehat{\calL_\e} f = \frac{\calL_\e f}{\|{\calL_\e f}\|_{1}}.\]

\begin{notation}
Given a compact metric space $(N,d)$ and $0<\alpha<1$ we denote by $\mathcal C^{\alpha}(N)$ the set of $\alpha$-H\"older functions $f:N \to \R$ and consider the $\alpha$-H\"older norm
\[\|f\|_{\mathcal C^\alpha} = \sup_{x\in N} |f(x)| + \sup_{x\neq y}\frac{|f(x) -f(y)|}{d(x,y)^{\alpha}}.\] 

\end{notation}
To obtain a quasi-stationary density for the conditioned process on each component $R^{i}_\delta$ we apply the Schauder-Tychonoff fixed point theorem (see, e.g.~\cite[Theorem 2.2.3]{Viana2016}) to the operator $\widehat{\calL_\e}$ acting on a suitable space $C_\beta$.

\begin{theorem}\label{thm:pianigiani}
    Consider an $\alpha$-H\"older weight function $\phi:R^i_\delta\to \R$ and suppose that $T$ satisfies Hypothesis~\ref{hyp:H1} \textcolor{black}{on $R^i$}. Let $\delta>0$ be small enough. Then, for every $\e >0$ small enough there exists a measure $\mu_\e(\dx)$ on $R^i_\delta$ such that:
        \begin{enumerate}[label = (\arabic*)]
        \item\label{it:piani1} $\mu_\e$ is \textcolor{black}{a} quasi-stationary measure of the $e^\phi$-weighted Markov process $X_n^\e$ on $R^i_\delta$ \textcolor{black}{with growth rate given by
        $\lambda_\e = r(\mathcal L_\e:L^1(R_\delta^i,\rho)\to L^1(R^i_\delta,\rho))),$}
        \item\label{it:piani2} $\mu_\e$ is absolutely continuous with respect to $\rho$, and
        \item\label{it:piani3} defining $m_\e \coloneqq \mu_\e(\d x)/\rho(\d x)$, there exists $C>0$ such that $\|m_\e\|_{\mathcal C^\alpha}\leq C$ and $m_\e(x)>0$ for every $x\in R_\delta^i.$
    \end{enumerate}
    \label{thm:me}
\end{theorem}
\begin{proof}
Given $\beta>0$ consider the set
\[C_\beta \coloneqq 
\left\{ f\in L^1(R^i_\delta,\rho)\ \middle\vert \begin{array}{l}
    \int f \, \d \rho =1,\ f>0, \text{ and }\frac{f(x)}{f(y)}\leq e^{\beta d(x,y)^\alpha}\  \text{if}\ x,y\\
    \text{lie in the same connected component of }R^i_\delta
  \end{array}\right\}. \]

We divide the proof into 3 steps, \textcolor{black}{first establishing the existence of the density $m_\e \in C_\beta$.}

\begin{step}[1]\label{step:contraction}
{T}here exists $\beta>0$ such that $\widehat {\calL_\e} (C_\beta)\subset C_\beta. $
\end{step}
\begin{stepproof}[Proof of Step~\ref{step:contraction}]
First of all, observe that if $\delta>0$ is small enough $f>0, f \in C_\beta,$ implies $\mathcal L_\e f>0$ for every $\e >0.$ Define  $\psi: R_\delta^i\to \mathbb R$ as $\psi \coloneqq \phi - \log|\det\, \d T|$ and let
\[ D \coloneqq \sup_{x\neq y}\frac{ \left|\psi(x)-\psi(y)\right|}{d(x,y)^\alpha}<\infty,\]
Recall from Lemma~\ref{lem:preims_list} (3) that if $\e>0$ is small enough, then given $x,y$ in the same connected component $C$ of $R^i_\delta$ we have
\[\#\{T^{-1}(x)\cap \Lambda_\delta\} =\#\{T_\w^{-1}(y)\cap \Lambda_\delta\}\]
for every $\w\in \Omega_\e.$ Suppose that $\#\{T^{-1}(x)\cap \Lambda_\delta\} = \ell.$ Let $h_1,\ldots,h_\ell:[-\e,\e]^m\times C \to R_\delta^i$ be the pre-image functions (inverse branches) defined in Lemma~\ref{lem:invbranches}. Given $f\in C_\beta$, we have that
\begin{align*}
\mathcal L_\e f(x) &= \mathbb E_\e\left[\sum_{T_\w (z) = x} \frac{e^{\phi(z)} f(z)}{|\mathrm{det}\, \d T_\w (z)|}\right]= \mathbb E_\e\left[\sum_{j=1}^{\ell} e^{\psi\circ h_j(\w,x)} f\circ h_j(\w,x)\right] \\
&= \mathbb E_\e\left[\sum_{j=1}^{\ell} e^{\psi\circ h_j(\w,x) - \psi\circ h_j(\w,y)} \frac{f\circ h_j(\w,x)}{f\circ h_j(\w,y)} f\circ h_j(\w,y) e^{\psi\circ h_j(\w,y)}    \right] \\
&\leq \left(\sup_{i\in\{1,\ldots, \ell\}}e^{ (\beta+ D) \dist (h_j(\w,x),h_j(\w,y))^\alpha }\right)\mathcal L_\e f(y)\\
&\leq e^{ \sigma^\alpha (\beta + D) \dist (x,y)^\alpha} \mathcal L_\e f(y),
\end{align*}
with $\sigma$ from Lemma~\ref{lem:preims_list}. Therefore, if $f\in C_\beta$ then $ \widehat{\calL_\e} f \in C_{\sigma^\alpha (\beta +D)}.$ Taking $\beta > D \sigma^\alpha/(1-\sigma^\alpha)>0$ we {conclude} Step~\ref{step:contraction}.
\end{stepproof}
\begin{step}[2]\label{step:qsm}
For every $\e>0$ small, there exists $m_\e \in C_\beta$ such that $\widehat{\calL_\e}m_\e = m_\e.$
\end{step}
\begin{stepproof}[Proof of Step~\ref{step:qsm}]
Observe that $C_\beta$ is pre-compact and convex in $L^1(R_\delta^i,\rho)$. From the Schauder fixed-point theorem, there exists $m_\e$ lying in the closure of $C_\beta$ such that $\widehat{\calL_\e}m_\e = m_\e$, which implies that $\mathcal L_\e m_\e = \lambda_\e m_\e$ for $\lambda_\e = \|\mathcal L_\e m_\e\|>0.$ We claim that $m_\e \in  C_\beta.$ Suppose for a contradiction that $m_\e \in \overline{ C_\beta}^{{}_{L^1(R_\delta^i,\rho)}}\setminus C_\beta.$ \textcolor{black}{Since $m_\e \in \overline{C_\beta}^{{}_{L^1(R_\delta^i,\rho)}}$, there exists a sequence $m^{(n)}_\e \in C_\beta, n \in \mathbb N,$ such that $m^{(n)}_\e \xrightarrow[]{n\to\infty}m_\e\ \text{in }L^1(R_\delta^i,\rho).$ Moreover, from~\cite[Corollary 4.10]{Salamon2016}, there exists a subsequence $\{m^{(n_k)}_\e\}_{k\in\mathbb N}\subset \{m^{(n)}_\e\}_{n\in\mathbb N}$ such that $m^{(n_k)}_\e \to m_\e$ pointwise $\rho$-almost surely. Therefore, for $\rho$-almost every $x,y$ lying in the same connected component of $R_\delta^i$ 
\begin{align}
    m_\e(x) = \lim_{k\to\infty} m_\e^{(n_k)} (x) \leq \lim_{k\to\infty }e^{\beta d(x,y)} m_\e^{(n_k)}(y) \leq  e^{\beta d(x,y)} m_\e(y).\label{eq:me}
\end{align}
The above equation implies that we may choose a representative of $m_\e$ in $L^1(R_\delta,\rho)$ which is continuous and therefore we obtain that the above equation holds for each $x,y$ lying in the same connected component of $R_\delta^i$. From \eqref{eq:me} and $m_\e\not\in C_\beta$ it must exist $x\in R^i_\delta$ such that $m_\e(x) = 0.$} \textcolor{black}{It follows that} $m_\e(y) = 0$ for every $y$ in the same connected component $C_x$ of $x$ in $R^i_\delta$, \textcolor{black}{as $x$ and $y$ may be interchanged in \eqref{eq:me}}. Hence, for every $y\in C_x$ and $n >0$,
\[ 0 = m_\e(y) =\frac{1}{\lambda_\e^n}\mathbb E_\e\left[\sum_{T^n_\w (z) = y} \frac{e^{S_n\phi(\w,z)} \1_{R_\delta^i}(z)m_\e(z)}{|\mathrm{det}\, \d T^n_\w (z)|}\right],\]
where $S_n \phi (\w,z) = \sum_{i=0}^{n-1} \phi\circ T^i_\w (z).$

This implies that $m_\e$ vanishes in the connected components of points in $T^{-n}(y)\cap R_\delta^i,$ for every $y\in C_x$. Since there exists $z\in R^i$ such that $\{T^n(z)\}_{n\in\mathbb N}$ is dense in $R^i$, it follows that $m_\e \equiv 0$, which is a contradiction.
\end{stepproof}

\begin{step}[3]\label{step:LFeller} \textcolor{black}{We prove that if $f\in L^1(R_\delta^i, \rho)$ then $ \mathcal L_\e f \in \mathcal C^0(R_\delta^i)$.}
\end{step}
\begin{stepproof}[\textcolor{black}{Proof of Step~\ref{step:LFeller}}]
\textcolor{black}{
Let $x\in R_\delta^i$ and $r>0$ small enough. As in Step~\ref{step:contraction}, let the function
$$h_1,\ldots,h_\ell : [-\e,\e]^m \times B_r(x)\cap R_\delta^i \to R_\delta^i$$
be such that for each $j\in\{1,\ldots, k\}$
\begin{itemize}
    \item $T_{\omega} \circ h_{j}(\omega,z) = z$, for every $(\omega_0,z)\in [-\e,\e]^m\times B_r(x)$ for each $i\in\{1,\ldots, k\}$;
    \item for any $z\in B_r(x)$, $T_\omega^{-1}(z) = \{h_1(\omega,x),\ldots, h_\ell(\omega,x)\}.$
\end{itemize}
Moreover, since 
$$ 0 = \partial_\omega(z) = \partial_\omega (T_\omega \circ h_j(\omega,z)) = \partial_\omega T_\omega (h_j(\omega,x)) + \d  T_\omega(h_j(\omega,x)) \partial_\omega h_j(\omega,x),$$
we have that
$$\partial_\omega h_j(\omega,x) = -[\d  T_\omega(h_j(\omega,x))]^{-1}\partial_\omega T_\omega (h^{(i)}(\omega,x)),$$
which is well defined due to \eqref{eq:expansion}. Since $\partial_\omega T_\omega$ is full rank, we obtain that $\partial_\omega h_j$ is also full rank.}

\textcolor{black}{Given $f\in L^1(R_\delta^i,\rho)$ we have that defining $\psi(\omega,x) = \phi(x) - \log |\det \d T_\omega(z)|$
\begin{align*}
   \mathcal L_\e f(x) &= \mathbb E_\e\left[\sum_{T_\w (z) = x} \frac{e^{\phi(z)} f(z)}{|\mathrm{det}\, \d T_\w (z)|}\right]= \mathbb E_\e\left[\sum_{j=1}^{\ell} e^{\psi(\omega, h_j(\w,x))} f\circ h_j(\w,x)\right]\\
   &=\sum_{i=1}^\ell \frac{1}{(2\e)^m}\int_{[-\e,\e]^m}  e^{\psi(\omega, h_j(\w,x))} f\circ h_j(\w,x) \d \omega,
\end{align*}
Repeating verbatim the computations in the proof of Proposition~\ref{prop:kernel1}, noting that this argument uses only the facts that the map $F$ (which in the
present setting is \(h_j\)) has full rank in the \(\omega\)-variable and that $f\in L^1(R_\delta^i,\rho)$, we obtain the desired result.}
\end{stepproof}

\textcolor{black}{We may now conclude the proof of the theorem.}
    

\textcolor{black}{Items~\ref{it:piani2} and~\ref{it:piani3} follow directly from Step~\ref{step:qsm}. To prove item~\ref{it:piani1}, define $\mu_\e(\d x) = m_\e (x) \d x$, then for any $f:R^i_\delta\to \R$ bounded and measurable
\begin{align*}
\int_{R^i_\delta} \mathcal P_\e f(x) \mu_\e (\d x) &=   \int_{R^i_\delta} \mathcal P_\e f(x) m_\e (x) \d x \\&= \|\mathcal L_\e m_\e\|_1 \int_{R^i_\delta} f(x) \widehat{\mathcal L_\e} m_\e (x) \d x  = \lambda_\e \int_{R^i_\delta} f(x) m_\e(x) \dx,
\end{align*}
where $\lambda_\e = \|\mathcal L_\e m_\e\|_1.$}  

\textcolor{black}{From Step~\ref{step:LFeller} we have that $\mathcal L_\e (L^1(R^i_\delta, \rho) \subset \mathcal C^0(R^i_\delta)$. Moreover $\mathcal L_\e : L^1(R^i_\delta, \rho) \to \mathcal C^0(R^i_\delta)$ is a positive linear operator between Banach spaces, therefore $\mathcal L_\e$ is bounded (see~\cite[Theorem II.5.3]{Schaefer76}, in particular there exists $K>0$
\begin{align}
    \sup_{\|f\|_1 =1}\|\mathcal L_\e f\|_\infty \leq K. \label{eq:k1}
\end{align}}
\textcolor{black}{From Step~\ref{step:qsm} we have that $\inf m_\e (x) = c >0$. Observe that for any $f\in L^1(R^i_\delta, \rho),$ then for every $x\in R_\delta^i$
\begin{align}
    \mathcal L_\e f(x)\leq \frac{\|\mathcal L_\e f\|_\infty}{c} m_\e(x).\label{eq:k2}
\end{align}
Therefore, using \eqref{eq:k1} and \eqref{eq:k2} 
\begin{align*}
 r(\mathcal L_\e) &= \lim_{n\to\infty} \sup_{\|f\|_1 = 1} \|\mathcal L^n_\e f\|_1^{1/n}  = \lim_{n\to\infty} \sup_{\|f\|_1 = 1, f\geq 0} \|\mathcal L^n_\e f\|_1^{1/n} \\ 
 &=\lim_{n\to\infty} \sup_{\|f\|_1 = 1,f\geq 0} \left|\int_{R_\delta^i}\mathcal L^n_\e f \, \rho (\d x)\right|^{1/n} \\
 &\leq \lim_{n\to\infty} \sup_{\|f\|_1 = 1,f\geq 0} \left|\int_{R_\delta^i}   \frac{\|\mathcal L_\e f\|_\infty}{c } \mathcal L^{n-1}_\e m_\e(x) \,  \d x\right|^{1/n} \\
 &= \lim_{n\to\infty} \sup_{\|f\|_1 = 1,f\geq 0} \lambda_\e^{\frac{n-1}{n}}  \left|\int_{R_\delta^i}   \frac{\|\mathcal L_\e f\|_\infty}{c } m_\e(x) \,  \d x\right|^{1/n}\\
 &\leq \lim_{n\to\infty} \sup_{\|f\|_1 = 1,f\geq 0} \lambda_\e^{\frac{n-1}{n}}  \left|\int_{R_\delta^i}   \frac{K}{c} m_\e(x) \,  \d x\right|^{1/n} = \lambda_\e.
\end{align*}
Since $\lambda_\e \in \sigma(\mathcal L^1\to L^1)$ we have that $\lambda_\e = r(\mathcal L_\e)$, which implies item~\ref{it:piani1}.}
\end{proof}

For {the remainder of this} section, let $m_\e \in C_\beta$ denote the unique function such that $\calL_\e m_\e = \lambda_\e m_\e$ and let $\mu_\e$ be the unique quasi-stationary measure, i.e.~be such that $m_\e = \mu_\e(\dx)/\rho(\dx)$, as of Theorem~\ref{thm:pianigiani}.

\subsection{Analysis of the operator \texorpdfstring{$\calP_\e:L^\infty(R_\delta^i) \to L^\infty(R_\delta^i)$}{P\_ε}} We now study the adjoint operator of $\mathcal L_\e$ to obtain the eigenfunction $g_\e$, associated with the maximal eigenvalue $\lambda_\e$ from the previous result, and its properties. We also construct the unique quasi-ergodic measure of the $e^\phi$-weighted Markov process $X_n^\e$ on $\{g_\e > 0\}$.

\begin{lemma} \label{lem:positive}
\textcolor{black}{Assume that $T$ satisfies Hypothesis~\ref{hyp:H1} on $R^i$}, let $\delta, \epsilon > 0$ be small enough and let $\lambda_\e$ be the eigenvalue associated with $m_\e$ from Theorem~\ref{thm:pianigiani}. Let $g_\e \in \ker(\calP_\e -\lambda_\e \Id)\cap \mathcal C^0_+(R^i_\delta)$, then $R^i\subset \{g_\e>0\}$. 
\end{lemma}

\begin{proof}
We divide the proof into two steps \textcolor{black}{and note that the fact that $\ker(\calP_\e -\lambda_\e \Id)\cap \mathcal C^0_+(R^i_\delta) \neq \emptyset$ will be shown in Theorem~\ref{thm:span_g_eps}}. First, we check that $g_\e$ is positive on dense orbits of $T$ and, second, construct a neighbourhood of $R^i$ where $g_\e$ is positive. It is clear that a dense orbit exists since $T:R^i \to R^i$ is topologically transitive by Lemma~\ref{lem:dyndec}. Let $K_0 = \min_{x\in R_\delta^i} e^{\phi(x)} >0.$

\begin{step}[1]\label{step:positive1}
If $\{T^j (x_0)\}_{j\in \N}$ is dense in $R^i,$ then $g_\e (x_0) >0.$ \label{step:gex}
\end{step}
\begin{stepproof}[Proof of Step~\ref{step:positive1}]
Recall that $g_\e \in \mathcal C_+^0(R_\delta^i)$. Assume that $g_\e (x_0) =0,$ then for every $n \in \N$,
\[
    0 = g_\e (x_0) = \frac{1}{\lambda_\e^n}\mathcal P_{\e}^n g_\e(x_0) \geq  \frac{K_0^n}{\lambda_\e^n}\E_\e[ g_\e \circ T^n_\w(x_0) \cdot  \1_{R_\delta^i}\circ T^n_\w(x_0)]. 
\] 
Combining this with the submersion theorem applied to $\partial_\w T_\w$ (see~\cite[Theorem 4.12]{Lee2013}) and the fact that $R_\delta^i$ is compact, we obtain that there exists $r_0>0$ such that $\left.g_\e\right|_{B_{r_0}(T^n(x_0))} =0$ for every $n\in\N$. Since $\{T^j(x_0)\}_{j \in \N}$ is dense in $R^i$, there exists a neighbourhood $U\supset R^i$ such that $\left.g_\e\right|_U = 0.$

Let $T_{\cap R_\delta^i}(U) \coloneqq T(U \cap R_{\delta}^i)$. Recall that, from \eqref{eq:Taylor}, there exists $N\in\mathbb N$ such that $T^N_{\cap R_\delta^i}(U) \supset R_\delta^i$.
Take $y\in R_\delta^i$. Then, there exists $z\in U$ such that $T^N(z) = y$ and $T^j(z) \in U$ for every $j\in\{1,\ldots, N\}.$ Since \[0  = g_\e (z) = \frac{1}{\lambda_\e^N} \mathcal P^N_\e g_\e(z) \geq   \frac{K_0^N}{\lambda^N_\e}\E_\e[ g_\e \circ T^N_\w(z) \cdot  \1_{R_\delta^i}\circ T^N_\w(z)],\]
continuity of $g_\e$ and the submersion theorem applied to $\partial_\w T_\w$ yields $g_\e (y) = 0$. This contradicts $g_\e \neq 0$.
\end{stepproof}

\begin{step}[2]\label{step:positive2}
{T}here exists an open set $B\supset R^i$, such that $g_\e(x)>0$ for every $x\in B.$
\end{step}
\begin{finalstepproof}
Set $B \coloneqq \{x\in R_\delta^i;\ \exists\ \w_0 \in (-\e/2,\e/2)^m\ \text{s.t. }T_{\w_0}(x)\in R^i\}.$
From the submersion theorem, we have that given $x\in B$ and $\w_0 \in (-\e/2,\e/2)^m$ such that $T_{\w_0}(x)\in R^i$ we obtain that there exists $r_1>0$ such that
\[\bigcup_{\w\in \Omega_\e} T_\w (x) \supset B_{r_1}(T_{\w_0} (x)) \ \text{for some }r_1>0. \]

Let $x_0 \in R^i$ be such that $\{T^j(x_0)\}_{j\in\mathbb N}$ is dense in $R^i$, then there exists $N_0$ such that $T^{N_0}(x_0) \in B_{r_1}(T_{\w_0} (x)).$ From Step~\ref{step:gex}, $g_\e(T^{N_0}(x_0)) >0.$ Continuity of $g_\e$ then implies
\[0<\frac{K_0}{\lambda_\e} \E_\e [g_\e \circ T_\w (x) \cdot \1_{R_\delta^i} \circ T_\w(x)]\leq \frac{1}{\lambda_\e} \mathcal P_\e g_\e (x)  = g_\e(x).\]
This {concludes} Step~\ref{step:positive2} and {proves} the lemma.
\end{finalstepproof}
\begingroup
\renewcommand{\qedsymbol}{\empty}
\end{proof}
\endgroup

\begin{theorem}\label{thm:span_g_eps}
Consider the operator $\calP_\e:L^\infty(R^i_\delta) \to L^\infty(R^i_\delta)$. Then, $\ker(\calP_\e - \lambda_\e  \Id)= \s(g_\e)$ for some $g_\e \in \mathcal C^0_+(R_\delta^i).$
\end{theorem}
\begin{proof}
Let $g_\e\in \ker (\calP_\e -\lambda_\e)$. Since $\calP_\e$ is strong Feller, then $g_\e\in \mathcal C^0(R_\delta^i,\mathbb C)$. Moreover, since $\calP_\e (\mathcal C^0(R_\delta^i))\subset\mathcal C^0(R_\delta^i),$ it is clear that $\mathrm{Re}(g_\e),\mathrm{Im}(g_\e)\in \ker (\calP_\e - \lambda_\e).$ 
Given $g_\e\in \ker (\calP_\e -\lambda_\e )\cap \mathcal C^0(R_\delta^i)$, we claim that $g^{\pm}_\e \in \ker (\calP_\e-\lambda_\e \Id).$ Indeed, observe that (see~\cite[Propositions 3.1.1 and 3.1.3]{Lasota1994})
\[g_\e^\pm = \left(\frac{1}{\lambda_\e} \calP_\e g_\e\right)^\pm \leq \frac{1}{\lambda_\e} \calP_\e g_\e^\pm,\]
where $g_\e^\pm = \max\{0, \pm g_\e\}. $ Therefore,
\begin{align*}
  0 &= \int_{R^i_\delta}\left( \frac{1}{\lambda_\e }\calP_\e |g_\e|(x) - |g_\e|(x) \right)  \mu_\e(\d x)\\
  &= \int_{R^i_\delta} \left(\frac{1}{\lambda_\e }\calP_\e g_\e^+(x) - g_\e^+(x) \right)\mu_\e(\d x) +  \int_{R^i_\delta}\left( \frac{1}{\lambda_\e }\calP_\e g_\e^-(x) - g_\e^-(x) \right) \mu_\e(\d x).
\end{align*}
From Theorem~\ref{thm:me}, $\supp \mu_\e = R_\delta^i,$ therefore $\calP_\e g_\e^\pm = \lambda_\e g_\e^\pm.$

Take $g_1,g_2 \in \ker (\calP_\e - \lambda_\e)\cap \mathcal C_+^0( R^i_\delta).$ From Lemma~\ref{lem:positive}, we have $g_1, g_2 > 0$ on $R^i$.
Choose $t_0>0$ such that $t_0 =\inf\{t;\, g_1(x) - tg_2(x)<0\text{ for some }x\in R^i\}.$

Since $g_1 - t_0 g_2 \in \ker (\calP_\e - \lambda_\e\Id)$, then $(g_1 - t_0g_2)^+ \in  \ker (\calP_\e - \lambda_\e).$ However, from the choice of $t_0$ and Lemma~\ref{lem:positive} we obtain that $(g_1 - t_0g_2)^+=0$. From the minimality of $t_0$, it follows that $g_1(x) = t_0 g_2 (x)$ for every $x\in R^i.$ {Observe that $(g_1 - t_0g_2)^+=0$ yields that $t_0 g_2 \geq g_1.$} Therefore $t_0 g_2 -g_1 \in \ker (\calP_\e - \lambda_\e)\cap \mathcal C^0_+(R^i_\delta)$ and $\left.\left(t_0 g_2 -g_1\right)\right|_{R^i}=0,$ implying that $t_0 g_2 -g_1 = 0$.
\end{proof}

\begin{theorem}\label{thm:unique-qsm}
   \textcolor{black}{There exists a unique quasi-stationary measure $\mu$ on $R_\delta^i$ for $X^\e_n$  such that $\Lambda \subset \supp \mu.$ In fact, this quasi-stationary measure is the fully-supported measure constructed in Theorem~\ref{thm:pianigiani}.}
\end{theorem}
\begin{proof}
\textcolor{black}{Let $\mu$ be a quasi-stationary measure on $R_\delta^i$ with growth rate $\lambda$. We divide the proof into two cases: $(i)$ $\lambda\neq \lambda_\e$ and $(ii)$ $\lambda = \lambda_\e$.}

\textcolor{black}{To prove $(i)$, observe that $g_\e\in \mathrm{ker}(\mathcal P_\e -\lambda_\e) \cap \mathcal C^0_+(R_i^\delta)$ and $\{g_\e > 0\} \supset \Lambda$, from Lemma\ref{lem:positive} and Theorem~\ref{thm:span_g_eps}. It follows that
$$\int_{R_\delta^i} g_\e (x) \mu(\d x) = \frac{1}{\lambda_\e} \int_{R_\delta^i}  \mathcal P_\e g_\e (x) \mu (\d x) = \frac{\lambda}{\lambda_\e} \int_{R_\delta^i} g_\e (x) \mu(\d x).$$
Since $\lambda\neq \lambda_\e$ we obtain that $\int_{R_\delta^i} g_\e (x) \mu(\d x) =0$ which implies that $\supp \mu \cap \Lambda = \emptyset.$}

\textcolor{black}{To prove $(ii)$, since $\mathcal P_\e(x,\d y) \ll \rho(\d y)$ for every $x\in R_\delta^i$ (recall Lemma~\ref{lem:SF}) we have that $\mu \ll \rho$. Define $$m = \frac{\mu(\d x)}{\rho(\d x)}\in L^1(R_\delta^i ,\rho).$$ 
Since $\mu(\d x) = m (x) \d x$ we have that $\mathcal L_\e m = \lambda_\e m$, however $\lambda_\e$ is a simple eigenvalue of $\mathcal L_\e$, this follows from the fact that $\mathcal L_\e^* = \mathcal P_\e$, $\mathcal P_\e$ is a quasi-compact operator (since $\mathcal P_\e^2$ is compact, see Lemma~\ref{lem:SF} in the Appendix) and $\lambda_\e = r(\mathcal P_\e)$ is a simple eigenvalue (see Theorem~\ref{thm:span_g_eps}).}
\end{proof}

For {the remainder of this} section, let $g_\e \in \mathcal C^0_+(R^i_\delta)$ denote the unique function such that $\calP_\e g_\e = \lambda_\e g_\e$ and normalise{d} so that $\int g_\e \d\mu_\e = 1$, as of Theorems~\ref{thm:span_g_eps} and~\ref{thm:unique-qsm}. We summarise some relevant properties of $\calP_\e :L^\infty(R_\delta^i) \to L^\infty(R^i_\delta)$ that have been shown above:
\begin{enumerate}
    \item $\calP_\e: L^\infty(R_\delta^i,\rho)\to L^\infty(R_\delta^i,\rho)$ is a strong Feller operator,
    \item $\dim\ker (\calP_\e - \lambda_\e) =1$, where $\lambda_\e = r(\calP_\e)$ is the spectral radius,
    \item there exists $\mu_\e \in \mathcal M_+(R_\delta^i)$ and $g_\e\in \mathcal C^0_+(R_\delta^i),$ such that $\calP^*_\e \mu_\e = \lambda_\e \mu_\e$ and $\calP_\e g_\e = \lambda_\e g_\e$, and
    \item $\mu_\e \ll \rho$ and $\supp \mu_\e = R_\delta^i$.
    \item \textcolor{black}{measure $\mu_\e$ is the unique quasi-stationary measure of $X_n^\e$ on $R_\delta^i$ such that $\Lambda\subset \supp \mu_\e$.}
\end{enumerate}
In particular, this implies that $\calP_\e$ satisfies Hypothesis~\ref{hyp:A} in Appendix~\ref{sec:appendix}. The lemma below is a consequence of the properties just listed and Theorems~\ref{thm:QED} and~\ref{thm:QED2}, whose proof is deferred to the Appendix in order not to break the flow of the text.

\begin{lemma}\label{lem:uniqueqed}
The measure $\nu_\e(\d x) \coloneqq g_\e(x) \mu_\e(\d x)$ is the unique quasi-ergodic measure of the $e^\phi$-weighted Markov process $X_n^\e$ on $\{g_\e >0\}$. If we further assume that $T:R^i\to R^i$ is topologically mixing, then $\nu_\e$ is also a quasi-ergodic measure of the $e^\phi$-weighted Markov process $X_n^\e$ on $R^{i}_\delta$.
\end{lemma}

\begin{proof}
It is clear from the properties of $\calP_\e$ listed above that it satisfies Hypothesis~\ref{hyp:A} (see Appendix~\ref{sec:appendix}). Hence, Theorem~\ref{thm:QED} implies that $\nu_\e$ is the unique $e^\phi$-weighted quasi-ergodic measure for $X_n^\e$ on $\{g_\e >0\}$.

To finish the proof of the theorem, it remains to be shown that if $T$ is topologically mixing, then $\nu_\e$ is a $e^\phi$-weighted quasi-ergodic measure for $X_n^\e$ on $R^i_\delta$. Since $T$ is topologically mixing, then $X_n^\e$ is aperiodic in $R^i_\delta$ and $\{g_\e >0\}$. Let $k_\e \coloneqq \#(\sigma_{\mathrm{per}}(\frac{1}{\lambda_\e}\mathcal P_\e)\cap \mathbb S^1 )$. From Lemma~\ref{lemma:point_periph_spec}, $k_\e <\infty$. Moreover, from Proposition~\ref{prop:nerdola} and Lemma~\ref{lem:ciclying} we obtain that there exist sets $C_i \subset \{g_\e >0\}$, $i\in \{0,1,\ldots,k_\e-1\}$ such that $C_0 \sqcup C_1 \sqcup \ldots \sqcup C_{k_\e-1} = \{g_\e >0\},$
and $\{\mathcal P_\e\1_{C_i} >0 \}\subset C_{i-1\ (\mathrm{mod}\ k_\e)},$ for every $i\in\{0,1,\ldots,k_\e -1\}$. Since $T$ is assumed to be topologically mixing on $R^i$, $X_n^\e$ is aperiodic and thus $k_\e=1$. Finally, from Theorem~\ref{thm:QED2} we obtain that $\nu_\e$ is a quasi-ergodic measure of the $e^\phi$-weighted Markov process $X_n^\e$ on $R^i_\delta$.
\end{proof}

\subsection{The limit \texorpdfstring{$\e \to 0$}{ε->0}. Proof of the main (local) result} We conclude this section with the main results concerning the stochastic stability of equilibrium states on each repeller $R^i_\delta$ and their limiting behaviour as $\e \to 0$.

\begin{notation} 
    Recall that given a suitable function $f$ we denote the action of the (deterministic) Ruelle-Perron-Frobenius operator $\calL$ for the potential $\psi = \phi -\log|\det \d T|$ {\cite[Chapter 12]{Viana2016}} by\footnote{\textcolor{black}{Observe that the potential in Theorem~\ref{thm:ruele} is simply $\phi$.}}
    \[\calL: f \mapsto \sum_{T(y) = x}\frac{e^{\phi(y)}f(y)\1_{R_\delta^i}(y)}{|\det \d T(y)|} = \sum_{T(y) = x}e^{\psi(y)}f(y)\1_{R_\delta^i}(y) ,\]
    when this is well-posed (see Theorem~\ref{thm:ruele}). In particular, this is the case for any function supported on $R_\delta^i$.
\end{notation}

\begin{lemma}\label{lem:Linftyconv}
    Consider a sequence $\{m_\e\}_{\e >0} \subset  \mathcal C^\alpha(R_\delta^i)$, with $\|m_\e\|_{\mathcal C^\alpha} \leq C$ for every $\e >0$. Then, ${\|\calL_\e m_\e - \calL m_\e\|_{L^\infty} \xrightarrow[]{\e \to 0} 0.}$
\end{lemma}
\begin{proof}
Using the usual bounds, we have that
    \begin{align*}
    |\calL_\e m_\e(x) &-\calL m_\e(x)| = \left|\E_\e\left[\sum_{T_\w(y)=x} \frac{e^{\phi(y) }m_\e(y)}{|\det \d T_\w (y)|} - \sum_{T(y) =x}\frac{e^{\phi(y)}m_\e(y)}{|\det \d T(y)|}\right]\right|\\
        &= \left|\E_\e\left[\sum_i \frac{e^{\phi\circ h_i(\w,x)}m_\e\circ h_i(\w, x )}{|\det \d T_\w\circ h_i(\w, x)|} - \frac{e^{\phi\circ h_i(0,x)}m_\e\circ h_i(0, x)}{|\det \d T\circ h_i(0, x)|}\right]\right|\\
        &\leq \sum_i \E_\e\left[\left|\frac{e^{\phi \circ h_i(\w,x)}m_\e(h_i(\w, x))- e^{\phi \circ h_i(0,x)}m_\e(h_i(0, x)}{|\det \d T_\w(h_i(\w, x))|}\right|\right. \\
        &\left. \qquad\qquad + |e^{\phi \circ h_i(\w,x)}m_\e\circ h_i(0,x)|\left|\frac{1}{|\det \d T_\w(h_i(\w, x))|} -\frac{1}{|\det \d T(h_i(0, x))|}\right|\right]\\
        &\leq N \max_i K_iC (\sup|D_\w h_i|\e)^\alpha + C K_i' \sup|D_\w h_i| \e \limsto 0, \quad \text{as }\e \to 0,
    \end{align*}    
where $N = \sup_{(x,\w) \in R_\delta^i\times  \Omega_\e}\# (T_\w^{-1}(\{x\}\cap R_\delta^i)<\infty$, the $K_i$ provide a bound for the term $|\det \d T_\w(h_i(\w, x))^{-1}|$, $C\sup|D_\w h_i|^\alpha \e^\alpha$ are a H\"older-like bound for the difference \[|e^{\phi\circ h_i(\w,x)}m_\e\circ h_i(\w, x)-e^{\phi\circ h_i(0,x)} m_\e\circ h_i(0, x)|,\] and so is $K_i' \sup|D_\w h_i| \e$ for $({|\det \d T_\w(h_i(\w, x))|}^{-1} -{|\det \d T(h_i(0, x))|}^{-1}|)^{-1}$.
\end{proof}

\begin{proposition}\label{prop:meconv}
Let $m_\e :R_\delta^i\to \mathbb R$ be the functions given by Theorem~\ref{thm:me}. There exist $\lambda_0>0$ and $m_0\in \mathcal C^0(R_\delta^i)$, such that $\lambda_\e \xrightarrow[]{\e\to 0}\lambda_0$ and $\left.m_\e\right|_{R^i} \xrightarrow[]{\e\to 0} m_0$ in $\mathcal C^0(R^i)$, with $\mathcal L m_0 = \lambda_0 m_0$.
\end{proposition}
\begin{proof}
 Since $\|m_\e\|_{\mathcal C^\alpha}\leq C,$ there exists $\{\e_n\}_{n\in\mathbb N},$ such that $\e_n \to 0$ and $m_0\in \mathcal C^0(R_\delta^i)$, such that $\|m_{\e_n} - m_0\|_{\mathcal C^0(R_\delta^i)} \to 0.$ We can assume without loss of generality, by restricting to a subsequence if necessary, that $\lambda_{\e_n} \to \lambda_0 \geq 0,$ \textcolor{black}{which can be done since $\lambda_\e \leq \|\mathcal L_\e\|_{\mathcal C^\alpha}$ and the latter is uniformly bounded for $0<\e<\e_0$.}
From Lemma~\ref{lem:Linftyconv} we obtain that
\begin{align*}
    \lambda_0 m_0 = \lim_{n\to\infty} \lambda_{\e_n} m_{\e_n} = \lim_{n\to\infty} \calL_{\e_n} m_{\e_n} = \lim_{n\to\infty} \calL m_{\e_n} = \calL m_0. 
\end{align*}

In the following, we show that $\lambda_0>0$. Since $\int_{R^i_\delta} m_{\e_n}(x)\rho(\d x) = 1$ for every $n\in\mathbb N$, then by the Lebesgue-dominated convergence theorem $\int_{R^i_\delta}m_0(x)\rho(\d x) =1$. Therefore, there exists $x_0\in R^i_\delta$ such that $m_0(x_0) >0.$ Let $C_{x_0}\subset R^i_\delta$ be the connect component of $x_0$ in $R^i_\delta$. From the proof of Theorem~\ref{thm:me} (Step~\ref{step:qsm}), we obtain that for every $n\in \mathbb N$, 
$e^{- \beta d(x_0,y) }m_{\e_{n}}(x_0) \leq m_{\e_{n}}(y),$
for every $y\in C_{x_0}$. Therefore, taking $n\to \infty$ we obtain that $0 < m_0(y)$, for every $y\in C_{x_0}.$ In particular $\left.m_0\right|_{R^i} \neq 0.$ Assume for a contradiction that $\lambda_0=0$. Then, for every $x\in R^i$ \[\sum_{T(y) = x} \frac{e^{\phi(y)} m_0(y)}{|\det \d T(y)|} = \calL m_0 = 0.\]
Since $R^i \subset T^{-1}(R^i)$ the above equation implies that $\left.m_0\right|_{R^i} =0$, which is a contradiction. Therefore $\lambda_0>0$.

Since there exists a unique $m_0 \in \mathcal C^0(R^i)$ such that $\calL m_0(x) = \lambda_0 m_0(x)$ for every $x\in R^i$ and $\lambda_0>0$ ~\cite[Chapter 12]{Viana2016}, the proposition follows.
\end{proof}

\begin{proposition}\label{prop:atra}
    Let $g_\e \in \mathcal C^0_+(R_\delta^i)$ be the functions given by Theorem~\ref{thm:span_g_eps} and consider $\lambda_0>0$ as in Proposition~\ref{prop:meconv}. There exists a probability measure $\gamma$ on $R^i_\delta$ such that ${g_{\e}(x) \dx \xrightarrow{\e\to0}\gamma(\d x)}$ in the {weak}$^*$ topology of $\mathcal M(R_\delta^i)$. Moreover, $\gamma$ is the unique conformal measure for $T$ on $R^i$ for the potential $\phi -\log|\det \d T|$, i.e.~$\gamma$ is the unique probability measure on $R^i$ such that $\mathcal L^* \gamma = \lambda_0 \gamma$.
\end{proposition}
\begin{proof}
    Let $\gamma$ be an accumulation point of $\{g_\e(x) \d x\}_{\e>0}$ in the weak$^*$ topology of $\mathcal M(R_\delta^i)$, i.e.~there exists a sequence $\{g_{\e_n}(x)\d x\}_{n\in \N}$ such that $\e_n \xrightarrow[]{n\to\infty}0$ and ${g_{\e_n}(x) \dx \xrightarrow[]{n\to\infty}\gamma(\d x)}$ in the weak$^*$ topology. 
    We first check that $\gamma$ is a conformal measure on $R_\delta^i$. Indeed, for a test function $f \in \mathcal C^\alpha (R^i_\delta)$ we have:
    \begin{align*}
        (\calL^*\gamma)(f) &= \int_{R_\delta^i}\calL f \d\gamma= \lim_{n\to \infty} \int_{R_\delta^i}\calL f(x) g_{\e_n}(x)\dx\\
        \left(\substack{\text{Lem.~\ref{lem:Linftyconv}}}\right)&= \lim_{n\to\infty} \int_{R_\delta^i}\calL_{\e_n}f(x) g_{\e_n}(x)\dx =\lim_{n\to \infty} \int_{R_\delta^i}f(x) \calP_{\e_n} g_{\e_n}(x)\dx\\
        &=\lim_{n\to \infty} \lambda_{\e_n} \int_{R_\delta^i}f(x) g_{\e_n}(x)\dx =\lambda_0\int_{R_\delta^i} f(x)\gamma(\dx) = \lambda_0 \gamma(f).
    \end{align*}
    We claim that $\supp \gamma \subset R^i$. \textcolor{black}{From item~\ref{it:dyndecshort2} in Lemma~\ref{lem:dyndec_short} and items~\ref{it:preims1} and~\ref{it:preims2} in Lemma~\ref{lem:preims_list}, we obtain that}
\begin{align*}
    1&= \gamma(R_\delta^i) = \frac{1}{\lambda_0} \int_{R_\delta^i} \calL \1_{R_\delta^i} (x) \gamma(\d x) = \frac{1}{\lambda_0} \int_{R_\delta^i} \sum_{T(y)=x} \frac{e^{\phi(y)}\1_{R_\delta^i}(y)}{|\det \d T(y)|} \gamma(\d x)\\
    &= \frac{1}{\lambda_0} \int_{R_\delta^i} \sum_{T(y)=x} \frac{e^{\phi(y)}\1_{R_{\sigma_0\delta}^i}(y)}{|\det \d T(y)|} \gamma(\d x)=  \frac{1}{\lambda_\e} \int_{R_\delta^i} \calL \1_{R_{\sigma_0\delta}^i}(x) \gamma(\d x) = \gamma(R_{\sigma_0\delta}^i).
\end{align*}
Repeating this argument $n$ times we obtain that $\gamma(R_{\sigma_0^n \delta}^i) = 1$ and the claim follows by taking $n\to \infty$. Since there exists a unique measure $\gamma$ in $R^i$ such that $\calL^* \gamma = \lambda_0 \gamma$ (see~\cite[Chapter 12]{Viana2016}), we conclude that $g_\e(x) \d x \xrightarrow{\e\to0} \gamma(\d x)$ in the {weak}$^*$ topology.
\end{proof}

\begin{proposition}{\label{prop:48}}
\textcolor{black}{Assume that $T$ satisfies Hypothesis~\ref{hyp:H1} on $R^i$}. Let $\nu_\e$ be the unique quasi-ergodic measure of the $e^\phi$-weighted Markov process $X_n^\e$ on $\{g_\e >0\}$. Then, $\nu_\e \xrightarrow[]{\e \to 0}\nu_0(x) \coloneqq m_0(x) \gamma(\d x),$ in the weak${}^*$ topology.
Moreover, $\nu_0$ is the unique $T$-invariant equilibrium state for the potential $\phi -\log|\det \d T|$ in $R^i$.
\end{proposition}
\begin{proof}
From Lemma~\ref{lem:uniqueqed}, we have that $\nu_\e(\d x) = g_\e(x) \mu_\e(\d x)$ is the unique quasi-ergodic measure of the $e^\phi$-weighted Markov process $X_n^\e$ on $\{g_\e>0\}$ such that $R^i\subset \supp g_\e$. Since $m_\e \xrightarrow[]{\e\to 0} m_0$ in $\mathcal C^0(R_\delta^i)$ and $g_\e (x) \d x \xrightarrow[]{\e \to 0} \gamma(\d x)$ in the weak$^*$ topology then, $\nu_\e\xrightarrow[]{\e \to 0}\nu_0$ in the weak$^*$ topology. The final part of the proposition follows from well-known results in the thermodynamic formalism for expanding maps (see~\cite[Chapter 12]{Viana2016}).
\end{proof}

We close this section proving Theorem~\ref{thm:local}.


\begin{proof}[Proof of Theorem~\ref{thm:local}]
    Item~\ref{it:local-1} follows from Theorems~\ref{thm:pianigiani} and \ref{thm:unique-qsm}, item~\ref{it:local-2} follows from Theorem~\ref{thm:unique-qsm} and Proposition~\ref{prop:meconv}, item~\ref{it:local-3} follows from Theorem~\ref{thm:span_g_eps}, item~\ref{it:local-4} follows from Lemma~\ref{lem:uniqueqed}, and items~\ref{it:local-5} and~\ref{it:local-6} follow from Proposition~\ref{prop:48}. Thus, only the last statement is left to check.

If $\nu_0^\phi$ is mixing for the map $T:R^i\to R^i$, then $T:R^i\to R^i$ is topologically mixing since $\supp \nu_0^\phi = R^i$. From Lemma~\ref{lem:uniqueqed} we obtain that $\nu_\e^\phi$ is a quasi-ergodic measure of $e^\phi$-weighted Markov process $X_\e^\phi$ on $R^i_\delta$ and the result follows.
\end{proof}

\begin{corollary}
    Assume Hypothesis~\ref{hyp:H1} and that $\left.T\right|_{R^i}$ is topologically mixing. Let $\delta >0$ be small enough. For every $\e >0$ sufficiently small, let $\nu_\e(\d x)$ be the unique quasi-ergodic measure of the $e^\phi$-weighted Markov process $X_n^\e$ on $R^i_\delta$ such that $R^i \subset \supp \nu_\e$. Then, $\nu_\e(\d x) \xrightarrow[]{\e \to 0}\nu_0(\d x)$ in the weak$^*$ topology. Finally, $\nu_0$ is the unique $T$-invariant equilibrium state for the potential $\phi -\log|\det \d T|$ on $R^i$.
\end{corollary}

\begin{proof}
   From Lemma~\ref{lem:uniqueqed} we have that $g_\e(x)\mu_\e(\d x)$ is a quasi-ergodic measure of the $e^\phi$-weighted Markov process $X_n^\e$ on $R_\delta^i$. Together with Theorem~\ref{thm:local} we obtain the result.
\end{proof}

\section{The global problem}\label{sec:global}

We will prove conditioned stochastic stability of equilibrium states on the global repeller $\Lambda$ by studying the quasi-ergodic measure of the $e^\phi$-weighted Markov process $X_n^\e$ on $\Lambda_\delta$ and absorbed in $\partial \coloneqq U \cup (E \setminus M_\delta)$ for some open set $U\subset M$.
As in Section~\ref{sec:local}, let us fix once and for all an $\alpha$-H\"older weight function $\phi: M_\delta \setminus U \to \R$. Moreover, we assume that $(T,\phi,\Lambda)$ satisfies Hypothesis~\ref{hyp:H2}, with $\Lambda$ as in equation~\eqref{eq:lambda}.
 
We start by arguing that restricting the study of quasi-ergodic measures on $\Lambda_\delta$ is sufficient to characterise those on $M_\delta \setminus U$. Then, we decompose $\Lambda_\delta$ into transient and recurrent subsets, the latter being those that contain the original repellers $R^i$. In particular, we show that all the relevant information for the global dynamics follows from the recurrent subset containing the repeller $R^0$ of maximal growth rate. The stochastic stability of global equilibrium states is then inferred via the stochastic stability of equilibrium states around $R^0$.

\begin{proposition}\label{prop:restrictQSMtoLambda}
   Assume that $(T,\phi,\Lambda)$ satisfies Hypothesis~\ref{hyp:H2}. Let $\delta>0$ be sufficiently small and $\mu_\e$ be a quasi-stationary measure of the $e^\phi$-weighted Markov process $X_n^\e$ on $M_\delta\setminus U$. Then, for sufficiently small $\e>0$,
    \begin{enumerate}
    \item\label{it:restrict1} $\mu_\e \ll \rho,$
    \item\label{it:restrict2} $\supp \mu_\e \cap \Lambda_\delta \neq \emptyset$, and
    \item\label{it:restrict3} $\left.\mu_\e\right|_{\Lambda_\delta},$ after normalisation, is a quasi-stationary measure of the $e^\phi$-weighted Markov process $X_n^\e$ on $\Lambda_\delta$.
\end{enumerate}
\end{proposition}
\begin{proof}
Observe that~\ref{it:restrict1} follows directly from the fact that $\calP_\e (x,\d y) \ll \rho(\d y)$ for every $x\in M_\delta\setminus U$.

To show item~\ref{it:restrict2}, arguing for a contradiction, suppose that $\supp\mu_\e \cap \Lambda_\delta = \emptyset$. We claim that there exists $N \in \N$ and $\e>0$ small enough such that for every $x\in \overline{M_\delta\setminus \Lambda_\delta}$ there exists $i\in \{0,1,\ldots, N\}$ such that $T^i_\omega(x) \in U\cup (E\setminus M_\delta)$ for every $\omega\in \Omega_\e$, or in other words, $\tau(x,\omega) \leq N$ for every $\omega\in \Omega_\e$. This is sufficient to prove~\ref{it:restrict2} since, if true, any measurable set $A \subset M_\delta \setminus U$ would be assigned measure zero:
\begin{align*}
     \mu_\e(A)&= \frac{1}{\lambda_\e^N}\int_{M_\delta \setminus U} \calP ^{N}_\e (x, A)\mu_\e(\d x)  \\
    \left(\substack{\text{assumed}\\\supp \mu_\e \cap \Lambda_\delta = \emptyset}\right) &= \frac{1}{\lambda_\e^N}\int_{M_\delta \setminus \Lambda_\delta} \calP_\e^N(x,A) \mu_\e(\dx) \\
     &=\frac{1}{\lambda_\e^N}\int_{M_\delta \setminus \Lambda_\delta} \E_\e\left[ e^{\sum_{i=0}^{N-1} \phi \circ T_\w^i (x) } \1_A \circ T^N_\w(x) \1_{\{\tau(\w,x)>N\} }\right]\mu_\e(\d x) = 0,
\end{align*}
   which is a contradiction. To verify the claim, choose $y \in \overline{M_\delta \setminus \Lambda_\delta}$. Then there exists $n(y) \in \mathbb N$ such that 
$T^n(y) \in U \cup (E \setminus M_\delta)$. Since this is an open set, by continuity of $T$ and $\mathcal C^2$ closeness of the perturbation, there exists $r(y) > 0$ and $\e(y)>0$ such that $T^{n(y)}_\w(B_{r(y)}(y)) \subset U \cup (E \setminus M_\delta)$ for all $\w \in \Omega_{\e(y)}$. Consider a finite open cover of $M_\delta \setminus \Lambda_\delta$ with such balls around $n$ points $y_1, \dots, y_n$ with respective radius $r(y_1), \dots, r(y_n)$. Setting $N = \max \{n(y_1), \dots, n(y_n)\}$, and  $\e = \min\{\e(y_1),\ldots,\e(y_n)\}$ the claim follows.

Finally, we show~\ref{it:restrict3}. Since $T^{-1}(\Lambda)\cap M = \Lambda,$ from the same proof of  Lemma~\ref{lem:preims_list} items~\ref{it:preims2} and~\ref{it:preims3} we obtain that $T_\w^{-1}(\Lambda_\delta) \cap M_\delta \subset \Lambda_\delta$ for every $\w\in\Omega_\e$. Let $A$ be a measurable subset of $\Lambda_\delta$, then
   \begin{align*}
       \int_{\Lambda_\delta} \calP_\e (x,A) \mu_\e(\d x) &=\int_{\Lambda_\delta}  e^{\phi(x)} \E_\e[\1_A \circ T_\w (x)] \mu_\e(\d x)\\&=\int_{M_\delta\setminus U} e^{\phi(x)} \mathbb E_\e[\1_A \circ T_\w (x)] \mu_\e(\d x)\\
       &= \int_{M_\delta\setminus U} \calP_\e(x,A) \mu_\e(\d x) = \lambda_\e \mu_\e(A),
   \end{align*}
   so $\left.\mu_\e\right|_{\Lambda_\delta}$ normalised is a quasi-stationary measure of the $e^\phi$-weighted Markov process $X_\e^\phi$ on $\Lambda_\delta$.
\end{proof}

\begin{proposition}
 Assume that $(T,\phi,\Lambda)$ satisfies Hypothesis~\ref{hyp:H2}. Consider the operator $\calP_\e: L^\infty(M_\delta\setminus U) \to L^\infty(M_\delta\setminus U) $. If ${g \in L_+^\infty(M_\delta\setminus U)}$ is such that $\calP_\e g = \lambda_\e g$, then $\left.g \right|_{M_\delta\setminus \Lambda_\delta}=0$.\label{prop:eigenvector0}
\end{proposition} 
\begin{proof}
 The proof of Proposition~\ref{prop:restrictQSMtoLambda} yields that for $\e>0$ small enough there exists $N$ such that $T_\w^N(x) \in U$ for every $x \in M_\delta \setminus \Lambda_\delta$ and $\w \in \Omega_\e$. Therefore,
 \[\calP^N_\e(x,M_\delta \setminus U) =0\ \text{for every }x\in M_\delta\setminus \Lambda_\delta.\]
 It follows that for every $x\in M_\delta\setminus \Lambda_\delta$ and $N > 0$, 
 \[0\leq g(x) = \frac{1}{\lambda_\e^N} \calP_\e^N g(x)\leq \frac{\|g\|_\infty}{\lambda_\e^N} \calP_\e^N (x,{M_\delta \setminus U})=0,\]
 verifying the claim.
\end{proof}

As a result of Propositions~\ref{prop:restrictQSMtoLambda}
and~\ref{prop:eigenvector0}, it is natural to redefine the operator $\calP_\e$ as 
\begin{align*}
   \calP_\e:L^\infty(\Lambda_\delta) &\to L^\infty (\Lambda_\delta) \\*
   f &\mapsto e^{\phi } \mathbb E_\e [ f\circ T_\w\cdot \1_{\Lambda_\delta}\circ T_\w], 
\end{align*} and denote by $\lambda_\e = r(\calP_\e)$ its spectral radius. Moreover, observe that
\begin{align*}
    \mathcal L_\e:L^1(\Lambda_\delta)&\to L^1(\Lambda_\delta)\\*
    f &\mapsto \mathbb E_\e \left[\sum_{T_\w (y) = x} \frac{ e^{\phi(y) } f(y) \1_{\Lambda_\delta}(y)}{|\det \d T_\w(y)|}\right]
\end{align*}
is well defined and that $\mathcal L_\e^* = \calP_\e.$

\subsection{Recurrent and transient regions}\label{sec:graph}
In this section, we represent the relevant dynamical behaviour of the absorbing Markov process $X_n^\e$ for every $\e > 0$ via a graph whose vertices are the connected components of $\Lambda_\delta$. This approach resembles the graphs constructed via chain recurrence and filtration methods for classical dynamical systems (see~\cite{Crovesier2015, deLeo2021a, deLeo2021b}). Later, we use this construction to characterise the support of the {relevant quasi-stationary measure of the $e^\phi$-weighted Markov process $X_n^\e$}.

Given $\e>0,$ we define an equivalence relation $\sim_\e$ on the set of connected components
\[\Gamma_\delta \coloneqq \{C \subset \Lambda_\delta; \ C\ \text{is a connected component of }\Lambda_\delta\}\] 
as follows: for any $C_1,C_2\in \Gamma_\delta$, we say that $C_1\sim_{\e} C_2$ if
\begin{itemize}
    \item $C_1=C_2$, or
    \item both sets are reachable from each other, i.e.~for every $i,j\in\{1,2\}$, there exist sets $W_0, W_1,\ldots, W_n, W_{n+1} \in\Gamma_\delta$ such that 
    \(\min_{\ell\in\{0,\ldots,n\}}\sup_{x\in W_\ell}\calP_\e (x, W_{\ell+1}) >0, \)
    with $W_0 = C_i$ and $W_{n+1} = C_j$.
\end{itemize}

\begin{proposition}\label{prop:stab}
Assume that $(T,\phi,\Lambda)$ satisfies Hypothesis~\ref{hyp:H2}. The set of equivalence classes $\Gamma_\delta/\sim_\e$ stabilises as $\e\to 0,$ i.e.~there exist $C_1,\ldots,C_n \in \Gamma_\delta$ such that for every $\e$ small enough we have that \[\Gamma_\delta/\sim_\e= \{[C_1],\ldots, [C_n]\},\]
where $[C_i]$ represents the equivalence class of the element $C_i.$
\end{proposition} 

\begin{proof}
The cardinality of $\Gamma_\delta$ is finite and observe that if $0<\e_1<\e_2$, then $C_1 \sim_{\e_1} C_2$ implies $C_1 \sim_{\e_2} C_2$. This ensures that $\Gamma_\delta/\sim_\e$ stabilises as $\e\to 0.$
\end{proof}

\begin{definition}\label{not:mi}
Given $\delta >0$ small enough, let $C_1,\ldots,C_n \in \Gamma_\delta$ be the sets given in Proposition~\ref{prop:stab}. Define 
\[M_i \coloneqq \bigcup_{C\in [C_i]}C,\]
i.e.~$M_i$ is the (disconnected) region spanned by all elements in the class $[C_i]$. Then:
\begin{itemize}
    \item If there exists $j \in \{1, \ldots, k\}$ such that $R^j_\delta \subset M_i$, we say that $M_i$ is a \emph{recurrent region}.
    \item If there are no sets $R_\delta^j$ intersecting $M_i$, we say that $M_i$ is a \emph{transient region}.
\end{itemize}
\end{definition}

\begin{lemma}
  All regions $M_i$ can be classified as either recurrent or transient.
\end{lemma}
\begin{proof}
    Assume that $M_i$ is not a transient region so that there exists $R^j_\delta$ such that such that $R_\delta^j \cap M_i \neq \emptyset$. Then, there exists a connected component $C \Subset R^j_{\delta}$ such that $C \subset M_i$. Since $T$ is topologically transitive on $R^j$ we obtain that $R^j_\delta \subset M_i$ and therefore $M_i$ is recurrent.
\end{proof}

\begin{proposition}\label{prop:GraphDec}
    Let $M_t$ be a transient region, then there exists $N \in \N$ such that for all $x \in M_t$, $\calP^n_\e(x, M_t) = 0 \text{ for all } n \geq N.$ 
\end{proposition}

\begin{proof}
We begin by showing that there exists $N\in\N$ such that for every $x\in M_t$, either:
\begin{itemize}
    \item $T^n(x) \in  \mathrm{Int}\left(\bigcup M_{r}\right)$ for some $n \leq N$, where $\bigcup M_{r}$ is the union of all recurrent regions, or
    \item $T^n(x) \not\in \Lambda_\delta$ for some $n \leq N$.
\end{itemize}

Let $x\in \widetilde{\Lambda} \cap M_t,$ where  
\[\widetilde{\Lambda} \coloneqq \{x\in M_\delta;\ \text{there exists }n\in\mathbb N\ \text{such that }T^n(x) \in R\}.\] 
There exists an open neighbourhood $U_x$ of $x$, such that $T^{n_x}(U_x) \subset \mathrm{Int}\left( \bigcup M_{r} \right)$, the union of recurrent regions. Since $\widetilde\Lambda \cap M_t$ is compact, there exist points $x_1,\ldots, x_{s}$ with respective open neighbourhoods $U_{x_1}, \ldots, U_{x_s}$ such that
\[\widetilde\Lambda \cap M_t \subset \bigcup_{j = 1}^{s}U_{x_j}.\]
Set $N = \max\{n_{x_1},\ldots, n_{x_s}\}.$

On the other hand, observe that for every $y\in M_t\setminus \widetilde\Lambda \subset M_\delta \setminus \widetilde\Lambda$, it follows from $T$ satisfying Hypothesis~\ref{hyp:H2} and~\cite[Theorem 11.2.14]{Viana2016} that there exists $n_y$ such that $T^{n_y}(y) \notin \Lambda_\delta.$ From continuity there exists an open neighbourhood $V_y$ of $y$ such that $T^{n_y}(V_y)\cap \Lambda_\delta =\emptyset.$ Since $M_t\setminus B$ is compact, there exist $y_1,\ldots, y_{m}$ with respective open neighbourhoods $V_{y_1}, \ldots, V_{y_m}$ such that
\[M_t\setminus B \subset \bigcup_{i= 1}^{m}V_{y_i}.\]
Set $N = \max\{n_{y_1},\ldots,n_{y_m}\}.$ From continuity of $(x,\w)\to T_\w (x)$ we obtain that for every $x\in M_t$, either $T^n_\w(x) \in \bigcup M_{r}$, for every $\w\in \Omega_\e$ and some $n \leq N$; or $T^n_\w(x) \notin \Lambda_\delta$, for every $\w\in \Omega_\e$ and some $n \leq N$.
In the first case, allowing return to $M_t$ would join the equivalence classes $[C_t]$ of $M_t$ with $[C_r]$ for some recurrent region $M_r$, contradicting transience. In the second case, once the process escapes $\Lambda_\delta$ it is killed. Thus, $\mathcal P_\e^n(x,M_t) = 0$ for every $n>N$.
\end{proof}

Proposition~\ref{prop:GraphDec} naturally motivates the following definition.

\begin{definition}
Fix $\e>0$ such that the conclusions of Proposition~\ref{prop:stab} hold. Let $M_1,\ldots, M_n$ be the sets introduced in Definition~\ref{not:mi}. We define the directed graph $\mathscr{G}_\delta =( V_\delta, E)$ in the following way:
\begin{itemize}
 \item the set of vertices $V_\delta$ is given by $V_\delta \coloneqq \{M_1,\ldots, M_n\}$,
    \item given $M_i,M_j\in V_\delta$ we say that the edge $M_i \to M_j$ is in $E_\e$ if $M_i\neq M_j$ and there exists $x\in M_{i}$ such that $\mathcal P_\e (x,M_j) >0.$ 
\end{itemize}

Observe that using the same argument as in Proposition~\ref{prop:stab}, the set of edges $E$ does not depend on $\e$ as long as this parameter is small enough.
\end{definition}

\begin{proposition}
 Given a transient region $M_t\in V_\delta$ there exists a path in $\scrG_\delta$ connecting $M_t$ to a recurrent region $M_r$. Moreover, the graph $\scrG_\delta$ is acyclic.    
\end{proposition}
\begin{proof}
To see the first part of the proposition, observe that there exists $x\in M_{t}\cap (\Lambda\setminus R)$. In this way, there exists $n\in\mathbb N$, such that $T^n(x) \in R.$ Defining $M_r$ as the unique recurrent region such that $T^n(x)\in M_r$, we obtain that there exists a path from $M_t$ to $M_r$ in the graph $\scrG_\delta$.

Finally, observe that if $\scrG_\delta$ had a cycle then this would contradict the maximality of the equivalence classes $[C_1],\ldots, [C_n]$.
\end{proof}

\subsection{Proof of the main (global) result}

Recall from Lemma~\ref{lem:dyndec} that $R = \sqcup_{i=1}^{k}R^i$. For every $i\in\{1,\ldots,k\}$, consider the (deterministic) operator
\begin{align*}
    \calL_i: \mathcal C^0(R^i)&\to \mathcal C^0(R^i)\\*
    f&\mapsto \sum_{T(y) = x}\frac{e^{\phi(y)} f(y)}{|\det \d T(y)|},
\end{align*}
and set $\lambda_i = r(\calL_i).$

\begin{notation}
Assume Hypothesis~\ref{hyp:H2}. Given a closed set $A\subset \Lambda_\delta$ we write:
\begin{itemize}
 \item $\calP_{A,\e}: L^\infty(A,\rho) \to L^\infty(A,\rho)$, $\calP_{A,\e} f = \calP_\e (\1_A \cdot f)$,
 \item $\calL_{A,\e}: L^1(A,\rho) \to L^1(A,\rho)$, $\calL_{A,\e} f =\calL_\e (\1_A \cdot f),$ and
 \item for each vertex $M_v$ of the graph $\mathscr G_\delta$ we define
        \begin{align*}
        \calL_{M_v}:  \mathcal C^0(M_v)  &\to \mathcal C^0(M_v) \\*
        f &\mapsto \sum_{T(y)=x}\frac{e^{\phi(y) }f(y)\1_{M_v}(y)}{|\det \d T(x) |}.
    \end{align*} 
Note from Lemma~\ref{lem:preims_list} that this linear operator is well-defined.
\end{itemize}
\end{notation}

\begin{lemma}\label{lem:specradius}
Given a recurrent region $M_v$ we have that \[r(\calP_{M_v, \e}) \xrightarrow[]{\e\to 0} \lambda_{M_v}\coloneqq\max\{ \lambda_i;\ i\in \mathcal I_{M_v}\}, \]
where $\mathcal I_{M_v} \coloneqq\{i\in\{1,\ldots, k\};\ R^i\subset M_v\}.$
\end{lemma}
\begin{proof}

We divide the proof into two steps.
\begin{step}[1]\label{step:ineq_leq}
 $\lambda_{M_v} \leq \liminf_{\e\to 0} r(\calP_{M_v,\e}).$
\end{step}
\begin{stepproof}[Proof of Step~\ref{step:ineq_leq}]
Observe that for every $i\in \mathcal I_{M_v}$ and every non-negative $f\in L^\infty (M_v)$,
\(\1_{R_\delta^i} \calP_{\e} (\1_{R_\delta^i}\cdot f) \leq \calP_{M_v,\e}f.\)
From Theorem~\ref{thm:local} and the above equation we obtain \[\lambda_i = \lim_{\e \to 0 }r(\calP_{R_\delta^i,\e}) \leq \liminf_{\e \to 0} r(\calP_{M_v,\e}),\] for every $i \in \mathcal I_{M_v}.$ 
\end{stepproof}
\begin{step}[2]\label{step:fallen}
 $\limsup_{\e\to\infty}  r(\calP_{M_v,\e}) \leq \lambda_{M_v}.$
\end{step}
\begin{finalstepproof}[Proof of Step~\ref{step:fallen}]
Repeating the same argumentation of Section~\ref{sec:local}, we obtain that:
\begin{enumerate}
    \item there exists $g_\e \in \ker(\calP_{M_v, \e} - r(\calP_{M_v, \e}) ) \cap \mathcal C^0_+(M_v)$ with $\int g_\e \d \rho =1$,
    \item $\ker(\mathcal L_{M_v, \e} - r(\calP_{M_v, \e}) )= \s(m_{\e})$ for some $m_\e \in \mathcal C^\alpha (M_v)$ and $m_\e(x)>0$ for every $x\in M_v$, and
    \item there exists a sequence $\{\e_n\}_{n\in \mathbb N}$ satisfying $\e_n\to 0,$ such that:
    \begin{itemize}
        \item $r(\calP_{M_v, \e_{n}})\xrightarrow[]{n\to\infty}  \lambda_0 = \limsup_{\e \to 0} r(\calP_{M_v,\e})$,
        \item $g_{\e_{n}}(x)\d x\xrightarrow[]{n\to\infty} \gamma(\d x)$ in the weak-$^*$ topology and $\mathcal L_{M_v}^* \gamma = \lambda_0 \gamma,$ and
        \item $m_{\e_{n}} \xrightarrow[]{n\to\infty} m$ in $\mathcal C^0(M_v)$ and $\mathcal L_{M_v} m = \lambda_0 m$.
    \end{itemize}
\end{enumerate}
 It is clear that $\gamma (M_v\cap \Lambda) = 1$. Since $\Lambda = \bigcup_{n\in \N} T^{-n} (R)$, there exists $N\in \N$ such that $\gamma(M_v\cap T^{-N}(R)) >0.$ This implies that
\begin{align*}
 0<\gamma(M_v\cap T^{-N}(R)) &= \frac{1}{\lambda_0^N}\int_{M_v\cap \Lambda} \mathcal L^N_{M_v}\1_{T^{-N}(R)}(x) \gamma(\d x) \\ 
 &= \frac{1}{\lambda_0^N}\int_{M_v\cap \Lambda}  \sum_{T^N(y) = x} \frac{e^{S_N\phi(y)}\1_{R}\circ T^N(y)}{|\det \d T^N(y)|} \gamma(\d x) \\
 & =\frac{1}{\lambda_0^N} \int_{M_v\cap R}  \sum_{T^N(y) = x} \frac{e^{S_N\phi(y)}\1_{R}\circ T^N(y)}{|\det \d T^N(y)|} \gamma(\d x),
\end{align*}
where $S_N \phi(x) = \sum_{i=0}^{N-1} \phi\circ T^i(x)$, therefore $\gamma(M_v\cap R) >0.$ In this way, there exists $R^j \subset M_v$ such that $\gamma(R^j) >0.$ Define $\gamma_j (\d x) \coloneqq \gamma( R^j\cap \d x).$ Given $f\in \mathcal C^0_+(R^j),$ we obtain that
\begin{align*}
\calL_j^*\gamma_j(f) &= \gamma_j(\mathcal L_j f) = \int_{R^j} \sum_{T(y) =x} \frac{e^{\phi(y)}\1_{R^j}(y) f(y)}{|\det \d T(y)|}\gamma (\d x) \\
&= \int_{\Lambda} \mathcal L(\1_{R^j} f) \gamma (\d x) = \lambda_0\gamma(\1_{R^j} f) = \lambda_0 \gamma_j (f)\leq r(\mathcal L_j) \gamma_j(f).
\end{align*}
Since $r(\calL_j) = \lambda_j$, this implies that 
\[\lambda_0 = \limsup_{\e \to 0}  r(\calP_{M_v,\e}) \leq \lambda_j \leq \lambda_{M_v},\]
and we conclude the proof.
\end{finalstepproof}
\begingroup
\renewcommand{\qedsymbol}{\empty}
\end{proof}
\endgroup

\begin{remark}
   Observe that from Theorem~\ref{thm:ruele}, item~\ref{it:hyp2-2} of Hypothesis~\ref{hyp:H2} is equivalent to the existence of $i\in\{1,\ldots,k\}$ such that $\lambda_i > \max_{j\neq i}\lambda_j.$   
\end{remark}

\begin{notation}
If $(T, \phi, \Lambda)$ satisfies Hypothesis~\ref{hyp:H2}, we define $\lambda_0 \coloneqq \max\{\lambda_i;\ i\in\{1,\ldots,k\}\}.$ 
Let $i_0\in\{1,\ldots,k\}$ be the unique natural number such that  $\lambda_{i_0}=\lambda_0.$ We denote by $M_0$ the unique recurrent region such that $R^0 \coloneqq R^{i_0}\subset M_0$.
\end{notation}

\begin{proposition}\label{prop:end}
Assume that $(T, \phi, \Lambda)$ satisfies Hypothesis~\ref{hyp:H2} and let $\e$ be small enough. If $g \in \mathrm{ker}(\calP_\e- \lambda_\e)$, then $\calP_{M_0,\e} (\1_{M_0}g) = \lambda_\e \1_{M_0} g.$  Moreover, for every vertex $M_v$ of $\scrG_\delta$ such that there exists a path from $M_0$ to $M_v$, we have that $\left.g\right|_{M_v} = 0.$ Also, if $\left.g\right|_{M_0} = 0$, then $g (x) = 0$ for every $x\in \Lambda_\delta.$
\end{proposition}
\begin{proof}
First, observe that such a $g$ exists from the Krein-Rutman Theorem~\cite[Theorem~4.1.4]{Meyer-Nieberg1991}.

Let \[V_g \coloneqq\{ M_i;\ \mbox{$M_i$ is a vertex of $\scrG_\delta$ and $M_i\cap \{g\neq 0\}\neq \emptyset$}\}\] and define $\scrG_g \coloneqq (V_g,E_g)\subset \scrG_\delta$ as the maximal subgraph which contains the vertices $V_g$. Since $\scrG_g$ is acyclic, there exists a terminal vertex $M_f \in V_g$, i.e.~no edge in $\scrG_g$ exits from $M_f$. We claim that $M_f = M_0$.

Observe that if $x\in M_f$ and $T_\w(x) \in \{g\neq 0\}$ for some $\w \in \Omega_\e$, then $T_\w(x) \in M_f$. Indeed, if there exists $M_v \in V_g$ such that $T_\w (x) \in M_v,$ then  $M_f \to M_v \in E_g$ but $M_f$ is a terminal vertex. 
This shows the second part of the proposition for $M_f$. It remains to verify that $M_f = M_0.$

We claim that $\calP_{M_f,\e} (\1_{M_f} g) = \lambda_\e \1_{M_f} g$. Indeed, for every $x\in M_f$ we obtain that
\begin{align*}
    \calP_{{M_f},\e}(\1_{M_f} g)(x) &= e^{\phi(x)}\mathbb E_\e [\1_{M_f}\circ T_\w (x)\cdot  g\circ T_\w(x)] \\
    &= e^{\phi(x)} \mathbb E_\e[\1_{M_f \cap \{g\neq 0\}}\circ T_\w (x) \cdot g\circ T_\w(x)]  \\
    &=  e^{\phi(x)}\E_\e[ \1_{\{g\neq 0\}}\circ T_\w (x) \cdot g\circ T_\w(x)] = \calP_\e g(x) = \lambda_\e g(x).
\end{align*}
Taking $\e \to 0$, from Lemma~\ref{lem:specradius} and item~\ref{it:hyp2-2} of Hypothesis~\ref{hyp:H2} we obtain that $M_f = M_0$.
\end{proof} 

\begin{proposition}\label{prop:begin}
Assume that $(T, \phi, \Lambda)$ satisfies Hypothesis~\ref{hyp:H2} and let $\e$ be small enough. We have that, if $m \in \mathrm{ker}(\calL_\e- \lambda_\e)\cap L^1_+(\Lambda_\delta)$, then $\calL_{M_0,\e} (\1_{M_0}m) = \lambda_\e \1_{M_0} m.$ Moreover, for every vertex $M_v$ of $\scrG_\delta$ such that there exists a path from $M_v$ to $M_0$, we have that $\left.m\right|_{M_v} = 0.$
\end{proposition}
\begin{proof}
Again, such an $m$ exists from the Krein-Rutman Theorem~\cite[Theorem~4.1.4]{Meyer-Nieberg1991}. Analogous to the previous proof, let \[V_m \coloneqq\{ M_i;\ \mbox{$M_i$ is a vertex of $\scrG_\delta$ and $M_i\cap \{m> 0\}\neq \emptyset$}\}\]
and define $\scrG_m\subset \scrG_\delta$ as the maximal subgraph which contains the vertices $V_m$. Since $\scrG_m$ is acyclic, there exists an initial vertex $M_s \in V_m$, i.e.~no edge in $\scrG_m$ ends in $M_s$. We claim that $M_s = M_0$. 

Observe that for every $x\in M_s$ and $\w \in \Omega_\e$, \[T_\w^{-1}(M_s)\cap \{m> 0\} = T_\w^{-1}(M_s)\cap M_s\cap \{m> 0\}.\]
This shows the second part of the proposition for $M_s$. It remains to show that $M_s = M_0$.

We claim that $\mathcal L_{{M_s},\e}(\1_{M_s} m) = \lambda_\e \1_{M_s} m.$ In fact, observe that for every $x\in M_s$
\begin{align*}
    \mathcal L_{{M_s},\e}(\1_{M_s} m)(x) &=\mathbb E_\e\left[ \sum_{T_\w(y)=x} \frac{e^{\phi(y)}\1_{M_s}(y) m(y)}{|\det \d T_\w (y)|}\right] = \lambda_\e \1_{M_s}(x) m(x).
\end{align*}
Hence, from the choice of $\e$ we obtain that $M_s=M_0$ and the result follows.
\end{proof} 

\begin{proposition}
Assume that $(T, \phi, \Lambda)$ satisfies Hypothesis~\ref{hyp:H2} and let $\e>0$ be small enough. There exists $g_\e \in \mathcal C^0_+(\Lambda_\delta)$ and $m_\e \in L^1_+(\Lambda_\delta)$ such that:
\begin{enumerate}[label = (\arabic*)]
    \item\label{it:propo-1} $\ker(\calP_\e-\lambda_\e)=\s(g_\e)$, 
    \item\label{it:propo-2} $\ker(\mathcal L_\e-\lambda_\e)=\s(m_\e)$,
    \item\label{it:propo-3} $g_\e(x) > 0$ for every $x \in R^0$, and
    \item\label{it:propo-4} $\1_{M_0}m_\e \in \mathcal C^\alpha(M_0)$ and $m_\e(x) >0$ for every $x\in M_0.$ 
    \item\label{it:propo-5} \textcolor{black}{$\mu_\e(\d x) = m_\e(x) \d x$ is the unique quasi-stationary measure such that $R^0\subset \supp \mu_\e$.}
\end{enumerate}
\label{prop:exist_right_left_efuncs}
\end{proposition}
\begin{proof}
From the same method provided in Theorem~\ref{thm:me}, there exists $\widetilde{m}_\e \in\mathcal C^\alpha(M_0)$ such that $\calL_{{M_0},\e}\widetilde m_\e = \lambda_\e \widetilde m_\e$ and $M_0=\{\widetilde{m}_\e>0\}$.

Given $g_\e \in \ker(\calP_\e - \lambda_\e)$, from Proposition~\ref{prop:end}, we have that $\calP_{{M_0},\e} (\1_{{M_0}} g_\e) = \lambda_\e \1_{{M_0}} g_\e.$ Since {$M_0=\{\widetilde{m}_\e>0\}$}, repeating the same argument as in Theorem~\ref{thm:span_g_eps}, we obtain that 
\begin{equation}\label{eq:plusminus}
    \1_{{M_0}}g_\e^{\pm} \in \ker(\calP_{M_0,\e} - \lambda_\e).
\end{equation}
We divide the remainder of the proof into three steps.

\begin{step}[1]\label{step:gR0}
    For every $\e>0$ sufficiently small, if $\widetilde g_\e \in \ker(\calP_{M_0,\e} - \lambda_\e)$, then $\widetilde g_\e \in \mathcal C^0(\Lambda_{\delta})$ and $\widetilde g_\e(x) >0$ for every $x\in R^0.$
\end{step}
\begin{stepproof}[Proof of Step~\ref{step:gR0}]
 Using the fact that $\calP_{M_0,\e}$ is strong Feller and equation~\eqref{eq:plusminus}, assume for a contradiction that there exists a sequence of positive numbers $\{\e_n\}_{n\in\N}$ such that $\e_n \to 0$, and for every $n \in\N,$ there exists a non-negative function $\widetilde{g}_{\e_n} \in \ker(\calP_{M_0,\e_n} - \lambda_{\e_n})$ such that $\widetilde{g}_{\e_n}(x_n) = 0$ for some $x_n\in R^0.$ 

From the same arguments presented in Steps~\ref{step:positive1} and~\ref{step:positive2} of Lemma~\ref{lem:positive} we have that if $\widetilde{g}_{\e_n}(x) = 0$ for some $x\in R^0$ then $\left.\widetilde{g}_{\e_n}\right|_{R_\delta^{0}} = 0.$ Again, as in the proof of Lemma~\ref{lem:specradius}, up to taking a subsequence of $\{\e_n\}_{n\in\mathbb N}$ we can assume that
    \begin{enumerate}
        \item $r(\calP_{M_0, \e_n})\xrightarrow[]{n\to\infty} \lambda_0$,
        \item $\widetilde{g}_{\e_{n}}(x)\d x\xrightarrow[]{n\to\infty} \gamma(\d x)$ in the weak-$^*$ topology and $\mathcal L_{M_0}^* \gamma = \lambda_0 \gamma,$ and
        \item $\widetilde{m}_{\e_{n}} \xrightarrow[]{n\to\infty} m_0$ in $\mathcal C^0(M_0)$ and $\mathcal L_{{M_0}} m_0 = \lambda_0 m_0$.
\end{enumerate}

Observe that $\gamma(R^0_\delta) =0$ by construction. Repeating the same computations in Step~\ref{step:fallen} of Lemma~\ref{lem:specradius} (now with $\Lambda$ instead of $M$) we obtain that there exists $R^j \subset M_0$ such that $\gamma(R^j) >0$ and $\mathcal L_j^*\gamma(R^j\cap \d x) =\lambda_0 \gamma(R^j\cap \d x)$, contradicting Hypothesis~\ref{hyp:H2} since $r(\calL_j)<\lambda_0$. Therefore, $\widetilde{g}_{\e_n}(x) >0$ for every $x\in R^0$ and $n\in\N$.
\end{stepproof}
\begin{step}[2]\label{step:gdim1} We show that $\dim \ker (\calP_\e - \lambda_\e) \leq 1$.
\end{step}
\begin{stepproof}[Proof of Step~\ref{step:gdim1}]
Let $g_1,g_2 \in \ker (\calP_\e - \lambda_\e)$. Observe that from the same proof of Theorem~\ref{thm:span_g_eps}, we obtain that there exists $t_0$ such that $\left.(g_1 - t_0 g_2)\right|_{R^0} =0$. Since $g_1 -t_0 g_2 \in \ker (\calP_\e - \lambda_\e),$ we have from Step~\ref{step:gR0} that $\1_{M_0} (g_1 -t_0 g_2) =0$. Finally, from Proposition~\ref{prop:end} we obtain that $g_1 -t_0 g_2 =0$.
\end{stepproof}


    

We may now conclude the proof of the proposition.
\textcolor{black}{From the Krein-Rutman theorem~\cite[Theorem~4.1.4]{Meyer-Nieberg1991} and the fact that $\lambda_\e >0$, we obtain that there exists $g_\e \in L^\infty_+(\Lambda_\delta)$ such that $\calP_\e g_\e = \lambda_\e g_\e$. Since $\calP_\e$ is strong Feller we obtain that $g_\e \in \mathcal C^0_+(\Lambda_\delta)$. This provides items~\ref{it:propo-1} and \ref{it:propo-3}. Using the fact that $\mathcal L_\e^*=\calP_\e$, we may choose $m_\e \in L^1(\Lambda_\delta)$ such that $\1_{M_0} m_\e = \widetilde m_\e$, which yields items~\ref{it:propo-2} and \ref{it:propo-4} following the arguments in the beginning of this proof (observe that in Step~\ref{step:gdim1} we apply Theorems~\ref{thm:pianigiani} and \ref{thm:unique-qsm} in place of Theorem~\ref{thm:span_g_eps}, and we use Proposition~\ref{prop:begin} instead of Proposition~\ref{prop:end}.). Finally, item~\ref{it:propo-5} follows from Theorem~\ref{thm:unique-qsm}.}
\end{proof}

\begin{theorem}\label{thm:main_qem}
Assume that $(T, \phi, \Lambda)$ satisfies Hypothesis~\ref{hyp:H2} and let $\e>0$ be small enough. Let $g_\e \in \ker(\calP_\e - \lambda_\e)$ and $m_\e \in \ker(\calL_\e - \lambda_\e)$ be non-negative functions. Then \[\nu_\e^\phi(\d x) = \frac{m_\e (x)g_\e(x) \rho(\d x)}{\int_{\Lambda_\delta} m_\e(y) g_\e(y)\rho(\d y)}\]is the unique quasi-ergodic measure of the $e^\phi$-weighted Markov process $X_\e^\phi$ on $\{m_\e >0\}\cap \{g_\e >0\}$. Moreover, $\nu_\e^\phi \to \nu_0^\phi$ as $\e \to 0$ in the weak-$^*$ topology, where $\nu_0^\phi$ is the unique equilibrium state for $T$ for the potential $\phi-\log |\det \d T|$ supported on $\Lambda.$
\end{theorem}

\begin{proof}
For every $\e>0$ small enough, choose $g_\e \in \mathcal C_+^0(\Lambda_\delta)$ and $m_\e\in L^1_+(\Lambda_\delta)$ satisfying the conclusions of Proposition~\ref{prop:exist_right_left_efuncs}. Following the same strategy as in the proof of Lemma~\ref{lem:uniqueqed} we obtain that
 \[\nu_\e^\phi(\d x) =  
 \frac{g_\e(x)m_\e(x)\rho(\d x)}{\int_{M_0} g_\e(x)m_\e (x) \rho(\d x)}\]
is a quasi-ergodic measure of the $e^\phi$-weighted Markov process $X_\e^\phi$ on $\{g_\e m_\e >0\}.$ From Propositions~\ref{prop:end} and~\ref{prop:begin} we obtain that $R^0 \subset \{g_\e m_\e >0\}\subset M_0$, $\mathcal P_{M_0,\e} \1_{M_0}g_\e = \lambda_\e g_\e$ and $\mathcal L_{M_0,\e} \1_{M_0,\e} m_{\e} = \lambda_\e \1_{M_0} m_\e$. Repeating the proof of Proposition~\ref{prop:48} and Theorem~\ref{thm:local} changing $R^i_\delta$ to $M_0$ we obtain the last part of the result.
\end{proof}

We close this section proving Theorem~\ref{thm:global}.


\begin{proof}[Proof of Theorem~\ref{thm:global}]
Items~\ref{it:thm-exp1} to~\ref{it:thm-exp5} follow directly from Propositions~\ref{prop:restrictQSMtoLambda} and~\ref{prop:exist_right_left_efuncs} and Theorem~\ref{thm:main_qem}.  

We divide the rest of the proof into six steps.

\begin{step}[1]\label{step:properties}
If $\nu_0^\phi$ is topologically mixing, then for every $\e>0$ small enough, the operator
\begin{align*}
    \overline{\mathcal P}_\e : \mathcal{C}^0(M_\delta \setminus U) &\to \mathcal{C}^0(M_\delta\setminus U)\\*
    f&\mapsto  e^{\phi(x)} \mathbb E_\e[f\circ T_\w (x) \cdot \1_{M_\delta\setminus U}\circ T_\w (x) ]
\end{align*}
satisfies the following properties:
\begin{enumerate}[label = (\arabic*)]
    \item\label{it:props1} $\overline{\calP}_\e$ is a strong Feller operator, therefore $\overline{\calP}_\e^2$ is a compact operator,
    \item\label{it:props2} $r(\overline{\calP}_\e) = r(\calP_\e)=\lambda_\e$,
    \item\label{it:props3} there exists $\overline{\mu}_\e \in \mathcal M(M_\delta\setminus U)$ a probability measure such that $\s\{\overline{\mu}_\e\} = \mathrm{ker}(\oPe^*-\lambda_\e)$ and such that $\left.\overline{\mu}_\e\right|_{\Lambda_\delta}/\overline{\mu}_\e(\Lambda_\delta) = \mu_\e$, where $\mu_\e (\dx) \coloneqq m_\e(x) \dx$ is given by Proposition~\ref{prop:exist_right_left_efuncs}, and
    \item\label{it:props4} $\s\{\overline{g}_\e\} = \mathrm{ker}(\oPe-\lambda_\e)$ where $\overline{g}_\e \coloneqq \1_{\Lambda_\delta} g_\e \in \mathcal C^0(M_\delta\setminus U)$, with $g_\e$ given by Proposition~\ref{prop:exist_right_left_efuncs} and $\int \overline{g}_\e \d \overline{\mu}_\e =1.$ 
\end{enumerate}
\end{step}
\begin{stepproof}[Proof of Step~\ref{step:properties}]
Observe that the strong Feller property of $\oPe$ follows by the same computations provided in Proposition~\ref{prop:kernel1}, showing~\ref{it:props1}. Item~\ref{it:props2} follows since $\oPe$ is strong Feller, so\[r( \overline{\mathcal P}_\e: \mathcal C^0(M_\delta\setminus U)\to  \mathcal C^0(M_\delta\setminus U)) = r( \mathcal P_\e: L^\infty (\Lambda_\delta,\rho) \to L^\infty(\Lambda_\delta,\rho))= \lambda_\e.\]
Items~\ref{it:props3} and~\ref{it:props4} are direct consequences of Propositions~\ref{prop:restrictQSMtoLambda},~\ref{prop:eigenvector0} and~\ref{prop:exist_right_left_efuncs}.
\end{stepproof}

\begin{step}[2]\label{step:PB} The operator $\frac{1}{\lambda_\e}\oPe$ is power-bounded, i.e.~$\sup_{n\in\mathbb N} \|\frac{1}{\lambda_\e^n}\oPe^n\|<\infty.$ 
\end{step}
\begin{stepproof}[Proof of Step~\ref{step:PB}]
Repeating the same argumentation as in the proof of Proposition~\ref{prop:restrictQSMtoLambda} item~\ref{it:restrict2}, there exists $N>0$ such that $\mathcal \oPe_\e^N f(x) = 0$ for every $x\in M_\delta\setminus \Lambda_\delta$ and $f\in\mathcal C^0(M_\delta\setminus U)$. In this way, for every $n>0$,
$$ \frac{1}{\lambda_\e^{n+N}}\oPe^{N+n} f = \1_{\Lambda_\delta}\frac{1}{\lambda_\e^{n+N}} \mathcal P_\e^n \left(\1_{\Lambda_\delta} \oPe^{N}  f\right).$$
Since $\frac{1}{\lambda_\e^n}\mathcal P_\e^n $ is power-bounded we obtain the result.
\end{stepproof}

\begin{step}[3]\label{step:absolute_is_eval}
Given a function $f\in \mathcal C^0(M_\delta\setminus U,\C)$ let us define $|f|\in \mathcal C^0(M_\delta\setminus U)$ as the function $x\mapsto \|f(x)\|_{\C}$. Let $\alpha \geq 0$ and $f_\e \in \mathcal C^0(M_\delta\setminus U, \C)$ be such that $\frac{1}{\lambda_\e}\oPe f_\e = e^{i\alpha} f_\e.$ Then for every $x\in \supp \overline{\mu}_\e$, we have that $|f_\e|(x) = \overline{g}_\e(x)\int |f_\e| \d \overline{\mu}_\e$ and $\int |f_\e| \d \overline{\mu}_\e>0$.  
\end{step}
\begin{stepproof}[Proof of Step~\ref{step:absolute_is_eval}]
We have that
\[|f_\e| = \left|e^{i \alpha} f_\e \right| = \left|\frac{1}{\lambda_\e} \oPe f_\e \right| \leq \frac{1}{\lambda_\e} \oPe |f_\e|,\]
therefore, for every $n\in\mathbb N$ we obtain
\[ |f_\e|  \leq \frac{1}{\lambda_\e} \oPe |f_\e| \leq \frac{1}{\lambda_\e^2} \oPe^2 |f_\e| \leq\ldots \leq \frac{1}{\lambda_\e^n} \oPe^n |f_\e|.\]
Since $\oPe^2$ is a compact operator and $\frac{1}{\lambda_\e}\oPe$ is power-bounded from Step~\ref{step:PB}, the above sequence is monotone and bounded. Hence, there exists $g\in \mathcal C^0(M_\delta\setminus U)$ such that $\frac{1}{\lambda_\e^n}\oPe^n |f_\e| \xrightarrow[]{n\to\infty} g\ \text{in }\mathcal C^0(M_\delta\setminus U).$  It follows that $g\in \ker (\oPe -\lambda_\e) = \s\{\overline g_\e\}$. From $0\leq |f_\e|\neq 0$, we obtain that there exists $a>0$ such that $g= a \overline g_\e$. 
Finally, since $|f_\e| \leq g$, both functions are continuous, and their integrals with respect to $\overline{\mu}_\e$ coincide, i.e.
\[\int_M |f_\e|\, \d\overline{\mu}_\e= \int_M g \, \d \overline{\mu}_\e = \int_M a \overline g_\e\, \d\overline{\mu}_\e,\]
it follows that $|f_\e|(x) = \overline g_\e(x) \int_M |f_\e| \d \overline{\mu}_\e $ for every $x\in\supp \overline{\mu}_\e$.
\end{stepproof}

\begin{step}[4]\label{step:spec_gap}
The operator $\oPe$ has the spectral gap property in $\mathcal C^0(M_\delta\setminus U)$, i.e.~there exists a $\oPe$-invariant closed space $W\subset \mathcal C^0(M_\delta\setminus U)$ such that $\mathcal C^0(M_\delta\setminus U) = \s\{g_\e\}\oplus W$ and $r(\left.\oPe\right|_W) <\lambda_\e$. 
\end{step}
\begin{stepproof}[Proof of Step~\ref{step:spec_gap}]
Since $\oPe^2$ is a compact operator and $\frac{1}{\lambda_e^n}\oPe^n$ is power-bounded, it is enough to show that $\sigma_{\mathrm{per}}(\frac{1}{\lambda_\e}\oPe)\cap \mathbb S^1 = \{1\}$ (see details in the proof of Lemma~\ref{lem:alg_mult}). Choose $\alpha \in [0, 2\pi)$ such that $e^{i\alpha}\lambda_\e \in \sigma(\oPe)$. Then there exists $f_\e \in \mathcal C^0(M_\delta\setminus U,\C)$ such that $\frac{1}{\lambda_\e} \oPe f_\e = e^{i\alpha} f_\e$. From Step~\ref{step:absolute_is_eval}, we can assume without loss of generality that $\int |f_\e|\d\overline{\mu}_\e = 1$. Using again Step~\ref{step:absolute_is_eval} and Propositions~\ref{prop:end} and~\ref{prop:begin} we have that, there exists a continuous function $\theta:\{g_\e >0\}\to \R$ such that $f_\e(x) = \overline{g}_\e(x) e^{i\theta (x)} =g_\e(x) e^{i\theta (x)}$ for every $x\in M_0$ and
\[\frac{1}{\lambda_\e}\mathcal P_{M_0,\e} (\1_{M_0} f_\e) = e^{i\alpha} \1_{M_0} f_\e.\]
In this way, for every $n\in\mathbb N$ and $x \in M_0$
\[ e^{i (\theta(x)+n\alpha) }g_\e(x) = \int_{M_0} e^{i\theta(y)} g_\e (y) \frac{1}{\lambda_\e^n} (\mathcal P^n_{M_0,\e})^*(\delta_x) (\d y),\]
which implies that
\[ g_\e(x) = \int_{M_0} e^{i(\theta(y) -\theta(x)-n\alpha)  } g_\e (y) \frac{1}{\lambda_\e^n} (\mathcal P^n_{M_0,\e})^*(\delta_x) (\d y).\]
Since
\[ g_\e(x) = \int_{M_0} g_\e (y) \frac{1}{\lambda_\e^n} (\mathcal P^n_{M_0,\e})^*(\delta_x) (\d y),\]
we obtain that  $e^{i(\theta(y) -\theta(x)-n\alpha) } =1$ for every $y\in \supp\{(\mathcal P^n_{M_0,\e})^*(\delta_x)\}\cap \{g_\e >0\}.$ By hypothesis, the measure $\nu_0^\phi$ is mixing for the map $T:R^0\to R^0$, so $T:R^0\to R^0$ is topologically mixing and therefore topologically exact. Hence, there exists $n_0 \in \N$ such that 
\[R^0 \subset \supp\{(\mathcal P^n_{M_0,\e})^*(\delta_x)\}\cap \{g_\e >0\},\ \text{for every }n>n_0.\]
This implies that $e^{i(\theta(x) - \theta(y) - n\alpha)}=1\ \text{for every } n>n_0$ and $x,y\in R^0$, so $\alpha =0$.
\end{stepproof}

\begin{step}[5]\label{step:phi_pen_qem}
 We show that $\nu_\e^\phi(\d x) = g_\e (x) \mu_\e(\d x) = \overline{g}_\e (x) \overline{\mu}_\e(\d x)/\int \overline{g}_\e(y) \overline{\mu}_\e(\d y)$  is a quasi-ergodic measure of the $e^\phi$-weighted Markov process $X_\e^\phi$ on $M_\delta\setminus U$.
\end{step}
\begin{stepproof}[Proof of Step~\ref{step:phi_pen_qem}]
From Step~\ref{step:properties} and Propositions~\ref{prop:restrictQSMtoLambda},~\ref{prop:eigenvector0} and~\ref{prop:exist_right_left_efuncs} it is clear that $\nu_\e^\phi(\d x) = g_\e (x) \mu_\e(\d x) = \overline{g}_\e (x) \overline{\mu}_\e(\d x)/\int \overline{g}_\e(y) \overline{\mu}_\e(\d y)$. From Steps~\ref{step:absolute_is_eval} and~\ref{step:spec_gap} we obtain that for every bounded and measurable function $h:M_\delta\setminus U \to \R$,
\[\frac{1}{\lambda^n} \oPe^n h \xrightarrow[]{n\to\infty} \overline g_\e \int_{M_\delta\setminus U} h(y) \overline{\mu}_\e(\d y) \ \text{in }\mathcal C^0(M_\delta\setminus U),\]
since $\oPe h\in \mathcal C^0(M_\delta \setminus U). $

Recall that $\tau^\phi = \min\{n; X_n^{\e} \in (E\setminus M_\delta)\cup U\}$ and by construction of the operator $\oPe$, for every $x\in \{g_\e>0\}\cap \supp \mu_\e = \{\overline g_\e>0\}\cap \supp \overline{\mu}_\e$ and for every $n \in \N$
\begin{equation}
\mathbb E_x^\phi \left[\frac{1}{n}\sum_{i=0}^{n-1}h\circ X_i^\e \, \bigg| \,  \tau^\phi >n\right] =  \frac{\lambda_\e^n}{\oPe^n \mathbbm 1_{M_\delta\setminus U} (x)} \frac{1}{n}\sum_{i=0}^{n-1}\frac{1}{\lambda_\e^i}\oPe^i\left(h \frac{1}{\lambda_\e^{n-i}}\oPe^{n-i}\mathbbm 1_{M_\delta\setminus U}\right)(x).
\end{equation}
Since $\frac{1}{\lambda_\e^n}\oPe^n \mathbbm 1_{M_\delta\setminus U} (x) \xrightarrow[]{n \to \infty} \overline{g}_\e(x),$
it is enough to show that
 \[\frac{1}{n}\sum_{i=0}^{n-1}\frac{1}{\lambda_\e^i}\oPe^i\left(h \frac{1}{\lambda_\e^{n-i}}\oPe^{n-i}\mathbbm 1_{M_\delta\setminus U}\right)(x) \xrightarrow[]{n\to\infty} \overline g_\e (x) \int h(y) \overline g_\e(y) \overline \mu_\e(\d y).\]
This holds true since
\begin{align*}
 \frac{1}{n}\sum_{i=0}^{n-1}\frac{1}{\lambda_\e^i}\oPe^i\left(h \frac{1}{\lambda_\e^{n-i}}\oPe^{n-i}\mathbbm 1_{M_\delta\setminus U}\right)(x) =&\frac{1}{n}\sum_{i=0}^{n-1}\frac{1}{\lambda_\e^i}\oPe^i\left(h \left( \frac{1}{\lambda_\e^{n-i}}\oPe^{n-i}\mathbbm 1_{M_\delta\setminus U} - \overline g_\e\right)\right)(x)\\*
 &+ \frac{1}{n}\sum_{i=0}^{n-1}\frac{1}{\lambda_\e^i}\oPe^i\left(h  \overline g_\e\right)(x),
\end{align*} 
and 
\[ \frac{1}{\lambda_\e^i}\oPe^i\left(h  \overline g_\e\right)(x) \xrightarrow[]{i\to\infty } \overline g_\e(x) \int h(y) \overline g_\e (y) \overline\mu_\e(\d y). \]
\end{stepproof}

We may now conclude the proof of the theorem. To do so, we need to show that if $\supp \nu_0^\phi = R^0 \subset \mathrm{Int}(M\setminus U)$, then $\nu_\e^\phi$ is a quasi-ergodic measure of the $e^\phi$-weighted Markov process $X_\e^\phi$ on $M\setminus U$. Redefine the operator $\oPe$ as
\begin{align*}
    \overline{\mathcal P}_\e : \mathcal{C}^0(M \setminus U) &\to \mathcal{C}^0(M\setminus U)\\*
    f&\mapsto  e^{\phi(x)} \mathbb E_\e[f\circ T_\w (x) \cdot \1_{M\setminus U}\circ T_\w (x) ].
\end{align*}
Observe that since $R^0 \subset \mathrm{Int}(M\setminus U)$, we can choose $\delta>0$ small enough such that $M_0 \subset M\setminus U$. Repeating Steps~\ref{step:properties},~\ref{step:PB},~\ref{step:absolute_is_eval} and~\ref{step:spec_gap} we obtain that
\begin{enumerate}[label = (\roman*)]
    \item $\overline{\calP}_\e$ is a strong Feller operator,
    \item $r(\overline{\calP}_\e) = r(\calP_\e)=\lambda_\e$,
     \item there exists a probability measure $\overline{\mu}_\e$ on $M_\delta\setminus U$ such that $\s\{\overline{\mu}_\e\} = \mathrm{ker}(\oPe^*-\lambda_\e)$ and $\left.\overline{\mu}_\e\right|_{M_0}/\overline{\mu}_\e(M_0) = \left.\mu_\e\right|_{M_0}/\mu_\e(M_0).$
    \item $\s\{\overline{g}_\e\} = \mathrm{ker}(\oPe-\lambda_\e)$ and $\1_{M_0} \overline{g}_\e = \1_{M_0} g_\e $, with $g_\e$ given by Proposition~\ref{prop:exist_right_left_efuncs} and $\int \overline{g}_\e \d \overline{\mu}_\e =1.$ 
    \item $\oPe:\mathcal C^0(M\setminus U)\to \mathcal C^0(M\setminus U)$ has the spectral gap property.
\end{enumerate}

As in Step~\ref{step:phi_pen_qem}, we obtain that $\nu_\e^\phi(\d x) = g_\e (x) \mu_\e (\d x) = \overline{g}_\e (x) \overline{\mu}_\e (\d x)$ is a quasi-ergodic measure of the $e^\phi$-weighted Markov process $X_\e^\phi$ on ${M\setminus U}$.
\end{proof}

\addtocontents{toc}{\protect\setcounter{tocdepth}{0}}

\section*{Acknowledgements}
    We would like to thank Matteo~Tabaro and Martin~Rasmussen for valuable comments and insightful discussions, as well as Mark~Demers and Alex~Blumenthal for their feedback and suggestions.
    The authors gratefully acknowledge support from the EPSRC Centre for Doctoral Training in Mathematics of Random Systems: Analysis, Modelling and Simulation (EP/S023925/1). MMC has been supported by an Imperial College President’s PhD scholarship. JSWL has been supported by the EPSRC (EP/Y020669/1) and thanks IRCN (Tokyo) and GUST (Kuwait) for their research support. 



\addtocontents{toc}{\protect\setcounter{tocdepth}{2}}

\bibliographystyle{abbrv}
\bibliography{refs}

\newpage

\input{appendix}

\end{document}

%% file: Settings/packages.tex
\usepackage{amssymb,amsthm,amsmath,mathrsfs}
\usepackage{mathtools}
\usepackage{bbm}
\usepackage{stackengine}
\usepackage{enumitem}
\usepackage{soul}
\usepackage{thm-restate}

\usepackage[foot]{amsaddr}

\usepackage{caption,color,graphicx}
\usepackage[dvipsnames]{xcolor}
\usepackage{soul}

\usepackage{hyperref}
\hypersetup{colorlinks  = true,
            urlcolor    = blue,
            linkcolor   = red,
            citecolor   = blue}

\usepackage{verbatim}

\usepackage{environ}
\mathtoolsset{showonlyrefs,showmanualtags}

\usepackage[centering, top=3cm, bottom=2.5cm]{geometry}

%% file: Settings/commands.tex
\newcommand{\N}{\mathbb{N}}

\newcommand{\R}{\mathbb{R}}
\newcommand{\C}{\mathbb{C}}
\newcommand{\E}{\mathbb{E}}
\renewcommand{\P}{\mathbb{P}}
\newcommand{\1}{\mathbbm{1}}

\newcommand{\supp}{\operatorname{supp}\,}

\newcommand{\Id}{}
\renewcommand{\d}{\operatorname{d}\!}
\newcommand{\dx}{\operatorname{d}\!x}

\newcommand{\Leb}{\operatorname{Leb}}

\newcommand{\calP}{\mathcal{P}}
\newcommand{\calL}{\mathcal{L}}
\newcommand{\e}{\varepsilon}
\newcommand{\w}{\omega}
\newcommand{\limsto}{\longrightarrow}

\newcommand{\s}{\operatorname{span}}

\newcommand{\scrG}{\mathscr G}

\newcommand{\oPe}{\overline{\mathcal P}_\varepsilon}

\newcommand{\dist}{\operatorname{dist}}

\newcommand{\vertiii}[1]{{\left\vert\kern-0.25ex\left\vert\kern-0.25ex\left\vert #1 \right\vert\kern-0.25ex\right\vert\kern-0.25ex\right\vert}}

\let\oldtocsection=\tocsection
\let\oldtocsubsection=\tocsubsection
\renewcommand{\tocsection}[2]{\hspace{0em}\oldtocsection{#1}{#2}}
\renewcommand{\tocsubsection}[2]{\hspace{1em}\oldtocsubsection{#1}{#2}}

\theoremstyle{plain}
\newtheorem*{theorem*}{Theorem}
\newtheorem{theorem}{Theorem}[section]
\newtheorem{lemma}[theorem]{Lemma}
\newtheorem{proposition}[theorem]{Proposition}
\newtheorem{corollary}[theorem]{Corollary}

\theoremstyle{definition}
\newtheorem{definition}[theorem]{Definition}
\newtheorem{remark}[theorem]{Remark}

\newtheorem{notation}[theorem]{Notation}

\newenvironment{stepproof}[1][]{%
  \renewcommand{\qedsymbol}{\stepqedsymbol}%
  \begin{proof}[#1]%
}{%
  \end{proof}%
}

\newenvironment{finalstepproof}[1][]{%
  \renewcommand{\qedsymbol}{\stepqedsymbol $\square$}%
  \begin{proof}[#1]%
}{%
  \end{proof}%
}

\newtheorem{manualtheorem}{Theorem}
\newenvironment{mytheorem}[1][]{%
    \renewcommand\themanualtheorem{#1}%
    \begin{manualtheorem}\it
}{\end{manualtheorem}}

\newtheorem{manualcorollary}{Corollary}
\newenvironment{mycorollary}[1][]{%
    \renewcommand\themanualcorollary{#1}%
    \begin{manualcorollary}\it
}{\end{manualcorollary}}
  
\newtheorem{manualstep}{\bf Step}
\newenvironment{step}[1][]{%
    \renewcommand\themanualstep{#1}%
    \begin{manualstep}\itshape
}{\end{manualstep}}

\newtheorem{manualhyp}{Hypothesis}
\newenvironment{hypothesis}[1][]{%
    \renewcommand\themanualhyp{#1}%
    \begin{manualhyp}\itshape
}{\end{manualhyp}}

\IfPackageLoadedTF{MnSymbol}{\newcommand{\coloneqq}{:=}}{}


%% file: Settings/authors.tex
\author[Bernat Bassols-Cornudella]{Bernat Bassols Cornudella$^{1}$}
\email{\href{mailto:bernat.bassols-cornudella20@imperial.ac.uk}{bernat.bassols-cornudella20@imperial.ac.uk}}

\author[Matheus M.~Castro]{Matheus M. Castro$^{1}$}
\email{\href{mailto:m.de-castro20@imperial.ac.uk}{m.de-castro20@imperial.ac.uk}}

\author[Jeroen S.W.~Lamb]{Jeroen S.W. Lamb$^{1,2,3}$}
\email{\href{mailto:jeroen.lamb@imperial.ac.uk}{jeroen.lamb@imperial.ac.uk}}
\address{$^{1}$Department of Mathematics, Imperial College London, London SW7 2AZ, UK}
\address{$^{2}$International Research Center for Neurointelligence, The University of Tokyo, Tokyo 113-0033, Japan}
\address{$^{3}$Centre for Applied Mathematics and Bioinformatics, Department of Mathematics and Natural Sciences, Gulf University for Science and Technology, Halwally 32093, Kuwait}

%% file: appendix.tex
\appendix 
\renewcommand{\theequation}{A.\arabic{equation}}
\section{Quasi-ergodic measures for a class of strong Feller Markov chains}\label{sec:appendix}

{In this appendix, we provide sufficient conditions for the existence and uniqueness of quasi-ergodic measures of $e^\phi$-weighted Markov processes. We prove Theorems~\ref{thm:QED} and~\ref{thm:QED2}, which are essential for the proof of Lemma~\ref{lem:uniqueqed}}.
The results below employ techniques of absorbing Markov processes theory~\cite{Castro2021, Champagnat2022, Castro2023} and Banach Lattice theory~\cite{Meyer-Nieberg1991}. 

{\color{black}A significant body of work addresses the existence and uniqueness of quasi-stationary and quasi-ergodic distributions for strong Feller Markov processes. Often, this is established under a combination of Lyapunov-type and irreducibility assumptions (see, for instance,~\cite{Champagnat2016,ChamognatVillemonais2023a,Guillin2024a,Guillin2024b} and the references therein). The assumptions in this appendix are of a different nature: they are purely spectral, formulated in terms of a spectral gap for the annealed Koopman operator $\mathcal P$, and in the discrete time setting. In particular, we do not construct Lyapunov functions and do not impose any global transitivity or irreducibility conditions on the dynamics. These hypotheses are deliberately tailored to the specific deterministic–random dynamical systems considered in this paper.}

Let $M$ be a compact metric space, and consider an absorbing Markov process $X_n$ on $E=M\sqcup \partial$ absorbed at $\partial$. For every $x \in M$ and function $f \in L^1(M, \mu)$, we denote by $\E_x[f \circ X_1]$ the expected value of the observable $f$ after one iterate of the process starting at $X_0 = x$. Define the annealed Koopman operator as
\begin{align*}
    \calP: L^\infty(M,\mu)&\to L^\infty(M,\mu)\\*
    f&\mapsto e^{\phi(x)} \mathbb E_x[f\circ X_1 \cdot \mathbbm 1_M \circ X_1].
\end{align*}
Throughout this section, we assume that $\mu$ is a probability measure on $M$ and $\phi: M\to \R$ is a continuous function. The assumptions on $\mathcal P$ exploited in this appendix are:
\begin{hypothesis}[HA]\label{hyp:A}\leavevmode
\begin{enumerate}
    \item $\mathcal P$ is strong Feller, i.e.~given $f\in L^\infty(M,\mu)$ then $\mathcal P f \in \mathcal C^0(M)$,
    \item $\dim \ker (\calP - \lambda) =1$, where $\lambda = r(\calP)$,
    \item there exists $\mu \in \mathcal M_+(M)$ and $g\in \mathcal C^0_+(M),$ such that $\calP^* \mu = \lambda \mu$ and $\calP g = \lambda g$ and $\int g\, \d \mu = \mu(\{g>0\})$, and
    \item $\supp \mu = M$.
\end{enumerate}
\end{hypothesis}

\begin{notation}
Given $n\in\mathbb N$ and $x\in M$ we write $\mathcal P^n(x,\d y)$ for the unique measure on $M$ such that $\mathcal P^n(x,A) = \mathcal P^n\1_A (x)$ for every measurable set $A\subset M$. Observe that $\mathcal P^n(x,\d y)$ is well defined since $\mathcal P(L^\infty(M,\mu))\subset \mathcal C^0(M).$
\end{notation}

\subsection{Spectral properties of \texorpdfstring{$\mathcal P$}{P}}\label{app:props_P}

{We begin by recalling a classical lemma in the theory of Markov processes and prove a series of results characterising the spectrum of $\mathcal P$.}

\begin{lemma}[{\cite[Chapter 1, Lemma 5.10 and 5.11]{Revuz1984}}]\label{lem:SF}
The operator $\mathcal P^n:L^\infty(M,\mu)\to L^\infty(M,\mu)$ is compact for every $n>1.$ 
\end{lemma}
\begin{proof}
    We follow closely the proofs in~\cite[Chapter~1, Lemmas~5.10 and 5.11]{Revuz1984}. Let $\{f_i\}_{i\in\mathbb N}\subset L^\infty(M,\mu)$ be a sequence of functions such that $\|f_i\|_{L^\infty(M,\mu)} \leq 1$. From the Banach–Alaoglu theorem~\cite[Theorem~3.15]{Rudin1991-book}, there exists a subsequence $\{f_{i_k}\}_{k\in\mathbb N} \subset \{f_i\}_{i\in \N}$ and a function $f\in L^\infty(M,\mu)$ such that $f_{i_k}\xrightarrow[]{k\to\infty} f$ in the weak-$^*$ topology of $L^\infty(M,\mu)$.
    For every $\ell \in \mathbb N$, observe that $(\calP^\ell)^*\delta_x(\cdot) \coloneqq \mathcal P^\ell(x, \cdot) \ll \mu$ since for a measurable set $A$ such that $\mu(A) = 0$, we have that
    \[\left|(\mathcal P^\ell)^*\delta_x(A)\right| = \left|\mathcal P^\ell\1_A(x)\right| \leq \|\mathcal P^\ell\1_A\|_{\infty} \leq \|\mathcal P\| \|\1_A\|_{L^\infty(M,\mu)} = 0,\]
    since $\1_A(x) = 1$ only for $x$ in a $\mu$-null measure set.
    Moreover, the Radon-Nykodim derivative $g_x$ of $(\mathcal P^\ell)^*\delta_x$ with respect to $\mu$ satisfies
    \[\|g_x\|_{L^1(M, \mu)} = \int g_x(y)\mu(\d y) = \int (\mathcal P^\ell)^*\delta_x(\d y) \leq \sup_{x\in M\setminus \partial}\mathcal  P^\ell \1(x) < \infty,\]
    namely $g_x \in L^1(M,\mu)$. Thus,
    \begin{align}\lim_{k \to \infty}\calP^\ell f_{i_k}(x) &= \lim_{k \to \infty} \int \calP^\ell f_{i_k}(y)\delta_x(\d y) = \lim_{k \to \infty} \int f_{i_k}(y) (\calP^\ell)^*\delta_x(\d y)\notag \\
        &= \lim_{k \to \infty} \int f_{i_k}(y)g_x(y)\mu(\d y) = \int f(y) g_x(y) \mu(\d y) \notag \\
        &= \int f(y)(\calP^\ell)^*\delta_x(\d y) = \mathcal P^\ell f(x).\notag
    \end{align}
For every $m\in \N$, define the bounded and measurable function $h_m:M\to \R$ as
$$h_m(x) \coloneqq \sup_{j>m} |\calP f_{i_j}(x) -\calP f (x)|.$$
We obtain that, for every $n, k\in \N$ and $x\in M$, $|\mathcal P^n f_{i_k}(x) - \mathcal P^n f(x)| \leq \mathcal P^{n -1} h_{k}(x).$ Finally, since for every $n>1$ the sequence $\{\mathcal P^{n-1} h_k\}_{k\in\mathbb N}$ is a monotonically decreasing sequence of continuous functions converging pointwise to $0$, Dini's lemma~\cite[p.~199]{Friedman1971-book} yields $\mathcal P^n f_{i_k} \xrightarrow[]{k\to\infty} \mathcal P^n f$ in $L^\infty (M,\mu)$.
\end{proof}

\begin{lemma}\label{lemma:point_periph_spec}
Let $\lambda= r(\mathcal P)$ denote the spectral radius of $\mathcal P.$ Then, there exists $k\in\N$ such that $\sigma_{\mathrm{per}}(\mathcal P) = \left\{\lambda e^{2\pi i j/k}\right\}_{j=0}^{k-1}$ where $\sigma_{\mathrm{per}}(\mathcal P)\coloneqq\{\alpha \in \C;\, \|\alpha\|_{\mathbb C} = r(\mathcal P)\ \text{and }\mathrm{ker}(\mathcal P - \alpha) \neq \{0\}\}$ denotes the point peripheral spectrum of $\mathcal P.$
\end{lemma}

\begin{proof}
We divide the proof into three steps:
\begin{step}[1]\label{step:app_abs}
{If} $f \in \ker(\mathcal P -\lambda e^{i\beta})$ for some $\beta>0$ then $|f| \in \s\{g\},$ where $|f|:M\to \R_+$, $|f|(x) = \|f(x)\|_{\mathbb C}$ and $\|\cdot\|_\C$ denotes the complex norm.
\end{step}
\begin{stepproof}[Proof of Step~\ref{step:app_abs}]
Since $\mathcal P$ is a positive operator $|f| = |e^{i \beta} f| = \frac{1}{\lambda
}|\mathcal P f|\leq \frac{1}{\lambda}\mathcal P|f|.$ Moreover,
$$0\leq \int_{M} \frac{1}{\lambda} \mathcal P|f| - |f| \, \d \mu =  \int_{M} |f| \d \mu - \int_M |f| \d\mu = 0.$$
Since $\supp \mu = M$ and $|f|$ is continuous, then $|f| \in \ker (\calP -\lambda) = \s \{g\}.$
\end{stepproof}

\begin{step}[2]\label{step:app_phases}
{I}f $e^{i\beta_1},e^{i\beta_2}\in\sigma_{\mathrm{per}}(\frac{1}{\lambda}\mathcal P)$ for some $\beta_1,\beta_2 >0$ then $e^{i (\beta_1 + \beta_2)} \in\sigma_{\mathrm{per}}(\frac{1}{\lambda}\mathcal P)$.
\end{step}
\begin{stepproof}[Proof of Step~\ref{step:app_phases}]
Given $j\in\{1,2\}$, let $f_j \in \mathrm{ker}(\mathcal P-\lambda e^{i\beta_j})$. From Step~\ref{step:app_abs} and rescaling $f_j$, if necessary, there exists a measurable function $\theta_j:M\to \mathbb R$ such that $f(x) = e^{i\theta_j(x)}g(x)$. 

Hence, for every $x\in M$
$$e^{i\beta_j} f_j(x) = e^{i\beta_j} \left(e^{i\theta_j(x)} g(x)\right) =\frac{1}{\lambda}\calP \left( e^{i\theta_j(x)} g\right)(x) =\frac{1}{\lambda} \int_{M} e^{i\theta_j(y)} g(y) \calP(x,\d y), $$
implying that
$$g(x) = \frac{1}{\lambda} \int_{M} e^{i(\theta_j(y)-\theta_j(x) - \beta_j)} g(y) \calP(x,\d y).$$
Since $g(x)\geq 0$ and $\lambda g(x) =\int_M g(y) \calP(x,\d y),$ we obtain that $e^{i(\theta_j(y)-\theta_j(x) - \beta_j)}=1$, for $\mathcal P(x,\cdot)$-almost every $y\in \{g>0\}.$

Finally, observe that {by} defining $h(x) \coloneqq e^{i ( \theta_1(x) + \theta_2(x))} g(x)$ we obtain that
\begin{align*}
    \mathcal Ph (x) &=\int_{M} e^{i(\theta_1(y)+ \theta_2(y))} g(y) \calP(x,\d y) \\&
    = \int_{M} e^{i(\theta_1(x) + \beta_1 + \theta_2(x) + \beta_2)} g(y) \calP(x,\d y) = e^{i (\beta_1+\beta_2)}\lambda h(x),
\end{align*}
which implies that $e^{i (\beta_1+\beta_2)}\in \sigma_{\mathrm{per}}(\frac{1}{\lambda}\mathcal P)$.
\end{stepproof}

We may now conclude the proof of the lemma. From Step~\ref{step:app_phases}, it is enough to show that $\sigma_{\mathrm{per}}(\mathcal P)$ is finite. Lemma~\ref{lem:SF} implies that $\mathcal P^2$ is a compact operator and therefore $\sigma_{\mathrm{per}}(\mathcal P^2)$ is finite. Finally, since
$\left\{\lambda^2; \lambda\in \sigma_{\mathrm{per}}(\mathcal P)\right\} \subset \sigma_{\mathrm{per}}(\mathcal P^2),$
we obtain that $\sigma_{\mathrm{per}}(\mathcal P)$ is also finite.
\end{proof}

From here onwards, let $k\in \N$ {be fixed as in Lemma}~\ref{lemma:point_periph_spec}.

\begin{lemma}\label{lem:powercompact}
   The sequence $\{\frac{1}{\lambda^n}\mathcal P^n:\mathcal C^0(M) \to \mathcal C^0(M)\}_{n\in\N}$ is bounded, i.e.~${\sup_{n\in \N}\|\frac{1}{\lambda^n}\mathcal P^n\| <\infty.}$
\end{lemma}

\begin{proof}
 We dive the proof into three steps. 

\begin{step}[1]\label{step:11} We show that $\sigma_{\mathrm{per}}(\frac{1}{\lambda}\mathcal P) =\sigma_{\mathrm{per}}(\frac{1}{\lambda^{k+1}}\mathcal P^{k+1})  $.
\end{step}
\begin{stepproof}[Proof of Step~\ref{step:11}]
 Let us consider $\beta \in (0,2\pi)$ such that $ e^{i\beta}\in \sigma_{\mathrm{per}}(\frac{1}{\lambda^{k+1}}\mathcal P^{k+1})\setminus \sigma_{\mathrm{per}}(\frac{1}{\lambda}\mathcal P)$. Since $\calP^{k+1}$ is a compact operator, there exists $f\in \ker(\mathcal P^{k+1} - \lambda^{k+1}e^{i\beta})$. Observe that for every $j\in\{0,1,\ldots,k\},$ we obtain that
 $$0 = \left( \frac{1}{\lambda^{k+1}}\mathcal P^{k+1} - e^{i\beta}\right)f =  \left( \frac{1}{\lambda}\mathcal P - e^{\frac{i\beta}{k+1}+\frac{2\pi i j}{k+1} }\right)\sum_{\ell=0}^{k} {e^{\frac{i \beta (k-\ell) }{k+1}+\frac{2\pi i j (k-\ell)}{k+1} } }\frac{1}{\lambda^\ell}\calP^{\ell}f. $$

From Step~\ref{step:app_phases} of Lemma~\ref{lemma:point_periph_spec}, we have $\gamma \coloneqq e^{\frac{i\beta}{k+1}+\frac{2\pi i j}{k+1}}\not\in \sigma_{\mathrm{per}}(\frac{1}{\lambda}\mathcal P)$ for every $j\in\{0,1,\ldots,k\},$ as otherwise we would have $\gamma^{k+1}=e^{i\beta} \in \sigma_{\mathrm{per}}(\frac{1}{\lambda}\calP).$ Moreover, summing over $j= 0, \dots, k$ yields
\begin{align*}
    0 &= \sum_{j = 0}^{k} \sum_{\ell=0}^{k} {e^{\frac{i \beta (k-\ell) }{k+1}+\frac{2\pi i j (k-\ell)}{k+1} } }\frac{1}{\lambda^\ell}\calP^{\ell}f = \sum_{\ell = 0}^{k} e^{\frac{i\beta (k-\ell)}{k+1}}\sum_{j=0}^{k}e^{\frac{2\pi j (k-\ell)}{k+1} }\frac{1}{\lambda^\ell}{\calP^\ell}f \\ & = \sum_{j = 0}^{k}\frac{1}{\lambda^k}\calP^kf= (k+1)\frac{1}{\lambda^k}\calP^k f,
\end{align*}
where the third equality follows from the only nonzero term $\ell = k$. Applying $\frac{1}{\lambda}{\calP}$ on both sides and since $f \in \ker\{\calP^{k+1} - \lambda^{k+1}\}$ we obtain $f= 0$, a contradiction. The other direction is trivial.
\end{stepproof}


\begin{step}[2]\label{step:12} We show that $\ker (\mathcal P^{k+1} - \lambda^{k+1}) = \ker (\mathcal P - \lambda)=\s\{g\}.$
\end{step}
\begin{stepproof}[Proof of Step~\ref{step:12}]
 It is clear that $\mathrm{ker}(\mathcal P-\lambda) \subset\mathrm{ker}(\mathcal P^{k+1}-\lambda^{k+1}).$ In the following, we show the reverse inclusion. Let $f\in \mathrm{ker}(\mathcal P^{k+1}-\lambda^{k+1})$. For every $j\in \{0,1,\ldots,k\}$, consider the functions
$$h_j \coloneqq \sum_{\ell=0}^{k} e^{\frac{-2 \pi i j \ell}{k+1}}\frac{1}{\lambda^\ell} \mathcal P^\ell f.$$
Note that 
\[\sum_{j = 0}^{k} h_j = \sum_{j = 0}^{k}\sum_{\ell=0}^{k} e^{\frac{-2 \pi i j \ell}{k+1}}\frac{1}{\lambda^\ell} \mathcal P^\ell f = \sum_{j = 0}^{k} f = {(k+1)}f.\]
Since $\calP^{k+1} f = \lambda^{k+1} f,$ we obtain that $\mathcal P h_j = \lambda e^{2\pi i j/(k+1)} h_j$
for every $j \in \{0,1, \dots, k\}$. From Lemma~\ref{lemma:point_periph_spec}, we have that $\lambda e^{2\pi i j/(k+1)} \not\in \sigma_{\mathrm{per}}(\mathcal P),$ therefore $h_j = 0$ for every $j\in \{1,\ldots, k\}.$ Therefore, $f = h_0 /(k+1) \in \ker{(\mathcal P - \lambda)}$.
\end{stepproof}

\begin{step}[3]\label{step:13} There exists a decomposition $\mathcal C^0(M) = \bigoplus_{j=0}^{k-1} \mathrm{ker}\left(\mathcal P^{k+1} - \lambda^{k+1} e^{2 \pi i j/k} \right) \oplus W_0, $ where $W_0$ is $\mathcal P^{k+1}$-invariant subspace of $\mathcal C^0(M)$ and $r(\calP^{k+1}\left.\right|_{W_0})<\lambda^{k+1}$. In particular, $\{\frac{1}{\lambda^n}\mathcal P^n\}_{n\in\mathbb N}$ is bounded.
\end{step}
\begin{stepproof}[Proof of Step~\ref{step:13}]
Recall from Lemma~\ref{lem:SF} {that}, $\mathcal P^{k+1}$ is a compact linear operator. Moreover, from Step~\ref{step:11} we obtain that $\sigma_{\mathrm{per}}(\frac{1}{\lambda}\mathcal P) = \sigma_{\mathrm{per}}(\frac{1}{\lambda^{k+1}}\mathcal P^{k+1})$. From~\cite[Theorems 8.4-3 and 8.4-5]{Kreyszig1978} and Lemma~\ref{lemma:point_periph_spec} we obtain that there exist non-zero $r_0,r_1,\ldots,r_{k-1}\in \N$ such that
\begin{equation}\label{eq:decspec}
    \mathcal C^0(M) =\bigoplus_{j = 0}^{k-1} \ker(\calP^{k+1} - \lambda^{k+1} e^{2\pi i j / k} \Id)^{r_j} \oplus W_0, 
\end{equation}
where $r_j = \inf\{m>0;\, \ker(\mathcal P^{k+1}- \lambda e^{2 \pi i j/k} )^{m+n} = \ker(\mathcal P^{k+1}- \lambda e^{2 \pi i j/k} )^{m}, \ \text{{for all} }n\in\N \},$ for each $j = 0,\dots, k-1$, and $W_0$ is $\mathcal P^{k+1}$-invariant satisfying $r(\calP^{k+1}\left.\right|_{W_0})<\lambda^{k+1}$. We show that $r_0 = r_1 = \ldots = r_{k-1}=1$. Using once again that $\mathcal P^{k+1}$ is a compact operator, we obtain from the Krein-Rutman theorem~\cite[Theorem 4.1]{Grobler1995} that the spectral radius $\lambda^{k+1} =r(\calP^{k+1})$ is a pole of maximal order in the spectral circle, i.e.~$r_0\geq \max\{r_1,\ldots,r_{k-1}\}.$ Suppose that $r_0>1,$ then there exists $f\in \mathcal C^0(M)$ such that $g= (\mathcal P^{k+1} - \lambda^{k+1})f.$ Integrating both sides with respect to $\mu$ yields $\int g \d\mu = 0$, so $r_0= r_1 = \ldots = r_{k-1} = 1$.
\end{stepproof}
This finishes the proof of Lemma~\ref{lem:powercompact}.
\end{proof}

\begin{lemma}\label{lem:alg_mult} 
There exists a decomposition  $\mathcal C^0(M) = \bigoplus_{j=0}^{k-1} \ker(\mathcal P -\lambda e^{\frac{2\pi i j}{k} }) \oplus W,$ where $W$ is a $\mathcal P$-invariant space $r\left(\left.\mathcal P\right|_{W}\right)<\lambda$, and $\dim \ker(\calP -\lambda e^{\frac{2\pi i j}{k}}) = 1$ for every $j\in\{0,1\ldots,k-1\}.$ 
\end{lemma}
\begin{proof}
From Lemma~\ref{lem:SF} and Lemma~\ref{lem:powercompact}, $\mathcal P^2$ is a compact linear operator and $\sup_{n\geq 0}\|\frac{1}{\lambda^n}\mathcal P^n\|<\infty$. Therefore, from~\cite[\emph{An extension of Frechet-Kryloff-Bogoliouboff's theorem}]{Yosida1} (see also~\cite[Théorèm above Définition 1.5]{Brunel} and~\cite[Equation (8) in the proof of Theorem 4]{Yosida2}), there exists a $\mathcal P$-invariant space  $W\subset \mathcal C^0(M)$ such that  $r\left(\left.\mathcal P\right|_{W}\right)<\lambda$ and   $$\mathcal C^0(M) = \bigoplus_{j=0}^{k-1} \ker\left(\mathcal P -\lambda e^{\frac{2\pi i j}{k} }\right) \oplus W.$$ 
Since $\frac{1}{\lambda}\calP$ is power-bounded, i.e.~$\sup_{n\geq 0}\left\|\frac{1}{\lambda^n}\calP^n\right\|<\infty$,~\cite[Theorem 5.1]{Gluck2016} implies that
\begin{align*}
    \mathrm{dim} \, \mathrm{ker}(\mathcal P-\lambda e^{2\pi i/k}) &\leq  \mathrm{dim} \,\mathrm{ker}(\mathcal P-\lambda e^{2\pi i 2/k}) \leq \ldots \\
    &\leq \mathrm{dim} \, \mathrm{ker}(\mathcal P-\lambda e^{2\pi i (k-1)/k}) \leq \mathrm{dim}\, \mathrm{ker}(\mathcal P-\lambda) =1,
\end{align*}
which concludes the proof.
\end{proof}

\subsection{Cyclic properties of \texorpdfstring{$\mathcal P$}{P}}

Consider ${\mathcal P}$ acting only on continuous functions $\calP :\mathcal C^0(M)\to \mathcal C^0(M).$ From Lemma~\ref{lem:alg_mult}, we know that
\( \mathcal C^0(M) = \ker (\calP^k -\lambda^k) \oplus W,\)
where $W$ is a $\calP^k$-invariant Banach space such that $r\left(\left.\calP^k\right|_{W}\right)< \lambda^k$ and $\dim\, \mathrm{ker}(\mathcal P^k - \lambda^k) = k.$

\begin{proposition}\label{prop:nerdola}
 There exist $k$ non-negative linearly independent eigenfunctions $g_0, \ldots, g_{k-1}\in \mathcal C^0_+(M) \cap \ker(\mathcal P ^k- \lambda^k)$ such that $\s_{\mathbb C} (\{g_i\}_{i=0}^{k-1}) = \ker ( \calP^k - \lambda^k  \Id)$ and  $\int g_i \d \mu =\mu(\{g>0\})$ for every $i\in\{0,1,\ldots,k-1\}.$ Moreover, these can be chosen such that the sets $C_i \coloneqq \{g_i>0\}$ are pairwise disjoint.
\end{proposition}
\begin{proof}
Observe that since $\lambda^k>0$ and $\calP(\mathcal C^0(M))\subset \mathcal C^0(M),$ it follows that if $f\in \mathcal C^0(M,\mathbb C)$ satisfies $\calP^k f = \lambda^k f,$ then $\calP^k \mathrm{Re}(f) = \lambda^k \mathrm{Re} (f)$ and $\calP^k \mathrm{Im}(f) = \lambda^k \mathrm{Im}(f)$.

Recall that $\mu$ is a measure on $M$ satisfying $\calP^*\mu = \lambda \mu$ and $\supp \mu = M$. Note that the operator $\calP^k$ satisfies
$\int_M \frac{1}{\lambda^k}\calP^k f(x)\mu(\d x) = \int_M f(x)\d \mu$, for every $f\in\mathcal C^0(M).$
By the same techniques of Theorem~\ref{thm:span_g_eps} (see also~\cite[Propositions 3.1.1 and 3.1.3]{Lasota1994}),  {it follows} that if $f\in\mathcal C^0(M)$ is an eigenfunction of $\calP^k$ associated with the eigenvalue $\lambda^k$, then $f^+(x) \coloneqq \max\{0,f(x)\}$ and $f^-(x) = \max\{0,-f(x)\}$ are also eigenfunctions of $\calP^k$ associated with the eigenvalue $\lambda^k.$ This provides a set of $k$ linearly independent non-negative continuous functions, $\{h_i\}_{i= 0}^{k-1}$ that span $\ker(\mathcal P ^k - \lambda^k)$. {We are} left to check that these can be chosen with pair-wise disjoint support. 

{Without loss of generality, define $G\coloneqq\{h_0>0\}\setminus \{h_1 >0\} \neq \emptyset$ and $H \coloneqq \{h_0> 0 \}  \cap \{h_1 > 0\}$.} {W}e claim that $\1_{G}h_1$ and $\1_{H}h_1$ are also a eigenfunctions of $\calP^k$ associated with the eigenvalue $\lambda^k.$ Observe that this is enough to conclude the proof since we can choose $k$ functions $g_0,\dots, g_{k-1}$ of the set below which have disjoint supports
$$\left\{h_i \1_{\{(\sum_{j=0}^{k-1} t_j h_j)^{\pm}>0\}};\ t_0,\ldots,t_{k-1}\geq 0 \ \text{and }i\in\{0,\ldots,k-1\}\right\}\subset \mathrm{ker}(\mathcal P^k-\lambda^k).$$

We {organise} the remainder of the proof into three steps:

\begin{step}[1]\label{step:app_no_mass_exchange1}
 We show that $\1_G \mathcal P ^k \1_H = 0$.
\end{step}
\begin{stepproof}[Proof of Step~\ref{step:app_no_mass_exchange1}]
Let $x \in G$ and assume that $\mathcal P ^k \1_H (x) > 0$. Then, $0<\mathcal P ^k \1_H (x) = \int \1_H(y) \mathcal P^k (x, dy),$ and since $h_1 > 0$ on $H$, $0< \int \1_H(y) h_1(y)(\mathcal P^*)^k\delta_x(dy),$
implying that $\mathcal P^k \1_H h_1 (x)> 0$. Moreover,
\[h_1(x) = \frac{1}{\lambda^k}\mathcal P ^k h_1 (x) = \frac{1}{\lambda^k}\mathcal P ^k (\1_H h_1 + (1 - \1_H) h_1)(x)  \geq \frac{1}{\lambda^k} \mathcal P ^k (\1_H h_1)(x) > 0,\]
which contradicts $x \not\in \{h_1>0\}.$ In particular, $\1_H \mathcal P^k \1_H = \mathcal P^k \1_H$.
\end{stepproof}
\begin{step}[2]\label{step:app_no_mass_exchange2}
We show that  $\1_H \mathcal P^k \1_G = 0$.
\end{step}
\begin{stepproof}[Proof of Step~\ref{step:app_no_mass_exchange2}]
{From} Step~\ref{step:app_no_mass_exchange1}, {it follows that}
\[\begin{split}
    \1_H h_0 &= \1_H \frac{1}{\lambda^k}\mathcal P^k h_0 = \1_H\frac{1}{\lambda} \mathcal P^k (\1_H h_0 + \1_G h_0)= \frac{1}{\lambda^k} \mathcal P^k (\1_H h_0) + \1_H \frac{1}{\lambda^k} \mathcal P^k (\1_G h_0).
\end{split}\]
Integrating {either} side and using $\mu \in \ker( (\mathcal P^k)^* - \lambda^k)$, {we obtain}
\begin{align*}\int \1_H h_0 \d\mu &= \int \frac{1}{\lambda^k} \mathcal P^k (\1_H h_0) \d \mu +  \int \1_H \frac{1}{\lambda^k} \mathcal P^k (\1_G h_0) \d\mu \\
&=\int \1_H h_0 \d \mu +  \int \1_H \frac{1}{\lambda^k} \mathcal P^k (\1_G h_0) \d\mu,
\end{align*}
implying that $\int \1_H \frac{1}{\lambda^k} \mathcal P^k (\1_G h_0) \d\mu = 0,$ and so $\1_H \mathcal P^k \1_G = 0$. Moreover, $\1_G \calP^k\1_G = \calP^k\1_G.$
\end{stepproof}
\begin{step}[3]\label{step:app_no_mass_exchange3}
 $\1_G h_0$ and $\1_H h_0$ are eigenfunctions of $\mathcal P^k$ with eigenvalue $\lambda^k$.
\end{step}
\begin{stepproof}[Proof of Step~\ref{step:app_no_mass_exchange3}]
From Steps~\ref{step:app_no_mass_exchange1} and~\ref{step:app_no_mass_exchange2}, it follows that 
\[\begin{split}
\1_H h_0 + \1_G h_0 &= h_0 = \frac{1}{\lambda^k} \mathcal P^k h_0 = \frac{1}{\lambda^k} \mathcal P^k(\1_H h_0 + \1_G h_0) \\ 
&= \frac{1}{\lambda^k} \mathcal P^k(\1_H h_0) + \frac{1}{\lambda^k} \mathcal P^k(\1_G h_0) = \1_H \frac{1}{\lambda^k} \mathcal P^k(\1_H h_0) + \1_G \frac{1}{\lambda^k} \mathcal P^k(\1_G h_0).
\end{split}\]
Since $G \cap H = \emptyset$ and $\mathcal P$ is strong Feller, the claim is verified.
\end{stepproof}
This finishes the proof of Proposition~\ref{prop:nerdola}.
\end{proof}
 
\begin{lemma}
Let $\{g_i\}_{i=0}^{k-1}\subset \mathcal C^0_+(M)$ be as in Proposition~\ref{prop:nerdola}. Then, these can be relabelled so that $\frac{1}{\lambda}\calP g_{i} = g_{{i-1\, \text{(mod $k$)}}}$, for $i \in \{0,1,\ldots, k-1\}$. In particular, we have that $g = \frac{1}{k} \sum_{i=0}^{k-1} g_i.$
\label{lem:ciclying}
\end{lemma}

\begin{proof}   

We divide the proof into a first step and conclude.

\begin{step}[1]\label{step:cont-func} There exists {a} continuous function $\theta:\{g>0\} \to \{0,1/k,2/k,\ldots, (k-1)/k\},$ such that
\begin{enumerate}[label = (\arabic*)]
     \item\label{it:lem-contfunc1} for every $j\in\{0,1,\ldots, k-1\}$, $\left.\theta\right|_{\{g_j >0\}} = \theta_j$ is constant,
     \item\label{it:lem-contfunc2} the set $\{g_j\}_{j = 0}^{k-1}$ can be relabelled so that $\theta_j = j/k$. 
\end{enumerate}
\end{step}
\begin{stepproof}[Proof of Step~\ref{step:cont-func}]
From Step~\ref{step:app_abs} of Lemma~\ref{lemma:point_periph_spec}, there exists {a} function $\theta:\{g>0\}\to \R$ such that $e^{2\pi i\theta(x)} g \in \ker(\mathcal P-\lambda e^{2\pi i /k}).$ Observe that by multiplying $\theta$ by a complex constant, we can assume without loss of generality that there exists $x\in M$ such that $\theta(x) =0$.  
Since $e^{2\pi i\theta(x)} g, g\in \ker (\mathcal P^k - \lambda^k) = \s\{g_0,\ldots, g_{k-1}\},$
there exist $\alpha_0,\ldots,\alpha_{k-1}\geq 0$ and $\theta_0,\ldots,\theta_{k-1} \geq 0$ such that
$$g = \sum_{j=0}^{k-1} \alpha_j g_j\ \text{and }e^{2\pi i \theta }g = \sum_{j=0}^{k-1} \alpha_je^{2\pi i \theta_j} g_j.$$
Since $\{g_{j_1}>0\} \cap \{g_{j_2}>0\}=\emptyset$ if $j_1\neq j_2$, then $\theta (x) = \theta_j$ for every $x\in \{g_j>0\}$. This {proves} (1). 

Without loss of generality, we may assume that $\alpha_0 \neq 0$ and $\theta_0 = 0$. Let us fix $x \in \{g_0 > 0\}$. Then, 
\[ e^{2\pi i/k} g(x) = \frac{1}{\lambda}\calP (e^{2\pi i \theta } g)(x) = \int_{M}  e^{2\pi i \theta(y)} g(y) \frac{1}{\lambda}\mathcal P(x,\d y), \]
and therefore
$$  g(x) = \int_{M}  e^{2\pi i (\theta(y)- 1/k)} g(y) \frac{1}{\lambda}\mathcal P^*\delta_x(\d y). $$
Since $g(x) = \int g(y) \frac{1}{\lambda}\mathcal P(x,\d y),$ and $\theta$ is continuous we obtain
$$\theta(y) = \frac{1}{k}\ \text{for every }y\in \supp \calP(x,\d y)\cap \{g>0\}.$$
The same argument for $\mathcal P^n$ yields
\begin{align}
   \theta(y) = \frac{n}{k}\ \text{for every }y\in \supp \calP^n(x,\d y)\cap \{g>0\}.\label{eq:theta} 
\end{align}
Note that if $y \in \supp \calP^n(x,\d y)\cap \{g_j>0\}$ for some $j$, then $\theta_j = \theta(y) = n/k$. This implies that $\supp \calP^m(x,\d y) \cap  \{g_j > 0\} = \emptyset$ for any $m \neq {n\, \text{(mod $k$)}}$. Since $\supp \mathcal P^k(x,\d y) \cap \{g_0\}  \neq \emptyset$ and there are exactly $k$ functions $g_0, \ldots, g_k$, each must have a different phase $\theta_j$. After relabelling, we may assume that $\theta_j = j/k$, for $j \in \{0,1,\ldots, k\}$, showing $(2)$.
\end{stepproof}


We may now conclude the proof of the lemma. It follows immediately from equation~\eqref{eq:theta} that $\{\mathcal P g_j >0\} \subset \{g_{{j-1\, \text{(mod $k$)}}} > 0\}$. 
Moreover, since $e^{2\pi i \theta} g \in \ker(\calP - \lambda e^{2\pi i /k})$, we have
\[
\lambda e^{2\pi i/k} e^{2\pi i \theta} g = \lambda e^{2\pi i/k} \sum_{j = 0}^{k-1} \alpha_j  e^{2\pi ij/k} g_j = \mathcal P (e^{2\pi i \theta} g) = \sum_{j = 0}^{k-1} \alpha_j  e^{2\pi ij/k} \mathcal P g_j,
\]
so $\lambda \alpha_{j-1}g_{j-1}= \alpha_j \mathcal P g_j,$
with $-1 = k-1$. Integrating both sides with respect to $\mu$ yields $\alpha_{{j-1\, \text{(mod $k$)}}} = \alpha_{j}$, from which we conclude that $\alpha_0 = \alpha_1 = \ldots = \alpha_{k-1}$ and $\frac{1}{\lambda} \mathcal P g_j = g_{{j-1\, \text{(mod $k$)}}}$, for every $j \in \{0,1, \ldots, k-1\}$.
\end{proof}

The following corollary follows directly from Lemma~\ref{lem:ciclying}.

\begin{corollary}
\label{cor:dec}
Every function 
\(f_{\ell} \coloneqq \frac{1}{k}\sum_{j=0}^{k-1} e^{2\pi i j \ell/k  } g_j \)   
satisfies $\calP f_{\ell} =\lambda e^{2\pi i \ell/k } f_{\ell },$ i.e.~$\ker(\calP - \lambda e^{2\pi i \ell /k}) = \s( f_\ell)$, for $\ell \in \{0,1,\ldots, k-1\}$.
\end{corollary}

\subsection{Existence of quasi-ergodic measures} \label{app:g_positive}

Recall that $g \in \mathcal C^0_+(M)$ is the unique function satisfying $\calP g = \lambda g$, up to a multiplicative factor.

\begin{lemma}\label{lem:g>0}
$\calP \1_{\{g>0\}}\leq c \1_{\{g>0\}}$, for some constant $c>0$. 
\end{lemma}
\begin{proof}

Observe that for every $a>0$ we have $\mathcal P\1_{\{g>a\}} \leq \frac{1}{a}\mathcal Pg = \frac{\lambda}{a} g.$ Hence, $\{\mathcal P\1_{\{g>a\}} >0\} \subset \{g>0\}.$ Since $\1_{\{g>0\}} = \sum_{n=1}^\infty \1_{\{\|g\|_{\infty}/n \geq g > \|g\|_{\infty}/(n+1)\}}$, we obtain that
\begin{align*}
 \{\mathcal P \1_{\{g>0\}} >0\} &\subset  \bigcup_{n\in\mathbb N}\{\mathcal P\1_{\{g>1/n\}}>0\} \subset \{g>0\}.
\end{align*}
It follows that $\mathcal P \1_{\{g>0\}}\leq \|\mathcal P\|\1_{\{g>0\}}.$
\end{proof}

\begin{notation}
    We define the operator $\calP_g :L^\infty(\{g>0\},\mu) \to L^\infty(\{g>0\},\mu)$ as the operator ${\calP_g f \coloneqq \mathcal P(\1_{\{g>0\}} f)}$, \textcolor{black}{where we may extend the definition of $f$ to $\{g > 0\}$ by setting it to be zero where undefined.}
\end{notation}

\begin{corollary}
 The measure $\widetilde{\mu}(\d x)\coloneqq \mu(\dx \cap \{g>0\})/\mu(\{g>0\})$ satisfies $\calP_g^* \widetilde \mu = \lambda \widetilde \mu$.
\end{corollary}
\begin{proof}
From Lemma~\ref{lem:g>0} we have that for every $h\in L^\infty (\{g>0\}),$ $\calP_g h =  \mathcal P ( \1_{\{g>0\}} h ) = \1_{\{g>0\}} \mathcal P_g h$. Therefore,
\begin{align*}
 \int h \d\, (\calP_g^*\widetilde{\mu}) &= \int \calP_g h \d\widetilde \mu = \frac{1}{\mu(\{g>0\})} \int \1_{\{g >0\}} \calP_g h \d \mu = \frac{1}{\mu(\{g>0\})} \int \calP(\1_{\{g >0\}} h) \d \mu \\ 
 &= \frac{1}{\mu(\{g>0\})} \int\1_{\{g >0\}} h \d (\mathcal P^*\mu) = \frac{\lambda}{\mu(\{g>0\})} \int\1_{\{g >0\}} h \d \mu  = \lambda \int h \d \widetilde \mu.
\end{align*}
\end{proof}

Observe that since $\int g\, \d\mu = \mu(\{g>0\})$, we have that $\int g \,\d\widetilde \mu=1$. The above corollary implies that $\sigma_{\mathrm{per}}(\calP) = \sigma_{\mathrm{per}}(\calP_g)$ and 
\[\1_{\{g>0\}}\ker (\calP - \lambda e^{2\pi i j/k}) = \ker (\calP_g - \lambda e^{2\pi i j/k}),\ \text{for every }j\in\{0,\ldots,k-1\}.\] 

Since each $g_i$ defined in Lemma~\ref{lem:alg_mult} satisfies $C_i = \{g_i>0\}\subset \{g>0\},$ we can assume by abuse of notation that $g_i\in L^\infty(\{g>0\},\widetilde \mu)$. Moreover,
\begin{align}\label{eq:newdec}
    L^{\infty}(\{g>0\},\widetilde \mu) = \s(g_0,\ldots,g_{k-1}) \oplus V, 
\end{align}
where $V$ is $\calP$-invariant and $r(\left.\calP\right|_{V})< \lambda$.

\begin{lemma}
For every $i\in\{0,1,\ldots,k-1\}$ define $\widetilde \mu_i (\d x)= \widetilde \mu( C_i \cap \d x)$, where $C_i = \{g_i>0\}$. Then $\mathcal P_g^* \widetilde \mu_i = \lambda \widetilde \mu_{i+1\, (\mathrm{mod}\ k)}$. \label{lem:mcycle}
\end{lemma}
\begin{proof}
We claim that $v \in V$ if and only if $\int_{C_i} v \, \d \widetilde \mu = 0$ for every $i\in\{0,1,\ldots,k-1\}.$ If the claim holds, to conclude the proof of the lemma, take $f\in L^\infty(\{g>0\},\mu).$ Therefore, $f = \sum_{i=0}^{k-1} \alpha_i g_i + v$ with $v\in V$. From the proof of Lemma~\ref{lem:mcycle}, it follows that
$$\alpha_i = \int_{C_i} f \, \d\widetilde\mu = \int f\, \d \widetilde \mu_i,$$
and
\begin{align*}
    \int f\, \d\calP^* \widetilde \mu_i &= \int \calP f\, \d \widetilde \mu_i   = \sum_{j=0}^{k-1}\int \mathcal P(\alpha_j g_j)\, \d\widetilde \mu_i \\
    &=\sum_{j=0}^{k-1}\int \lambda \alpha_j g_{j-1}\, \d\widetilde \mu_i  =\lambda \alpha_{i+1} =  \lambda \int f\,\d\widetilde\mu_{i+1\ (\mathrm{mod}\ k)}. 
\end{align*}

Let us now show the claim. Suppose first that $v \in V$. We claim that $\1_{C_i} v  \in V$ for all $i  \in \{0,1,\ldots, k-1\}$. Indeed, if $\1_{C_i} v \not \in V$, then $ v = \alpha_i g_i + w + \sum_{j\neq i}\1_{C_j}v$ with $\alpha_i \neq 0$ and $w \in V$. Since, $C_i \cap C_j = \emptyset$ for all $j \neq i$, we get that $v \not \in V$. It follows that $$\left|\int_{C_i} v\,\d\widetilde\mu\right| =\left|\int_{M} \1_{C_i}v\,\d\widetilde\mu\right| =\left|\int_{M}\frac{1}{\lambda^n}\mathcal P^n( \1_{C_i}v )\,\d \widetilde\mu \right|\leq \left\|\left.\frac{1}{\lambda^n}{\mathcal P^n}\right|_{V}\right\|\|v\|\xrightarrow[]{n\to\infty}0.$$

On the other hand, assume that $\int_{C_i} v \,\d \widetilde\mu = 0$ for every $i\in\{0,1,\ldots,k-1\}.$ Write $v = \sum_{i=0}^{k-1}\alpha_i g_i + w$, with $w\in V$. Observing that $\int g_i \d\widetilde \mu = 1$, we have
$$\alpha_i = \int_{C_i} \alpha_i g_i \,\d\widetilde\mu = \int_{C_i} \left(\sum_{j=0}^{k-1}\alpha_j g_j + w\right) \,\d\widetilde\mu = \int_{C_i} v\, \d\widetilde \mu = 0.$$
We obtain that $\alpha_i = 0$ for every $i\in\{0,1,\ldots,k-1\},$ which implies $v\in V$.
\end{proof}

\begin{theorem}\label{thm:QED}
Assume that $\mathcal P$ satisfies Hypothesis~\ref{hyp:A}. Given a bounded and measurable function $h:\{g>0\}\to \mathbb R$ we have that for every $x\in \{g>0\}$
\[\frac{1}{\mathbb E_x[e^{S_n \phi}\1_{\{\tau>n\}}]}\mathbb E_x\left[e^{S_n\phi} \1_{\{\tau> n\}} \frac{1}{n}  \sum_{i = 0}^{n-1} h\circ X_i\right]\xrightarrow[]{n\to\infty} \frac{\int hg\, \d\mu}{\int g\,\d\mu},\]
where $\tau =\min \{n\in \mathbb N; X_n\not \in \{g>0\}\}$ and $S_n \phi = \sum_{i=0}^{n-1}\phi\circ X_i.$ In other words, there exists a unique quasi-ergodic measure for the $e^\phi$-weighted Markov process $X_n^\phi$ on $\{g>0\}$ \textcolor{black}{and it satisfies $\d \nu = g\d\mu/\int g \d \mu$.}
\end{theorem}

\begin{proof}
In this proof, we adopt the notation $g_{m} := g_{m\, \text{(mod $k$)}}$, $\widetilde\mu_{m} := \widetilde\mu_{m\, \text{(mod $k$)}}$ and $C_m= C_{m\, \text{(mod $k$)}}.$
Recall that $$\frac{\int h(x)g(x) \mu(\d x)}{\int g(x)\mu(\d x)} = \int h(x)g(x)\,\widetilde\mu(\d x). $$
Given $n\in\N$ and $x\in \{g>0\}$ define
\begin{align*}
    Q_{h}^n(x):= \left| \frac{1}{\mathbb E_x[e^{S_n \phi}\1_{\{\tau>n\}}]}\mathbb E_x\left[e^{S_n\phi} \1_{\{\tau> n\}} \frac{1}{n} \sum_{i = 0}^{n-1} h\circ X_i\right]- \int h(x)g(x) \widetilde\mu(\d x)\right|.
\end{align*}
Observe that to prove the theorem, it {suffices} to show that for every measurable and bounded non-negative $h:\{g>0\}\to \R$ we have
$$ \max\left\{Q_h^{nk + \ell}(x);\, \ell\in\{0,1,\ldots, k-1\}\right\}\xrightarrow[]{n\to\infty}0,$$
for every $x\in C_s$ where $s\in\{0,1,\ldots,k-1\}$. Moreover, using the Markov property for the process $X_n$ we may write
\[\E_x\left[e^{S_n\phi}\1_{\tau_n}h\circ X_i\right] = \calP^i_g\left(h\calP^{n-i}_g \1_{\{g>0\}}\right)(x).\]
We divide the {remainder} of the proof into three steps. 

\begin{step}[1]\label{step:app_cyclicity}
{F}or every bounded and measurable function $h:\{g>0\}\to\R$, $\ell,s\in \{0,1,\ldots,k-1\}$ and $x\in C_s$ we have
$$\lim_{n\to\infty}\frac{1}{\lambda^{nk+\ell}}\mathcal P_g^{n k +\ell} h(x)= g_{s}(x) \int h \, \d\widetilde\mu_{s+\ell}.$$
\end{step}
\begin{stepproof}[Proof of Step~\ref{step:app_cyclicity}]
From Step~\ref{step:app_abs} of Lemma~\ref{lem:mcycle} it is clear that 
$$h = \sum_{j=0}^{k-1} g_j \int_{C_j} h\, \d\widetilde \mu + v, $$
with $v\in V$. Since $\mathcal P_g^{nk +\ell} g_{j}(x) = \lambda^{nk +\ell} g_{j-\ell}(x)$, we obtain that
\begin{equation}\label{eq:limit}  
\frac{1}{\lambda^{nk+\ell}}\mathcal P_g^{nk+\ell}h = \sum_{j=0}^{k-1} g_{j-\ell} \int_{C_j} h\, \d\widetilde \mu + \frac{1}{\lambda^{nk+\ell}}\calP^{nk+\ell} v  \xrightarrow[]{n\to\infty} \sum_{j=0}^{k-1} g_{j-\ell}\int h \, \d\widetilde \mu_j.  
\end{equation}
Finally, if $x\in C_s,$ then
$$\lim_{n\to\infty} \frac{1}{\lambda^{nk+\ell}}\mathcal P_g^{nk+\ell}h(x) = g_s(x) \int h\, \d \widetilde\mu_{s+\ell},$$
which yields the claim.
\end{stepproof}


\begin{step}[2]\label{step:app_cyclicity_ergodic}
{F}or every non-negative bounded and measurable function $h:\{g>0\}\to\R$, $\ell,s\in \{0,1,\ldots,k-1\}$ and $x\in C_s$ we have
$$\lim_{n\to\infty}\frac{1}{nk+\ell} \sum_{i=0}^{nk+\ell-1}\frac{1}{\lambda^i}\mathcal P_g^{i} \left(h\frac{1}{\lambda^{nk+\ell - i }}\calP_g^{nk+\ell - i }\1_{\{g>0\}}\right)(x) = g_{s}(x)\widetilde\mu(C_{s+\ell}) \int h g \, \d \widetilde \mu.$$
\end{step}
\begin{stepproof}[Proof of Step~\ref{step:app_cyclicity_ergodic}]

We denote $\mathcal G:= \mathcal P_g/\lambda$ to simplify the notation and improve readability. Recall that $\1_{\{g>0\}} = \sum_{j=0}^{k-1} \widetilde \mu(C_j) g_j + v$, where $v\in V$. It follows that,
\begin{align}\label{eq:gh}
   \sum_{i=0}^{nk+\ell-1}\mathcal G^{i}&\left(h {\mathcal G}^{nk+\ell - i }\1_{\{g>0\}}\right)(x) =\\*
   &=\sum_{j=0}^{k-1}\widetilde \mu(C_j) \sum_{i=0}^{nk+\ell-1} \mathcal G^{i}\left( h \mathcal G^{nk+\ell - i} g_j\right)(x) + \sum_{i=0}^{nk+\ell-1} \mathcal G^{i} (h  \mathcal G^{nk+\ell - i} v)(x)\nonumber\\
   &=\sum_{j=0}^{k-1}\widetilde \mu(C_j) \sum_{i=0}^{nk+\ell-1} \mathcal G^{i}\left( h g_{j-\ell+i}\right)(x) + \sum_{i=0}^{nk+\ell-1} \mathcal G^{i} (h  \mathcal G^{nk+\ell - i} v)(x).
\end{align}
Observe that
\begin{align}\label{eq:gv0}
    \left|\frac{1}{nk+\ell}\sum_{i=0}^{nk+\ell-1} \mathcal G^{i} (h  \mathcal G^{nk+\ell - i} v)(x)\right| &\leq  \sup_{i\geq 0} \|\mathcal G^i\| \frac{1}{nk +\ell} \sum_{i=0}^{nk+\ell-1}\|h  \mathcal G^{nk+\ell - i} v\|_{\infty} \nonumber\\
    &\leq  \sup_{i\geq 0} \|\mathcal G^i\| \frac{1}{nk +\ell} \sum_{i=0}^{nk+\ell-1}\|h\|_{\infty} \left\|  \mathcal G^{nk+\ell - i} v\right\|_{\infty} \xrightarrow[]{n\to\infty} 0.
\end{align}
For the first term, observe that $hg_{j-\ell + i} = g_{j-\ell +i} \int h g_{j-\ell + i} \d\widetilde\mu + v$, for some $v \in V$. Therefore, since $x \in C_s$, 
\begin{align*}
\frac{1}{nk+\ell}&\sum_{j=0}^{k-1}\widetilde \mu(C_j) \sum_{i=0}^{nk+\ell-1} \mathcal G^{i}\left( h g_{j-\ell+i}\right)(x) = \\ &= \frac{1}{nk+\ell}\sum_{j=0}^{k-1}\widetilde \mu(C_j)\sum_{i=0}^{nk+\ell-1} \left(\int h g_{j-\ell + i} \d\widetilde\mu\right)g_{j - \ell}(x) + \mathcal G^i (v)(x)\\
&=\widetilde{\mu}(C_{s + \ell}) g_s(x) \frac{1}{nk + \ell}  \sum_{i = 0}^{nk + \ell -1}\int hg_{s + i}\d\widetilde{\mu} + \mathcal G^i (v)(x)\\&= \widetilde{\mu}(C_{s + \ell}) g_s(x) \frac{1}{nk + \ell}  \left(\sum_{i = 0}^{n-1}\sum_{r = 0}^{k-1}\int hg_{s + r}\d\widetilde{\mu} + \sum_{t = 0}^{\ell}\int hg_{s + t}\d\widetilde{\mu} + \sum_{j = 0}^{nk +\ell -1}\mathcal G^j (v)(x)  \right)\\
& = \widetilde{\mu}(C_{s + \ell}) g_s(x) \frac{1}{nk + \ell}  \left(nk\int hg\d\widetilde{\mu} + k\int hg\d\widetilde{\mu} + \sum_{j = 0}^{nk +\ell -1}\mathcal G^j (v)(x)  \right)\\
&\xrightarrow[]{n \to \infty} \widetilde\mu(C_{s + \ell})g_s(x)\int hg \, \d\widetilde\mu.
\end{align*}

\end{stepproof}

We may now conclude the proof of the theorem. From Steps~\ref{step:app_cyclicity} and~\ref{step:app_cyclicity_ergodic}, we obtain that
$$ \lim_{n\to\infty}\frac{\lambda^{nk+\ell}}{\calP_g^{nk+\ell}\1_{\{g>0\}}(x)} = \frac{1}{g_s(x) \widetilde \mu(C_{s+\ell})},$$
and
\begin{align*}
    \lim_{n\to\infty}\frac{1}{{nk+\ell}}\sum_{i=0}^{{nk+\ell}-1} \frac{1}{\lambda^i}\calP_g^i\left(h \frac{1}{\lambda^{{nk+\ell}-i}} \calP_g^{nk+\ell-i}\1_{\{g>0\}}\right)(x) 
    &= g_s(x) \widetilde \mu(C_{s+\ell}) \int h g\, \widetilde \mu.
\end{align*}

Therefore, $Q_{h}^{nk+\ell}(x) \xrightarrow[]{n\to\infty} 0$ for {all} $s,\ell\in \{0,1,\ldots,k-1\}$ and $x\in C_s,$ which concludes the proof of the theorem.
\end{proof}



\begin{theorem}\label{thm:QED2}
Assume that $\mathcal P$ satisfies Hypothesis~\ref{hyp:A} and $\sigma(\frac{1}{\lambda}\mathcal P)\cap \mathbb S^1 = \{1\}$. Then, given a bounded measurable function $h:M\to \mathbb R$, for every $x\in \{g>0\}$,
\[\frac{1}{\mathbb E_x[e^{S_n \phi}\1_{\{\tau>n\}}]}\mathbb E_x\left[e^{S_n\phi} \1_{\{\tau> n\}} \frac{1}{n} \sum_{i = 0}^{n-1}h\circ X_i\right]\xrightarrow[]{n\to\infty} \frac{\int hg\, \d\mu}{\int g\, \d\mu},\]
where $\tau {\coloneqq}\min\{n; X_n\not\in M\}$ and $S_n\phi = \sum_{i=0}^{n-1}\phi\circ X_i$. In other words, there exists a unique quasi-ergodic of the $e^\phi$-weighted Markov process $X_n^\phi$ on $M$.
\end{theorem}

\begin{proof}
Note that the spectral gap in the operator $\frac{1}{\lambda}\mathcal{P}$, along with its strong Feller property, ensures that for any bounded and measurable function $h: M \to \mathbb{R},$ it holds that
\begin{align}
    \sup_{x\in M}\left|\frac{1}{\lambda^n} \mathcal P^n h(x) - g(x) \int h\d \widetilde\mu\right| \xrightarrow{n\to\infty} 0.\label{eq:T1}
\end{align}

Repeating the proof of Step~\ref{step:app_cyclicity_ergodic} of Theorem~\ref{thm:QED} we obtain that
\begin{align}
     \sup_{x\in M}\left|\frac{1}{n} \sum_{i=1}^{n-1} \frac{1}{\lambda^i} \calP^i\left(h \frac{1}{\lambda^{n-i}}\calP^{n-i}\1_{M}\right)(x) - g(x) \int hg\, \d \widetilde\mu \right|\xrightarrow{n\to\infty} 0. \label{eq:T2}
\end{align}

Combining equations~\eqref{eq:T1}-\eqref{eq:T2} and the same computations in the proof of Theorem~\ref{thm:QED} we obtain the result.
\end{proof}